\documentclass[10.5pt]{article}
\usepackage{amsmath}
\usepackage{amsfonts}
\usepackage{graphicx}
\usepackage{amsmath}
\usepackage{amssymb}
\usepackage{latexsym}
\usepackage{amsmath,amsfonts,amssymb,amsthm,euscript,makeidx,color,mathrsfs}
\usepackage{hyperref}
\usepackage{cite} 
\usepackage{enumerate}

\hypersetup{hypertex=true,
colorlinks=true,
linkcolor=blue,
anchorcolor=blue,
citecolor=blue}

\def\sqr#1#2{{\vcenter{\vbox{\hrule height.#2pt
              \hbox{\vrule width.#2pt height#1pt \kern#1pt \vrule width.#2pt}
              \hrule height.#2pt}}}}
\def\signed #1{{\unskip\nobreak\hfil\penalty50
              \hskip2em\hbox{}\nobreak\hfil#1
              \parfillskip=0pt \finalhyphendemerits=0 \par}}
\def\endpf{\signed {$\sqr69$}}
\def\3n{\negthinspace \negthinspace \negthinspace }
\def\2n{\negthinspace \negthinspace }
\def\1n{\negthinspace }

\def\dbA{\mathbb{A}}

\def\dbC{\mathbb{C}}

\def\dbE{\mathbb{E}}

\def\dbG{\mathbb{G}}
\def\dbH{\mathbb{H}}

\def\dbN{\mathbb{N}}

\def\dbP{\mathbb{P}}
\def\dbQ{\mathbb{Q}}
\def\dbR{\mathbb{R}}
\def\dbS{\mathbb{S}}

\def\dbX{\mathbb{X}}
\def\dbY{\mathbb{Y}}
\def\dbZ{\mathbb{Z}}

\def\sL{\mathscr{L}}
\def\sM{\mathscr{M}}
\def\sN{\mathscr{N}}

\def\sP{\mathscr{P}}

\def\sR{\mathscr{R}}
\def\sS{\mathscr{S}}

\def\sX{\mathscr{X}}

\def\={\buildrel \triangle \over =}

\def\ds{\displaystyle}

\def\ns{\noalign{\ss}}

\def\e{\varepsilon}

\def\si{\sigma}
\def\t{\tau}
\def\f{\varphi}
\def\th{\theta}

\def\D{\Delta}
\def\Th{\Theta}

\def\cA{{\cal A}}
\def\cB{{\cal B}}
\def\cC{{\cal C}}
\def\cD{{\cal D}}

\def\cI{{\cal I}}
\def\cJ{{\cal J}}

\def\cL{{\cal L}}

\def\cP{{\cal P}}
\def\cQ{{\cal Q}}

\def\cT{{\cal T}}

\def\cX{{\cal X}}
\def\cY{{\cal Y}}

\def\BB{{\bf B}}

\def\BD{{\bf D}}

\def\BR{{\bf R}}

\def\no{\noindent}

\def\ss{\smallskip}
\def\ms{\medskip}

\def\q{\quad}
\def\qq{\qquad}

\def\esssup{\mathop{\rm esssup}}

\def\h{\widehat}
\def\wt{\widetilde}
\def\cd{\cdot}

\def\ae{\hbox{\rm a.e.{ }}}

\def\({\Big (}
\def\){\Big )}
\def\[{\Big[}
\def\]{\Big]}

\def\({\Big (}
\def\){\Big )}
\def\[{\Big[}
\def\]{\Big]}

\def\bde{\begin{definition}}
\def\ede{\end{definition}}
\def\be{\begin{equation}}
\def\bel{\begin{equation}\label}
\def\ee{\end{equation}}
\def\bt{\begin{theorem}}
\def\et{\end{theorem}}
\def\bc{\begin{corollary}}
\def\ec{\end{corollary}}
\def\bl{\begin{lemma}}
\def\el{\end{lemma}}
\def\bp{\begin{proposition}}
\def\ep{\end{proposition}}
\def\bas{\begin{assumption}}
\def\eas{\end{assumption}}
\def\br{\begin{remark}}
\def\er{\end{remark}}
\def\ba{\begin{array}}
\def\ea{\end{array}}
\def\bpf{\begin{proof}}
\def\epf{\end{proof}}
\def\ed{\end{document}}

\def\square#1{\vbox{\hrule\hbox{\vrule height#1%
     \kern#1\vrule}\hrule}}
\def\rectangle#1#2{\vbox{\hrule\hbox{\vrule height#1%
     \kern#2\vrule}\hrule}}

\font\tenbb=msbm10 \font\sevenbb=msbm7 \font\fivebb=msbm5

\newfam\bbfam
\scriptscriptfont\bbfam=\fivebb \textfont\bbfam=\tenbb
\scriptfont\bbfam=\sevenbb

\newtheorem{lemma}{Lemma}[section]
\newtheorem{remark}{Remark}[section]

\newtheorem{theorem}{Theorem}[section]
\newtheorem{corollary}{Corollary}[section]

\newtheorem{definition}{Definition}[section]
\newtheorem{proposition}{Proposition}[section]
\newtheorem{assumption}{Assumption}[section]

\makeatletter
   
   \@addtoreset{equation}{section}
\makeatother

\def\ges{\geqslant}
\def\les{\leqslant}

\begin{document}

\title{Causal feedback strategies for controlled stochastic Volterra systems: a unified treatment\thanks{
{\bf Funding}: This work was supported by the National Natural Science Foundation of China
(11971332, 11931011, 12371449), the Science Development Project of Sichuan University (2020SCUNL201).}}

\author{Jiayin Gong\thanks{ School of Mathematics, Sichuan University, Chengdu 610065, China. ({gongjiayin916@163.com}).}
\and Tianxiao Wang\thanks{ Corresponding author. School of Mathematics, Sichuan University, Chengdu 610065, China. ({wtxiao2014@scu.edu.cn}).}}

\maketitle

\begin{abstract}
This paper is concerned with a unified treatment of linear quadratic control problem for stochastic Volterra integral equations (SVIEs), motivated by the various approaches and scattered results in the existing papers.
A novel class of optimal causal feedback strategy is introduced and characterized by means of a new Riccati system whose solvability is carefully discussed. To this end, a fundamental function space and an appropriate multiplicative rule among functions are defined for the first time.
In contrast with the existing works, our unified treatment not only provides a new approach, but also extends or improves the
 known conclusions in stochastic differential equations, convolution SVIEs, stochastic Volterra integro-differential equations (VIDEs), deterministic VIEs,
deterministic VIDEs. In addition, we found that for stochastic VIEs, the memory effect in terms of state process in general can not be demonstrated in the feedback controls, even though it holds true for certain deterministic VIEs in the literature.

\end{abstract}



\section{Introduction}

Let $(\Omega,\mathcal{F},\mathbb{F},\dbP)$ be a complete filtered probability space on which a standard one-dimension Brownian motion $W$
is defined and $\mathbb{F}=\{\mathcal{F}_t\}_{t\geq0}$ is the natural filtration of $W$ augmented by all the $\dbP$-null sets in $\mathcal{F}$.
Given $\t\in[0,T]$, and $\f\1n\in\1nC([\t, T];\dbR^d)\1n\=\1n \{\f(\cd): [\t, T]\rightarrow\dbR^d\big|\f$ is continuous $\1n\}$,
we consider the following stochastic Volterra integral equations (SVIEs for short):
\vskip-6mm
\bel{E1}\ba{ll}
\ns\ds X(t)=\varphi(t)+\int_\t^t\big[A(t,s)X(s)+B(t,s)u(s)+b(t,s)\big]ds\\[-1mm]
\ns\ds \qq\q+\int_\t^t\big[C(t,s)X(s)+D(t,s)u(s)+\sigma(t,s)\big]dW(s),\ \ t\in[\t,T],
\ea\ee
\vskip-2mm
\no where $A,B,C,D$ are given deterministic matrix-valued functions, $b$ and $\si$ are vector-valued stochastic inhomogeneous terms.
$\varphi$ is called  \emph{free term}, $u$ is the \emph{admissible control process} in
\vskip-7mm
$$\ba{ll}
\ns\ds \mathcal{U}[\t,T]\1n\=\1n\Big\{\1n u\1n:\1n[\t,T]\1n\times\1n\Omega\1n\rightarrow \1n\mathbb{R}^l\big|u(\cdot)\mbox{ is }
 \mathbb{F}\mbox{-adapted, measureable},\   \1n
\dbE\2n\int^T_\t\2n|u(s)|^2ds\1n<\1n\infty\1n\Big\}.
\ea
$$
\vskip-3mm
\no We introduce the following cost functional:
\vskip-5mm
\bel{E1.2}\ba{ll}
 \ns\ds \cJ(\t,\varphi;u)\1n=\1n\dbE\bigg[\langle GX(T),X(T)\rangle\1n+
\3n\int_\t^T\3n\Big(\langle Q(s)X(s),X(s)\rangle\1n+\1n\langle R(s)u(s),u(s)\rangle \\
 \ns\ds \qq\qq+2\langle S(s)X(s),u(s)\rangle\1n+2\langle q(s) ,X(s)\rangle\1n+2\langle \rho(s) ,u(s)\rangle \Big) ds\bigg],
\ea\ee
\vskip-2mm
\no where $G, Q, S$ and $R$ are given deterministic matrix-valued functions, and $q$ and $\rho$ are vector-valued adapted processes.
 The LQ problem for SVIE can be stated as follows:

\medskip

\textbf{Problem (LQ-SVIE)}. Given $\t\1n\in\1n[0,T],\f\1n\in\1nC([\t, T];\dbR^d),$
 find a $\hat{u }\in \mathcal{U}[\t,T] $
s.t.
\vskip-5mm
\bel{E2}
\cJ(\t,\varphi;\hat{u } )=\inf_{u\in \mathcal{U}[\t,T]}\cJ(\t,\varphi;u )\equiv V(\t,\f).
\ee
\vskip-2mm

Any $\hat{u } \in\mathcal{U}[\t,T]$ satisfying (\ref{E2}) is called an \emph{optimal  open-loop control} of Problem (LQ-SVIE)
corresponding to $(\t,\varphi)$, $\widehat{X}(\cdot)\equiv X(\cdot;\t,\varphi,\hat{u } (\cdot))$ is called an
\emph{optimal open-loop state process}, $V$ is called the \emph{value function}.

\ss

The classical optimal control theory was originally developed to deal with ordinary differential equations (ODEs).
However, it is found that there are many physical phenomena (e.g., the hereditary property in Bellman--Cooke \cite{Bellman-Cooke-1963}, Bellman--Danskin \cite{Bellman-Danskin-1954}, Volterra \cite{Volterra-book-1930}) that cannot be adequately described by ODEs but by Volterra integral equations (VIEs for short).
VIEs were introduced by Italian mathematician Vito Volterra,
and have brought us new theories (e.g., Volterra operator theory in Gohberg--Krein \cite{Gohberg-Krein-book-2004}).
Up until now, it is still an active research topic, see, e.g., the recent monograph of Brunner \cite{Brunner-2017}.
As to the optimal control theory of VIEs, it seems that one of the earliest paper was given Friedman \cite{Friedman-1964} to our best knowledge. In the stochastic case, Yong \cite{Yong-2008} studied the maximum principle of optimal controls by introducing proper
backward stochastic Volterra integral equations (BSVIEs for short), while Wang-Yong \cite{Wang-Yong-2023-SICON} recently extended it into the general case with arbitrary non-empty control domain. Some relevant papers include
Hamaguchi \cite{Hamaguchi-2021, Hamaguchi-2023}, Shi et al \cite{Shi-Wang-Yong-2015}, Viens-Zhang \cite{Viens-Zhang-2019}, Wang \cite{Wang-2020}, Wang-Zhang \cite{Wang-Zhang-2017}, to mention a few.

\ss

At this moment, let us return back to the linear quadratic problem. Practically, people expect the
optimal control to have state feedback representation which is non-anticipating. In the case of ODEs or SDEs, this can be done by a proper Riccati system.
Nevertheless, in the SVIEs (or even VIEs) case, it becomes quite challenging to seek appropriate feedback optimal control (see e.g. Lindquist \cite{Lindquist-1973}). We believe one important reason lies in the lack of \it flow property for SVIEs or VIEs. \rm
In order to overcome this essential difficulty, let us look at the following auxiliary function $\sX$:
\vskip-7mm
\bel{E3}\ba{ll}
\ns\ds\sX(s,t)=\varphi(s)+\int_\t^t\big[A(s,r)X(r)+B(s,r)u(r)+b(s,r)\big]dr\\[-1.8mm]
\ns\ds \qq\qq+\int_\t^t\big[C(s,r)X(r)+D(s,r)u(r)+\si(s,r)\big]dW(r),\;\;\; \t\les t < s\les T.
\ea\ee
\vskip-2mm
\no It is non-anticipating in the sense that $\sX(\cd,t)$ is determined by the information of $X$, $u$ up to the current time $t$. In addition, the flow property can be reestablished via $\sX$ (see e.g. Viens-Zhang \cite{Viens-Zhang-2019}).
According to monograph Corduneanu \cite{Corduneanu-2002}, we can regard $\sX(\cd,t)$ as a \emph{causal operator/abstract Volterra operator} acting on the function space of $(X, u)$, and name it the \emph{causal state trajectory} in our LQ problems. In the deterministic case, this term was used in LQ problems, see e.g.
Lee-You \cite{Lee-You-1989-AMO, Lee-You-1990-SICON}, Pritchard-You \cite{Pritchard-You-1996}, and recent Han et al \cite{Han-Lin-Yong-2023}. In the stochastic case, it was useful in Hamaguchi-Wang \cite{Hamaguchi-Wang-I, Hamaguchi-Wang-II}, Wang \cite{Wang-2016-submitted}, Wang et al \cite{Wang-Yong-Zhou-2023-SICON} for the LQ problems and Coutin-Decreusefond \cite{Coutin-Decreusefond}, Viens-Zhang \cite{Viens-Zhang-2019} for other problems.

\ss
 To explain the motivations of the current study, let us give a closer revisit to the existing literature on LQ problem of VIEs (e.g. De Acutis \cite{Acutis-1985}, Han et al \cite{Han-Lin-Yong-2023}, Lindquist \cite{Lindquist-1973}, Lee-You \cite{Lee-You-1989-AMO, Lee-You-1990-SICON}, Pritchard-You \cite{Pritchard-You-1996}, Pandolfi \cite{Pandolfi-2018-IEEETAC},
 Shaikhet \cite{Shaikhet-1988})
and that of SVIEs (e.g. Abi Jaber et al \cite{Jaber-Miller-Pham-2021}, Chen-Yong \cite{Chen-Yong-2007}, Hamaguchi-Wang \cite{Hamaguchi-Wang-I, Hamaguchi-Wang-II}, Wang et al \cite{Wang-Yong-Zhou-2023-SICON}). Let us demonstrate the details from four standpoints.

$\bullet$ As we know, the difference between optimal control problem of Bolza type and that of Lagrange type lies in the terminal cost.
On the one hand, the previous one provides a general framework to deal with more practical problem such as mean-variance optimization problem or
expected terminal utility maximization problems. On the other hand, in the literature of ODEs/SDEs, it seems that the Bolza and the Lagrange type can be easily transformed into each other, while in the case of singular VIEs, it becomes nontrivial, see \cite{Han-Lin-Yong-2023}.
Back to the LQ problem of SVIEs, Hamaguchi--Wang \cite{Hamaguchi-Wang-I, Hamaguchi-Wang-II} recently investigated the particular Lagrange type and characterized optimal causal feedback strategy via the above $\sX(\cd,\cd)$.
Then one may ask\\
(Q1): \emph{for SVIEs, is it a trivial extension from the Lagrange type to the Bolza type in the characterization of causal feedback strategy?}

$\bullet$    Mathematical speaking, the matrix-valued Riccati equation should be more natural and easily acceptable for the finite-dimensional Problem (LQ-SVIE). However, Abi Jaber et al \cite{Jaber-Miller-Pham-2021} studied the LQ probelm of convolution SVIE, a special form of (\ref{E1}), based on an infinite-dimensional lifting approach and ended up with operator-valued Riccati equations for feedback controls.
Wang et al  \cite{Wang-Yong-Zhou-2023-SICON} introduced an operator-valued path-dependent Riccati
equation to characterize the feedback representation of the optimal control, and their controlled SVIE, another special form of (\ref{E1}), is an essential particular integro-differential type.
One may ask\\
(Q2): \emph{can we use the matrix-value Riccati equation to represent optimal (causal) feedback strategy like that in the SDEs case (e.g. Sun-Yong \cite{Sun-Yong-2020, Sun-Yong-2014})}?

$\bullet$ Based on the semigroup theory approach, De Acutis \cite{Acutis-1985}  studied the representation of an optimal feedback control for a deterministic VIDE.
Pandolfi \cite{Pandolfi-2018-IEEETAC} obtained system of Riccati differential equations by dynamic programming method and state space approach. We found that both papers obtained the optimal control of the type
\vskip-6mm
\bel{E2402}
 u(t)=p_0(t)X(t)+\int^t_\t p_1(s,t)X(s)ds,\ \ t\in[\t,T],
\ee
\vskip-2mm
\no where $(p_0, p_1)$ is the solution to certain Riccati equations and $X$ corresponds to the optimal state.
Notice that similar form of (\ref{E2402}) also appeared in Shaikhet \cite{Shaikhet-1988}.
 Recall that the solution of VIEs usually has the memory effect,
the (\ref{E2402}) indicates that such a property in terms of state can also be represented in the non-anticipated/causal control.
A this moment, it then becomes quite natural to ask\\
(Q3):
\emph{For Problem (LQ-SVIE), can we obtain the optimal causal feedback control without $\sX$ in (\ref{E3}), but in the memory manner of (\ref{E2402})?}

$\bullet$ If (\ref{E1}) reduces to deterministic VIE, a so-called \emph{projection causality approach} was introduced in Pritchard--You \cite{Pritchard-You-1996},
and the optimal control was represented in linear causal feedback way with the feedback strategy determined by a Fredholm integral equation.
%
Notice that this feedback representation is different from (\ref{E2402}), and it relies on the above casual state $\sX(\cd,\cd)$.
Such Fredholm idea also happened in Lee--You \cite{Lee-You-1989-AMO, Lee-You-1990-SICON}. On the other hand, if (\ref{E1}) reduces to a SDEs, there are abundant papers on representing optimal controls via Riccati equations (see e.g. Sun-Yong \cite{Sun-Yong-2020, Sun-Yong-2014}).
At this moment, one may ask\\
(Q4): \emph{For Problem (LQ-SVIE), if we use the projection causality idea to obtain the corresponding Riccati system,  then what are the relationships between the Riccati system here and the Fredholm equation in \cite{Pritchard-You-1996}? }

\ss

The purpose of this paper is to give a unified treatment to the Problem (LQ-SVIE), from which the above four different questions can be answered carefully. Inspired by the previous works (e.g. \cite{Hamaguchi-Wang-I, Hamaguchi-Wang-II, Pritchard-You-1996, Pandolfi-2018-IEEETAC, Wang-Yong-Zhou-2023-SICON}),
we introduce a new feedback representation of control process:
\vskip-6mm
\bel{E4}
 u(t)\1n=\1n\Theta_1(t)X(t)\1n+\2n
\int^T_\t \3n\Theta_2(s,t)\big(X(s)\otimes\sX(s,t)\big)ds\1n+\1n\Theta_3(t)\sX(T,t)\1n+\1nv(t),\; t\in[\t,T],
\ee
\vskip-2mm
\no  where $\Theta_1,\Theta_2,\Theta_3$ and $v$ are deterministic functions merely depending on the data of Problem (LQ-SVIE) and
\vskip-6mm
\bel{equ01}
X(s)\otimes\sX(s,t)\=X(s)I_{[\t,t)}(s)+\sX(s,t)I_{[t,T]}(s),\q t,s\in[\t,T].
\ee
\vskip-2.5mm
\no Clearly, the representation (\ref{E4}) of $u$ at current time $t$ does not
 involve future information of the corresponding state process $X$.
In this sense, we call $(\Theta_1,\Theta_2,\Theta_3,v)$  a \emph{causal feedback strategy} (see Definition \ref{def2.1} for a precise definition).
In the Bolza setting, \cite{Wang-Yong-Zhou-2023-SICON} discussed the optimal feedback control by
$\sX(\cd,t)$, but not the case of $X(\cd)I_{[\t,t](\cd)}$. Here we add this part inspired by \cite{Acutis-1985, Pandolfi-2018-IEEETAC, Shaikhet-1988} and the above (Q3).
The main result is to characterize $(\Theta_1,\Theta_2,\Theta_3,v)$ by
introducing and discussing appropriate Riccati system which to our best knowledge is quite new in the literature. To indicate the unified powerfulness, we give detailed comparisons with the existing literature.

\ss

Now we highlight the innovations of this paper, together with some interesting phenomena revealed under our unified framework.

$\bullet$
From the notion viewpoint, the unified treatment requires us to introduce a new class of causal feedback strategy (\ref{E4}) and its characterization via the a new Riccati system (\ref{E4.1}). To this end, we define \emph{suitable space} (Definition \ref{def2.9} ) and \emph{appropriate multiplicative rule} (Definition \ref{def2.10}) which appear for the first time
in the literature.

$\bullet$
From  the methodology standpoint, our unified treatment provides a new approach to study the LQ problem of SVIEs.
For example, to obtain the explicit representation of optimal causal feedback strategy and
the uniqueness of Riccati system, we introduce \emph{strongly} optimal causal feedback strategy and the \emph{extended}
cost functional (see Remark \ref{remark2.7}, Remark \ref{remark3.8}, Remark \ref{remark3.9} for details) which reduce to the Lagrange case in \cite{Hamaguchi-Wang-I, Hamaguchi-Wang-II} and is unnecessary in the SDEs case. This shows the nontrivial extension to the literature, which answer the above (Q2) in a proper manner.

$\bullet$
From the conclusion standpoint, we give a unified treatment to the existing literature in several cases, including SVIEs \cite{Hamaguchi-Wang-II}, SVIDE \cite{Wang-Yong-Zhou-2023-SICON}, convolution SVIE \cite{Jaber-Miller-Pham-2021}, VIDEs \cite{Pandolfi-2018-IEEETAC}, VIEs \cite{Pritchard-You-1996}, SDEs in e.g. \cite{Sun-Yong-2020}, to mentioned a few.
We provide explicit representation of causal feedback strategies via finite dimensional Riccati system, which coincides with the Fredholm system in \cite{Pritchard-You-1996}, clarifies the operator-valued counterparts in  \cite{Jaber-Miller-Pham-2021, Wang-Yong-Zhou-2023-SICON}, and extends that in e.g. \cite{Hamaguchi-Wang-II, Pandolfi-2018-IEEETAC, Sun-Yong-2020}. This corresponds to the above (Q1), (Q2), (Q4).

$\bullet$ The unified treatment helps us to reveal an interesting phenomena, i.e.,
optimal control (if it exists) generally does not have representation in the spirit of (\ref{E2402}). It
gives us an answer to the above (Q3), and is based on three observations.
First, in the deterministic case ($B\1n=\1nD\1n=\1n0$), Theorem \ref{the3} shows $\Th_2(s,t),s\les t,$ plays role only when $R(\cd)$ is singular.
 Second, in particular case (see Corollary \ref{corollary3.10}), we prove it is contradictory once the above optimal representation (\ref{E2402}) exists.
 Third, in contrast with $X(\cd)$ in  (\ref{E4}), the term $\sX(\cd,t)$ is more essential and appropriate in the sense that the later reduces to the former in particular case (see Subsection 5.5).
 %

%
 %
 %

\ss

The rest of the paper is structured as follows. The next section is devoted to some preliminaries. In
Section 3 we establish an equivalent characterization of the existence of the strongly optimal feedback strategy.
In Section 4 we prove the well-posedness of the new  Riccati system.
 In Section 5 we apply our result to some special cases. Section 6 concludes this paper.

\section{Preliminaries}

Let $\mathbb{R}^d$ and $\mathbb{R}^{d\times l}$ be the usual $d$-dimensional space of real numbers and the set of all
$(d\times l)$ real matrices, $\mathbb{S}^d$ be the set of all $(d\times d)$ real symmetric matrices, and $I_d$ be the $(d\times d)$-identity matrix.
$M^{\top}$ and $M^\dagger$ denote the transpose and the Moore-Penrose pseudoinverse of matrix $M$, respectively. $\sR(M)$ denotes the range.
We shall use $|\cdot|$ and $\langle\cdot\rangle$ to denote the Euclidean norm and product.
 We define a  triangle  region  and a square pyramid region:
\vskip-6mm
$$\ba{ll}
 \triangle_\ast[0,T]\1n=\1n\big\{(r,s)\1n\in\1n[0,T]^2\big|0\1n\les\1n s\1n< \1nr\1n\les\1n T\big\},\;
\Box_3[0,T]\1n=\1n\big\{(s_1,s_2,t)\1n\in\1n[0,T]^3\big|t\1n<\1ns_1\wedge s_2\big\}.
\ea
$$
\vskip-2mm
 For any $\t\in[0,T]$ and Euclidean space $\dbH$ which could be $\dbR^{d\times l}$, or $\dbS^{d}$, and so on,
 we define several spaces.
Let $ L^\infty(\t,T;\dbH)$ be the set of $\dbH$-valued essentially bounded measurable functions,
$L^2(\triangle_\ast[\t,T];\1n\dbH)$ the set of $\dbH$-valued and square-integrable deterministic functions,
$L^{2,2,1}(\Box_3[\t,T];\1n\dbH)$ the set of $f:\1n\Box_3[\t,T]\1n\rightarrow\1n\dbH$ such that
$\1n\int^T_\t\2n\int^T_\t\2n\big(\1n\int^{s_1\wedge s_2}_\t\1n|f(s_1,\1ns_2,\1nt)|dt\big)^2$ $\1nds_1\1nds_2\2n<\2n\infty\1n$,
$L^2_{\mathbb{F}}(\t,\1nT;\1n\dbH)$ the set of $\dbH$-valued, square-integrable and $\mathbb{F}$-progressively measurable processes,
$L^2_{\mathbb{F}}(\Omega;C([\t,T];\1n\dbH))$ the set of $f\1n: \1n[\t,T]\1n\times\1n \Omega\1n\rightarrow\1n\dbH$
such that $t\1n\mapsto\1n f(t,\omega)$  is continuous and $\dbE(\1n\sup\limits_{t\in[\t,T]}|f(t)|^2)\1n<\1n\infty$,
%
%
%
 $L^2_{\mathbb{F}}( \triangle_\ast[\t,T];\dbH)$ the set of $f\1n:\1n \triangle_\ast[\t,T]\times\Omega\rightarrow\dbH$ such that
 $f(t,\cdot)\1n\in\1n L^2_{\mathbb{F}}(\t,t;\dbH)$ for  $\1nt\1n\in\1n[\t,T]$ and $\dbE\2n\int^T_\t\2n\int^t_\t\1n|f(t,s)|^2dsdt<\infty$,
$L^{2}_{\mathbb{F},c}( \triangle_\ast[\t,T];\dbH)$ the set of $f\in L^2_{\mathbb{F}}( \triangle_\ast[\t,T];\dbH)$
such that $ s\1n\rightarrow\1n f(t,s) $ is uniformly continuous on $(\t,t)$ and $\dbE\2n\int^T_\t\1n\sup\limits_{s\in[\t,t]}|f(t,s)|^2dt\1n<\1n\infty$.
Inspired by \cite{Hamaguchi-Wang-I, Hamaguchi-Wang-II}, let $\sL^2(\triangle_\ast[\t,T];\dbH)$ be the set of
$f\1n\in\1n L^2(\triangle_\ast[\t,T];\dbH)$ satisfying $\esssup\limits_{t\in[\t,T]}\big($ $\1n\int^T_t\1n|f(s,t)|^2ds\big)^{\frac{1}{2}}\1n<\1n\infty$ such that for any $\varepsilon>0$, there exists a finite partition $\{U_i\}^m_{i=0}$ of $[\t,T]$ with $\t=U_0<U_1<\cdots<U_m=T$ satisfying
$\esssup\limits_{t\in(U_i,U_{i+1})}\big(\1n\int^{U_{i+1}}_t\1n|f(s,t)|^2ds\big)^{\frac{1}{2}}\1n<\1n\varepsilon,\;i\1n\in\1n\{0,1,\ldots,m-1\}.$

\ms
We impose the following assumptions throughout this paper:

\medskip

\no\textbf{(H1)}. The coefficients $A,C\1n:\1n \triangle_\ast[0,T]\1n\rightarrow\1n\mathbb{R}^{d\times d}$ and
$B,D\1n:\1n\triangle_\ast[0,T]\1n\rightarrow\1n\mathbb{R}^{d\times l}$ are bounded.
 $b, \si \1n:\1n\triangle_\ast[0,T]\times\Omega\1n\rightarrow\1n\mathbb{R}^{d}$ are bounded and $ s\1n\rightarrow\1n \big(b(t,s),\si(t,s)\big)$
is $\mathbb{F}$-progressively measurable on $(0,t)$.\\
\textbf{(H2)}. $ Q\2n\in\2n L^\infty(0,T;\1n\mathbb{S}^{d}),R\2n\in\2n L^\infty(0,T;\mathbb{S}^{l}),G\2n\in\2n\mathbb{S}^{d},$
 $ S\2n\in\2n L^\infty(0,T;\mathbb{R}^{l\times d})$,$q\2n\in\2nL^2_{\mathbb{F}}(0,T;\mathbb{R}^{d}),\rho\2n\in\2nL^2_{\mathbb{F}}(0,T;\mathbb{R}^{l})$.\\

For notational simplicity, we define   $\cI\=\big\{(\t,\f)\big| \t\in[0,T), \varphi\in C(\t,T;\dbR^d)\big\}$ and
\vskip-5.5mm
$$\ba{ll}
\ns\ds \mathcal{S}[0,T]\=L^\infty(0,T;\mathbb{R}^{l\times d})\times L^2([0,T]^2;\mathbb{R}^{l\times d})\times L^\infty(0,T;\mathbb{R}^{l\times d})\times\mathcal{U}[0,T].
\ea
$$
\vskip-4mm
\begin{definition}\label{def2.1}
\vskip-3mm
(i) \no Any  4-tuple $(\Theta_1,\Theta_2,\Theta_3,v)\in\mathcal{S}[0,T]$ is called a causal feedback strategy.\\
(ii) For any $(\Theta_1,\Theta_2,\Theta_3,v)\in\mathcal{S}[0,T]$ and $(\t,\varphi)\in\cI$,
the pair  $(X,\sX)\1n\in\1n L^2_{\mathbb{F}}(\t,T;$ $\mathbb{R}^{ d})\1n\times \1nL^2_{\mathbb{F},c}(\triangle_\ast[\t,T];\mathbb{R}^{ d})$
 is called the causal feedback solution of (\ref{E1}) if it satisfies
\vskip-5.5mm
\bel{E5}\left\{\2n\ba{ll}
\ds X(t)=\1n\varphi(t)\1n+\2n\int_\t^t\2n\big[A(t,r)X(r)\1n
+B(t,r)u^{\Theta,v}(r)+b(t,r)\big]dr\\[-1.5mm]
\ns\ds\qq\q+\2n\int_\t^t\2n\big[C(t,r)X(r)\1n
+\1nD(t,r)u^{\Theta,v}(r)\1n+\1n\si(t,r)\big]dW(r),t\1n\in\1n[\t,T],\\[-1.5mm]
\ns\ds \sX(s,t)\1n=\1n\varphi(s)\1n+\2n\int_\t^t\2n\big[A(s,r)X(r)\1n
+B(s,r)u^{\Theta,v}(r)+b(s,r)\big]dr\\[-1.5mm]
\ns\ds\qq+\2n\int_\t^t\2n\big[C(s,r)X(r)\1n
+\1nD(s,r)u^{\Theta,v}(r)\1n+\1n\si(s,r)\big]dW(r),(s,t)\1n\in\1n\triangle_\ast[\t,T],\\[-1.5mm]
\ns\ds u^{\Theta,v}(t)\1n\=\1n\Theta_1(t)X(t)\2n+\2n
\int^T_\t\3n \Theta_2(r,t)\big(X(r)\1n\otimes\1n\sX(r,t)\big)dr\2n+\2n\Theta_3(t)\sX(T,t)\2n+\2nv(t),\;t\1n\in\1n[\t,T].
\ea\right.
\ee
\vskip-2mm
\no Here, $X\otimes\sX$ is in (\ref{equ01}), $u^{\Theta,v}$ is called the outcome of $(\Theta_1,\Theta_2,$ $\Theta_3,v)$ at $(\t,\varphi)$.

\end{definition}
\begin{remark}\label{Remark-difference-strategy}
As to the above definition, let us point out two interesting facts. First, we observe  $(\Theta_1,\Theta_2,\Theta_3,v)\in\mathcal{S}[0,T]$ does not rely on $(\t,\varphi)\in\cI$
while $u^{\Th,v}$ does.
Second, there are some important differences from \cite{Hamaguchi-Wang-I, Hamaguchi-Wang-II} even when the terminal cost disappears.
In fact, our strategy $\Th_2$ is defined on $[0,T]^2$ instead of $\D_*[0,T]$. Therefore, the $u^{\Th,v}$ in (\ref{E5}) have additional dependence on the value of $X(\cd)$ in $[\t,t]$.
\end{remark}

The following lemma gives the existence and uniqueness of a causal feedback solution. Since its proof is almost the same as
 \cite[Theorem 2.4]{Hamaguchi-Wang-I}, we omit the details.
\vskip-4mm
\bl\label{lem2.2}
Let $(\rm H1)$-$(\rm H2)$ hold. For each causal feedback strategy $(\Theta_1,\Theta_2,\Theta_3,v)$ $\in\mathcal{S}[0,T]$ and
each input condition $(\t,\varphi)\in\cI$, the controlled SVIE (\ref{E1}) has a unique causal feedback solution
$(X,\sX)\1n\in\1n L^2_{\mathbb{F}}(\t,T;\mathbb{R}^{ d})\1n\times \1nL^2_{\mathbb{F},c}(\triangle_\ast[\t,T];\mathbb{R}^{ d})$.
\el

If $(X,\sX)$ is the causal feedback solution of (\ref{E1}), then it is easy to see
$u^{\Theta,v}$ belongs to $\mathcal{U}[\t,T]$.
At this moment, let us present the definition of optimal causal feedback strategy.

\begin{definition}\label{def2.6}
A 4-tuple $(\widehat{\Theta}_1,\widehat{\Theta}_2,\widehat{\Theta}_3,\widehat{v})\in\mathcal{S}[0,T]$ is called optimal causal feedback strategy of Problem (LQ-SVIE) if
\vskip-4.5mm
$$\ba{ll}
\cJ(\t,\varphi;u^{\widehat{\Theta},\widehat{v}})\les \cJ(\t,\varphi;u), \ \ \forall (\t,\varphi)\in\cI,\  \ u\in\mathcal{U}[\t,T].
\ea
$$
\vskip-6mm
\end{definition}
\vskip-2.5mm
Since the cost functional depends on the terminal state, we have to introduce some new concepts. Instead of the above \it original framework \rm (including notations $\cI$, $\cJ$, etc,) with optimal strategy, we will work in the \it extended framework \rm with \it stronger optimality \rm to obtain the explicit forms and uniqueness of causal feedback strategy (once it exists). We refer to the following Remark \ref{remark2.7}, Remark \ref{remark3.8} and Remark \ref{remark3.9} for more detailed explanations.
To begin with, inspired by \cite{Wang-2020}, we define
$$\widetilde{\cI}\=\big\{(\t,\varphi_1,\varphi_2)|\ \t\1n\in\1n[0,T),\varphi_1\1n\in\1n L^2(\t,T;\mathbb{R}^{d}),\ \varphi_2\in\mathbb{R}^{d}\big\}.$$
For any $(\t,\varphi_1,\varphi_2)\in\widetilde{\cI}$ and $(\Theta_1,\Theta_2,\Theta_3,v)\in\mathcal{S}[0,T]$,
consider the following system:
\vskip-4.5mm
\bel{E0601}\left\{\ba{ll}
\ns\ds X_1(t)=\1n\varphi_1(t)\1n+\2n\int_\t^t\2n\big[A(t,r)X_1(r)\1n
+B(t,r)u_1^{\Theta,v}(r)+b(t,r)\big]dr\\[-2mm]
\ns\ds\qq\qq+\2n\int_\t^t\2n\big[C(t,r)X_1(r)\1n
+\1nD(t,r)u_1^{\Theta,v}(r)\1n+\1n\si(t,r)\big]dW(r),t\1n\in\1n[\t,T),\\[-2mm]
\ns\ds \sX_1(s,t)\1n=\1n\varphi_1(s)\1n+\2n\int_\t^t\2n\big[A(s,r)X_1(r)\1n
+B(s,r)u_1^{\Theta,v}(r)+b(s,r)\big]dr\\[-2mm]
\ns\ds\qq\q+\2n\int_\t^t\2n\big[C(s,r)X_1(r)\1n+\1nD(s,r)
u_1^{\Theta,v}(r)\1n+\1n\si(s,r)\big]dW(r),\t\1n\les \1nt\1n< \1ns<\1n T,\\[-2mm]
\ns\ds \sX_2(t)\1n=\1n\varphi_2+\2n\int_\t^t\2n\big[A(T,r)X_1(r)\1n
+B(T,r)u_1^{\Theta,v}(r)+b(T,r)\big]dr\\[-1mm]
\ns\ds\qq\qq+\2n\int_\t^t\2n\big[C(T,r)X_1(r)\1n+\1nD(T,r)
u_1^{\Theta,v}(r)\1n+\1n\si(T,r)\big]dW(r),t\1n\in\1n[\t,T],\\[-2mm]
\ns\ds u_1^{\Theta,v}(t)\2n\=\2n\Theta_1(t)X_1(t)\2n+\3n
\int^T_\t\3n\Theta_2(r,t)\big(X_1(r)\1n\otimes\1n\sX_1(r,t)\big)dr\2n+\2n\Theta_3(t)\sX_2(t)\2n+\2nv(t),t\1n\in\1n[\t,T).
\ea\right.
\ee
\vskip-1mm
\no A method similar to \cite[Theorem 2.4]{Hamaguchi-Wang-I} ensures the existence and uniqueness of $(X_1,\sX_1,\sX_2)\in L^2_{\mathbb{F}}(\t,T;\mathbb{R}^{ d})\times L^2_{\mathbb{F},c}(\triangle_\ast[\t,T];\mathbb{R}^{ d})
\1n\times\1n L^2_{\mathbb{F}}(\Omega;C([\t,T];\mathbb{R}^{ d}))$.
We call $(X_1,\sX_1,\sX_2)$ \emph{the extended causal feedback solution} of controlled SVIE, $u_1^{\Theta,v}$
the \emph{extended outcome} of  $(\Theta_1,\Theta_2,$ $\Theta_3,v)$ at $(\t,\varphi_1,\f_2)$, and (\ref{E0601}) the \it extended closed-loop system. \rm

For any $(\t,\varphi_1,\varphi_2)\in\widetilde{\cI}$ and $u\in\mathcal{U}[\t,T]$, we further introduce
\vskip-6.5mm
\bel{E1.9}\ba{ll}
\ns\ds\widetilde{\cJ}(\t,\varphi_1,\varphi_2;u)\1n\=\1n\dbE\Big[\langle G\sX_2(T),\1n\sX_2(T)\rangle\1n+
\3n\int_\t^T\3n\Big(\langle Q(s)X_1(s),\1nX_1(s)\rangle\1n+\1n\langle R(s)u(s),\1nu(s)\rangle\\[-1mm]
 \ns\ds \qq\q+2\langle S(s)X_1(s),u(s)\rangle\1n+2\langle q(s) ,X_1(s)\rangle\1n+2\langle
 \rho(s),u(s)\rangle \Big) ds\Big],
\ea\ee
\vskip-3mm
\no where $X_1,\sX_1$ and $\sX_2$ are the solutions of the first, second and third equations in (\ref{E0601})
replacing $u_1^{\Theta,v}$ with any $u\in\mathcal{U}[\t,T]$, respectively.
We name (\ref{E1.9}) the \emph{extended cost functional} and define the corresponding \emph{extended value function}
\vskip-3mm
$$\ba{ll}
\wt V(\t,\f_1,\f_2) \=\inf\limits_{u\in \mathcal{U}[\t,T]}\widetilde{\cJ}(\t,\varphi_1,\varphi_2;u).
\ea$$
\vskip-1mm
\no Similar to Definition \ref{def2.6}, we give the following definition.
\begin{definition}
A 4-tuple
$(\widehat{\Theta}_1,\widehat{\Theta}_2,\widehat{\Theta}_3,\h v)\in\mathcal{S}[0,T]$ is called strongly optimal causal feedback
 strategy of Problem (LQ-SVIE) if
\vskip-4mm
$$\ba{ll}
\widetilde{\cJ}(\t,\varphi_1,\varphi_2;u_1^{\widehat{\Theta},\widehat{v}})
\les\widetilde{\cJ}(\t,\varphi_1,\varphi_2;u), \ \ \forall (\t,\varphi_1,\varphi_2)\in\widetilde{\cI},\;\; u\in\mathcal{U}[\t,T].
\ea
$$
\vskip-2mm
\end{definition}
\begin{remark}\label{remark2.7}
\vskip-2mm
At this moment, let us explain the terms \it extended \rm and \it strongly \rm in the above.

First, for any $(\t,\f)\in\cI$, by taking $\f_1\=\f$  and $\f_2\=\f(T)$, we see that the extended causal feedback solution in (\ref{E0601})
naturally reduces to causal feedback solution in (\ref{E5}). Similar principle also holds for the outcome, the closed-loop system, the cost functional and the value function.

Second, since the set $\wt \cI$ is \it larger \rm than $\cI$, it is clear that a strongly optimal causal feedback strategy  is also optimal in the spirit of Definition \ref{def2.6}. This explains the term \it strongly. \rm
If $G\1n=\1n0$, and $\f\in L^2(\t,T;\dbR^d)$, the two concepts of strongly optimal causal feedback strategy and optimal causal feedback strategy are equivalent.

Third, in the SDEs case when the free term $\f$ reduces to vector in $\dbR^d$, the above extended framework is not necessary. From this sense, we capture an essential difference caused by the extension from finite dimensional Euclidean space $\dbR^d$ to infinite dimensional space $C([0,T];\dbR^d)$.

\end{remark}

\ss

To introduce the Riccati type system, we define a suitable space and certain right (or left) multiplicative rule.
Both of them appear for the first time in the literature.
\begin{definition}\label{def2.9}
Let $\Upsilon[0,T]$ be the set of  $(\sP_1,\1n\sP_2,\1n\sP_3,\1n\sP_4)$ with $\sP_1,\1n\sP_2: [0,T]\rightarrow \dbS^d$,
$\sP_3:\triangle_*[0,T]\rightarrow \dbR^{d\times d}$
and  $\sP_4: \Box_3[0,T]\rightarrow \dbR^{d\times d}$ such that

$(\rm i)$  $\sP_1,\sP_2 \1n\in  \1n L^\infty(0,T;\mathbb{S}^d)$,
$\sP_3 \1n\in \1n\sL^2(\triangle_\ast[0,T];\mathbb{R}^{d\times d}),$  $\sP_4\1n\in \1nL^{2,2,1}(\Box_3[0,T];\mathbb{R}^{d\times d})$;

$(\rm ii)$  for a.e. $\1n(s_1,s_2)\1n\in\1n [0,T]^2$,  $t\1n\mapsto\1n \sP_4(s_1,s_2,t)$ is absolutely continuous on $(0,s_1\wedge s_2)$,
and for a.e. $\1n(s_1,s_2,t)\1n\in \1n\Box_3[0,T]$, it holds that $\sP_4(s_1,s_2,t)\1n=\1n\sP_4(s_2,s_1,t)^{\top}$.

\end{definition}
The above $\Upsilon[0,T]$ will be the solution space for our Riccati system. Next we define some multiplicative rules that are frequently used in the Riccati system.
\begin{definition}\label{def2.10}
Let $\sP=(\sP_1,\sP_2,\sP_3,\sP_4)\1n\in\1n\Upsilon[0,T]$. For each $M:\1n\triangle_*[0,T]\1n\rightarrow \dbR^{d_1\times d}$ and
$N:\triangle_*[0,T]\rightarrow \dbR^{d\times d_2}$, we define:
\vskip-7.5mm
$$\ba{ll}
\ns\ds (M\triangleleft \sP_{2,3})(t)\=
M(T,t)\sP_2(t)\1n+\2n\int_t^TM(s,t)\sP_3(s,t)^{\top}ds,\;\;t\in[0,T],\\[-0.6mm]
\ns\ds (\sP_{2,3}\triangleright N)(t)\=
\sP_2(t)N(T,t)+\int_t^T\sP_3(s,t)N(s,t)ds,\;\;t\in[0,T],\\[-0.6mm]
\ns\ds (M\triangleleft \sP_{1,3,4})(s,t)\2n\=\2n
M(s,t)\sP_1(s)\1n+\1nM(T,t)\sP_3(s,t)\1n+\2n\int_t^T\3n\1nM(r,t)\sP_4(r,s,t)dr,(s,t)\1n\in\1n \triangle_{*}[0,T],\\[-0.6mm]
\ns\ds (\sP_{1,3,4}\triangleright N)(s,t)\2n\=\2n
\sP_1(s)N(s,t)\1n+\1n\sP_3(s,t)^{\top}N(T,\1nt)\1n+\2n\int_t^T\3n\1n\sP_4(s,r,t)N(r,t)dr,(s,t)\1n\in \1n\triangle_{*}[0,T],\\[-0.6mm]
\ns\ds (M\triangleleft \sP\triangleright N)(t)\2n\=\2nM(T,t)\1n\sP_2(t)N(T,t)\1n+\2n\int_t^T\3n\Big[M(s,t)\sP_1(s)N(s,t)
\1n+\1nM(T,t)\sP_3(s,t)N(s,t)\\[-0.6mm]
\ns\ds\qq\qq\q+M(s,t)\sP_3(s,t)^{\top}N(T,t)\1n+\2n\int_t^T\3nM(s,t)\sP_4(s,\th,t)N(\theta,t)d\th\Big]ds,\;t\in[0,T].
\ea
$$
\vskip-3mm
\end{definition}
\vskip-3mm
\begin{remark}\label{rem2.9}
\vskip-6mm
Notice that if we take $M(s,t)\equiv M(t), N(s,t)\equiv N(t)$, then $``\triangleleft"$ and $``\triangleright"$ in Definition \ref{def2.10}
reduces to the classical multiplicative rule, e.g.,
\vskip-7.5mm
$$\ba{ll}
\ns\ds (M\triangleleft\sP\triangleright N)(t)\2n=\2nM(t)\Big\{\sP_2(t)\1n+\2n\int_t^T\3n\1n\Big(\sP_1(s)\1n+\1n
\sP_3(s,t)\1n+\1n\sP_3(s,t)\1n^{\top}\2n+\2n\int_t^T\3n\2n\sP_4(s,\th,t)d\th\Big)ds\Big\}N(t).
\ea
$$
In addition, if $(\sP_2,\sP_3)\equiv (0,0)$, then $(\1nM\2n\triangleleft \1n\sP_{1,3,4}\1n)$, $(\1n\sP_{1,3,4}\1n\triangleright \2nN\1n)$ and $(M\triangleleft\1n \sP\triangleright\1n N)$ reduce to multiplicative rules in \cite{Hamaguchi-Wang-I, Hamaguchi-Wang-II}.
\vskip-3mm

\end{remark}
\vskip-3mm
At this moment, it is time for us to present the following system:
\vskip-4mm
\bel{E4.1}\3n\3n\left\{\2n\ba{ll}
\ns\ds \1n\sP_1(t)=Q(t)+(C^{\top}\1n\triangleleft \sP \triangleright  C)(t)\\
\ns\ds\q+\big( S(t)^{\top} \2n+ (C^{\top} \1n\triangleleft \sP\triangleright  D)(t)\big)
\BR(t)^\dagger\big((D^{\top}\1n\triangleleft \sP \triangleright  C)(t)+S(t)\big),\;t\in[0,T],\\
\ns\ds \1n\sP_2(t)\2n=G\1n-\1n\int_t^T\2n(\sP_{2,3}\triangleright B)(s)\BR(s)^\dagger
(B^{\top}\1n\triangleleft \sP_{2,3})(s)ds,\;t\in[0,T],\\
\ns\ds\1n\sP_3(t,r)\2n=\2n(\sP_{2,3}\1n\triangleright A)(t)\1n-\1n(\sP_{2,3}\1n\triangleright B)(t)
\BR(t)^\dagger\big((D^{\top}\1n\triangleleft \sP \triangleright  C)(t)+S(t)\big)\\
\ns\ds\q-\2n\int_r^t\2n(\sP_{2,3}\1n\triangleright B)(\th)\BR(\th)^\dagger
(B^{\top}\1n\triangleleft \sP_{1,3,4})(t,\th)d\th,\;(t,r)\in\triangle_{*}[0,T],\\
\ns\ds \1n\sP_4(s,t,r)\1n=\1n\sP_4(t,s,r)^{\top} \\
\ns\ds\q =\1n(\sP_{1,3,4}\1n\triangleright \1nA)(s,\1nt)\2n-\2n(\sP_{1,3,4}\1n\triangleright\1n B)(s,\1nt)\BR(t)^\dagger
\big((D\1n^{\top}\3n\triangleleft \1n\sP\1n \triangleright \1n C)(t)+S(t)\big) \\
\ns\ds\q -\2n\int_r^t\2n(\sP_{1,3,4}\1n\triangleright B)(s,\th)\BR(\th)^\dagger
(B^{\top}\2n\triangleleft\sP_{1,3,4})(t,\th)d\th,\;0\1n\les \2nr\2n\les \2nt\2n\les \2ns\2n \les\2n T,
\ea\right.
\ee
\vskip-1.5mm
\no  %
\vskip-1.5mm
\no  where $\BR(\cd)\1n\=\1nR(\cd)\1n+\1n(D^{\top}\1n\triangleleft \sP \triangleright   D)(\cd)$.
In Subsection 5.1 and 5.2, we shall show that (\ref{E4.1}) reduces to the so-called \emph{Riccati-Volterra equation} in \cite{Hamaguchi-Wang-II}
and the classical Riccati equation in LQ problem for SDE, respectively.
Therefore, in the current paper we will name (\ref{E4.1}) a \emph{Riccati system}.
Following the terminology in \cite{Hamaguchi-Wang-II}, we give the concept
of the regular solution for (\ref{E4.1}). Due to notational simplicity, we denote
$\BD(\cd)\1n\=\1n S(\cd)+(D\1n^{\top}\2n\triangleleft\1n \sP\1n \triangleright \1n C)(\cd)$,
$\BB_1(\cd,\cd)\1n\=\1n(B^{\top}\2n\triangleleft\1n \sP_{1,3,4})(\cd,\cd)$,
$\BB_2(\cd)\1n\=\1n(B^{\top}\2n\triangleleft \1n\sP_{2,3})(\cd)$.
\begin{definition}\label{def2.11}
Let $\sP\1n=\1n(\sP_1,\1n\sP_2,\1n\sP_3,\1n\sP_4)\1n\in\1n \Upsilon[0,T]$ be a solution to the Riccati system (\ref{E4.1}).
It is called regular if

$(\rm i)$  $\BR(t)\ges 0$ for a.e. $t\in[0,T]$,

$(\rm ii)$ $\sR\big(\BD(t)\2n+\2n\BB_2(t)\big)\1n\subset\1n\sR(\BR(t))$ for a.e. $\1nt\2n\in\2n[0,\1nT]$
and $\sR(\BB_1(r,t))\2n\subset\2n\sR(\BR(t))$ for a.e. $\1n(r,t)\2n\in\2n\triangle_\ast[0,\1nT]$,

$(\rm iii)$  $(\BR^\dagger \BD,\BR^\dagger \BB_1,\BR^\dagger \BB_2)\in L^\infty(0,T;\mathbb{R}^{l\times d})\times L^2(\triangle_\ast[0,T];\mathbb{R}^{l\times d})
\times L^\infty(0,T;\mathbb{R}^{l\times d})$.

\end{definition}

\section{Strongly optimal strategy and Riccati equation}

This section is devoted to the equivalence (Theorem \ref{the3}) between the existence of strongly optimal causal feedback strategy and the regular solvability of
 Riccati system (\ref{E4.1}).
To achieve this goal,  we first give a representation of the extended cost functional (Lemma \ref{the3.2}).

To begin with, for any $(\t,\varphi_1,\f_2)\in\wt\cI$ and $(\Theta_1,\Theta_2,\Theta_3,v)\in\mathcal{S}[0,T]$, we introduce the following SVIE:
\vskip-6.5mm
\bel{E7}
\dbX(t)\2n=\2n\Phi^{\f}(t)\1n+\2n\int_\t^t \3n\big(\dbA(t,s)\dbX(s)+\Phi^{b}(t,s)\big)ds
\1n+\2n\int_\t^t\3n \big(\dbC(t,s)\dbX(s)+\Phi^{\si}(t,s)\big)dW(s),t\1n\in\1n[\t,T].
\ee
\vskip-1mm
\no Here, $\Phi^{\f},\Phi^{b},\Phi^{\si},\dbA,\dbC$ are defined by:
\vskip-5.5mm
\bel{E2.6}
\Phi^{\f}(t)\=\left(\begin{array}{ccccc}
\varphi_1(t)   \\[-1.8mm]
\ns\ds  \1n\Theta_1(t)\varphi_1(t)\1n+\2n\int^T_\t \2n\Theta_2(r,t)\varphi_1(r)dr\1n +\1n\Theta_3(t)\varphi_2\1n+\1nv(t)\\[-1.2mm]
\varphi_2
\end{array}
\right), \q t\in[\t,T],
\ee
\vskip-3mm
\no and
\vskip-7mm
\bel{AB12}\ba{ll}
\ns\ds \dbA(t,s)\=\left(\begin{array}{ccccc}
  A(t,s)  &  B(t,s)  & 0  \\
A_{12}(t,s)  & B_{12}(t,s) & 0\\
  A(T,s)  &  B(T,s)  & 0
\end{array}
\right),
\ \
\dbC(t,s)\=\left(\begin{array}{ccccc}
  C(t,s)  &  D(t,s)   & 0 \\
C_{12}(t,s)  & D_{12}(t,s) & 0\\
  C(T,s)  &  D(T,s)   & 0
\end{array}
\right),
\ea\ee
\vskip-2.5mm
\no and
\vskip-7mm
\bel{b12}\ba{ll}
\ns\ds \Phi^{b}(t,s)\=\left(\begin{array}{ccccc}
  b(t,s)    \\
b_{12}(t,s) \\
 b(T,s)
\end{array}
\right),
\ \
\Phi^{\si}(t,s)\=\left(\begin{array}{ccccc}
  \si(t,s) \\
\si_{12}(t,s)\\
\si(T,s)
\end{array}
\right),
\ea\ee
\vskip-1.5mm
\no where  for $f\=A,B,C,D,b,\si$,
\vskip-6mm
\bel{E3.3}
 f_{12}(t,s)\1n\=\1n\Theta_1(t)f(t,s)\1n+\2n\int^T_s\2n \Theta_2(r,t)f(r,s)dr\1n+\1n\Theta_3(t)f(T,s),\;\;(t,s)
 \1n\in\1n\triangle_\ast[\t,T].
\ee
\begin{remark}
\vskip-1.5mm
Using the Cauchy-Schwarz inequality and (H1), it is easy to see that
$\Phi^{\f}\1n\in\1n L^2_{\mathbb{F}}(\t,T;$ $\mathbb{R}^{2d+l})$, $\Phi^b,\Phi^{\si}\1n\in\1n L^2_{\mathbb{F}}(\triangle_\ast[\t,T];\mathbb{R}^{(2d+l)})$
and $\dbA,\dbC\1n\in\1n\sL^2(\triangle_\ast[\t,T];\mathbb{R}^{(2d+l)\times(2d+l)})$. Similarly to \cite[Theorem 2.4]{Hamaguchi-Wang-I},
we can show the well-posedness of above system (\ref{E7}). Here, for the sake of page limit, we omit the proof.
\end{remark}

\ss

Next we give a result which reveals the structure of the extended causal feedback solution.
Its proof is slight adjustment to \cite[Lemma A.4]{Hamaguchi-Wang-I}. We omit the details.
\bl\label{lem2}
Let $(\rm H1)$-$(\rm H2)$ hold and $(\Theta_1,\Theta_2,\Theta_3,v)\in\mathcal{S}[0,T]$ be fixed. For any $(\t,\varphi_1,\f_2)\in\wt\cI$,
$(X_1,\sX_1,\sX_2)\in L^2_{\mathbb{F}}(\t,T;\mathbb{R}^{ d})\times L^2_{\mathbb{F},c}(\triangle_\ast[\t,T];\mathbb{R}^{ d})\times L^2_{\mathbb{F}}(\Omega;C([\t,T];\mathbb{R}^{ d}))$ is a
extended causal feedback solution if and only if the following equalities hold:
\vskip-6mm
\bel{E6}\left\{\2n\ba{ll}
\ns\ds \left(\begin{array}{ccccc}
X_1(t)^{\top}\2n &   u^{\Theta,v}_1(t)^{\top}\2n& \sX_2(t)\1n^{\top}\2n
\end{array}\right)=\dbX(t)^{\top},\;\;\;t\in[\t,T],\\[-2mm]
\ns\ds \sX_1(s,t)=\varphi_1(s)+\int_\t^t\big((A(s,r)\enspace B(s,r)\enspace 0)\dbX(r)+b(s,r)\big)dr\\[-1.8mm]
\ns\ds \qq\q+\1n\int_\t^t\1n\big((C(s,r)\enspace  D(s,r)\enspace  0)\dbX(r)+\si(s,r)\big)dW(r),\;(s,t)\in\triangle_\ast[\t,T],
\ea\right.
\ee
\vskip-1.5mm
\no where $\dbX\in L^2_{\mathbb{F}}(\t,T;\mathbb{R}^{2d+l})$ is the solution of the SVIE (\ref{E7}).

\el

\ss

The above Lemma \ref{lem2} helps us to provide some new ideas different from  \cite[Theorem 4.1]{Hamaguchi-Wang-II} in the
representation of the extended cost functional (see Remark \ref{Distinction-H-W}).
 Next, for each $(\Theta_1,\Theta_2,\Theta_3)\1n\in\1n L^\infty(0,T;\mathbb{R}^{l\times d})\1n\times\1n L^2([0,T]^2;\mathbb{R}^{l\times d})\1n\times\1n
 L^\infty(0,T;\mathbb{R}^{l\times d})$,
 we introduce a deterministic system inspired by \cite[Theorem 4.1]{Wang-Yong-2023-SICON}:
\vskip-6mm
\bel{P12}\3n\1n\left\{\3n\1n\ba{ll}
\ns\ds P_1(t)\2n=\2n\dbQ(t)\2n+\3n\int_t^T\3n\1n\dbC(s,\1nt)\1n^{\top}\1nP_1(s)\dbC(s,\1nt)ds\2n
+\3n\int_t^T\3n\int_{t}^T\3n\1n\dbC(s,t)\1n^{\top}\1nP_2(\theta,s)\dbC(\theta,t)d\theta ds,t\1n\in\1n[0,\1nT],\\[-1.5mm]
\ns\ds P_2(r,t)=\dbA(r,t)^{\top}P_1(r)+\int_{t}^T\dbA(s,t)^{\top}P_2(r,s)ds\\[-0.5mm]
\ns\ds \qq\qq\3n =P_2(t,r)^{\top},\;(r,t)\1n\in\1n \triangle_{*}[0,T],
\ea\right.
\ee
\vskip-1.5mm
\no  where $\dbA, \dbC$ are defined by (\ref{E3.3})  and
\vskip-5mm
\bel{E3.1}
 \mathbb{Q}(t)\1n\= \1n\left(\begin{array}{ccccc}
  Q(t)+C(T,t)^{\top}GC(T,t) \1n &\1n   C(T,t)^{\top}GD(T,t)+S(t)^{\top}  \1n &\1n A(T,t)^{\top}G  \\
S(t)+D(T,t)^{\top}GC(T,t) \1n  & \1nR(t)+D(T,t)^{\top}GD(T,t) \1n&\1n B(T,t)^{\top}G \\
  GA(T,t) \1n &\1n GB(T,t) \1n &\1n 0
\end{array}
\right).
\ee
\vskip-1.5mm
\no By the Cauchy-Shwarz inequality, $(\rm H1)$ and $(\Theta_1,\Theta_2,\Theta_3,v)\in\mathcal{S}[0,T]$, it is easy to show that $\dbQ\1n\in\1n L^\infty(0,T;\1n\mathbb{S}^{2d+l})$. Similarly to the proof of \cite[Theorem 4.11]{Hamaguchi-Wang-I},
we can obtain the well-posedness of (\ref{P12}) in the sense that
$(P_1,\1nP_2)\1n\in\1n L^\infty(0,T;\mathbb{S}^{2d+l})\1n\times \1nL^2([0,T]^2;\mathbb{R}^{(2d+l)\times(2d+l)})$
with $P_2(s,t)\1n=\1nP_2(t,s)\1n^{\top}\1n$.

\ss

Furthermore, for each $(\Theta_1,\Theta_2,\Theta_3)\1n\in\1n L^\infty(0,T;\mathbb{R}^{l\times d})\1n\times\1n L^2([0,T]^2;\mathbb{R}^{l\times d})\1n\times\1n
 L^\infty(0,T;\mathbb{R}^{l\times d})$ and
any $(\xi,\gamma)\1n\in\1n L^2_{\mathbb{F}}(0,T;\mathbb{R}^{2d+l})\1n\times\1n L^2_{\mathbb{F}}( \triangle_\ast[0,T];\mathbb{R}^{2d+l})$,
 we introduce the following equation which is useful in the sequel:
\vskip-5mm
\bel{YZ1}\left\{\2n\ba{ll}
\ns\ds \dbY(t,s)\1n=\1n\dbY(t,t)\1n+\2n\int_s^t\2n\gamma(t,r)dr\1n-\2n\int_s^t\2n \dbZ(t,r)d W(r),\;\;(t,s)\1n\in\1n\triangle_\ast[0,T],\\[-0.5mm]
\ns\ds \dbY(t,t)\1n=\1n\xi(t)\1n+\2n\int_t^T\2n\big[\dbA(r,t)^{\top}\dbY(r,t)\1n+\1n\dbC(r,t)^{\top}\dbZ(r,t)\big] dr,\;t\1n\in\1n[0,T].
\ea\right.
\ee
\vskip-1.5mm
\no Then,  applying the result of \cite[Theorem 3.3]{Hamaguchi-Wang-I}, we have the following result.

\bl\label{duality}

Let $(\rm H1)$-$(\rm H2)$ hold and $(\Theta_1,\Theta_2,\Theta_3)\1n\in\1n L^\infty(0,T;\mathbb{R}^{l\times d})\1n\times\1n L^2([0,T]^2;\mathbb{R}^{l\times d})\1n\times\1n
 L^\infty(0,T;$ $\mathbb{R}^{l\times d})$ be fixed.
For any $(\xi,\gamma)\1n\in\1n L^2_{\mathbb{F}}(0,T;\mathbb{R}^{2d+l})\1n\times\1n L^2_{\mathbb{F}}( \triangle_\ast[0,T];\mathbb{R}^{2d+l})$ and
$(\t,\varphi_1,\varphi_2)\1n\in\1n\widetilde{\cI}$, (\ref{YZ1}) admits a unique adapted solution
$(\dbY, \dbZ)\1n\in\1n L^{2}_{\mathbb{F},c}(\triangle_\ast[0,T];\mathbb{R}^{2d+l})\1n\times\1n L^2_{\mathbb{F}}( \triangle_\ast[0,T];\mathbb{R}^{2d+l})$.
Moreover,  the following duality principle holds:
\vskip-6mm
$$\ba{ll}
\ns\ds \dbE\1n \int_\t^T\2n\Big(\langle\xi(t),\dbX(t) \rangle\1n+\2n\int_t^T\2n \langle\gamma(s,t),\cX(s,t) \rangle ds\Big)dt\\[-0.5mm]
\ns\ds =\dbE\1n \int_\t^T\2n\Big(\langle\dbY(t,\t),\Phi^{\f}(t) \rangle\1n
+\2n\int_t^T\2n \big[\langle\dbY(s,t),\Phi^b(s,t) \rangle\1n+\1n\langle\dbZ(s,t),\Phi^\si(s,t) \rangle\big] ds\Big)dt,
\ea$$
\vskip-1.5mm
\no where $\dbX$ solves (\ref{E7}) and
$\cX$ is an auxiliary function of $\dbX$ satisfying the following equation:
\vskip-6.5mm
$$\ba{ll}
\ns\ds \cX(t,r)\1n=\1n\Phi^{\f}(t)\2n+\2n\int_\t^r\3n\big( \dbA(t,s)\dbX(s)\1n+\1n\Phi^{b}(t,s)\big)ds
\1n+\2n\int_\t^r\3n \big(\dbC(t,s)\dbX(s)\1n+\1n\Phi^{\si}(t,s)\big)dW(s),(t,r)\2n\in \2n\triangle_{*}[\t,T].
\ea
$$
\vskip-2.5mm

\el


 We now state the representation result for the extended cost functional.

\bl\label{the3.2}
Let $(\rm H1)$--$(\rm H2)$ hold and $(\Theta_1,\Theta_2,\Theta_3,v)\1n\in\1n\mathcal{S}[0,T]$ be fixed.
Let $(P_1,P_2)\1n\in\1n L^\infty(0,T;$ $\mathbb{S}^{2d+l})\1n\times \1nL^2([0,T]^2;\mathbb{R}^{(2d+l)\times(2d+l)})$
 be the solution  of the system (\ref{P12}) and $(\dbY, \dbZ)\1n\in\1n L^{2}_{\mathbb{F},c}(\triangle_\ast[0,T];$ $\mathbb{R}^{2d+l})\1n\times\1n L^2_{\mathbb{F}}( \triangle_\ast[0,T];\mathbb{R}^{2d+l})$ be the adapted solution of equation (\ref{YZ1}) with  $(\xi,\gamma)\1n=\1n(\dbG_3, \dbG_4)$,
%
%
%
%
%
%
%
%
\no where  $(\dbG_3, \dbG_4)\1n\in\1n L^2_{\mathbb{F}}(0,T;\mathbb{R}^{2d+l})\1n\times\1n L^2_{\mathbb{F}}( \triangle_\ast[0,T];\mathbb{R}^{2d+l})$ is defined by
%
\vskip-5mm
\bel{G34}\ba{ll}
\ns\ds \dbG_3(r)\1n\=\1n\left(\2n\begin{array}{ccccc}
q(r)+ C(T,r)^{\top}G\si(T,r)  \\
\rho(r)+ D(T,r)^{\top}G\si(T,r)   \\
  Gb(T,r)
\end{array}
\2n\right)\2n+\2n\int_r^T\2n\dbC(s,r)^{\top}\2n P_1(s)\Phi^\si(s,r)ds\\[-0.5mm]
\ns\ds \qq\qq+\int_r^T\int_r^T\2n\dbC(\th',r)^{\top}\1nP_2(\theta',s)\Phi^\si(s,r)d\th'ds, \;\;r\2n\in\2n[0,T],\\[-0.5mm]
\ns\ds \dbG_4(\th, r)\1n\=\1nP_1(\th)\Phi^b(\th,r)\2n+\2n \int_r^T\3nP_2(\th, \theta')\Phi^b(\th',r)d\th',\;(\th, r)\1n\in\1n \triangle_\ast[0,T].
\ea
\ee
\vskip-2mm
\no Then for any $(\t,\1n\varphi_1,\1n\varphi_2)\1n\in\1n\widetilde{\cI}$ and  $\tilde{v}\1n\in\1n \mathcal{U}[\t,T]$,
 $\widetilde{\cJ}(\t,\1n\varphi_1,\1n\varphi_2;\1nu_1^{\Theta,v+\tilde{v}})$ admits the following representation:
\vskip-7mm
$$\ba{ll}
\ns\ds  \widetilde{\cJ}(\t,\varphi_1,\varphi_2;u_1^{\Theta,v+\tilde{v}}) \2n =\2n\varphi_2^{\top}G\varphi_2\1n
+\1n\dbE\1n\int_\t^T\2n\bigg\{ \wt\Phi^{\f}(t)^{\top} \Big(P_1(t)\wt\Phi^{\f}(t)\1n+ \3n \int_\t^T\2nP_2(s,t)\wt\Phi^{\f}(s)ds\Big)\\
\ns\ds\qq+\si(T,t)^{\top} G\si(T,t)\1n+ \2n \int_t^T\2n \Phi^{\si}(s,t)^{\top} \Big(P_1(s)\Phi^{\si}(s,t)\1n+ \2n \int_t^T\2nP_2(\th,s)\Phi^{\si}(\th,t)d \th\Big)ds\\
\ns\ds\qq+2\wt\Phi^{\f}(t)^{\top}\dbY(t,\t)\1n+ 2\1n \int_t^T\2n \Big( \Phi^{b}(s,t)^{\top}\dbY(s,t)\1n+ \1n \Phi^{\si}(s,t)^{\top}\dbZ(s,t) \Big)ds\bigg\}dt,
\ea
$$
\vskip-2mm
\no where $ \wt\Phi^{\f}(t)$ is defined by
\vskip-6mm
$$\ba{ll}
\wt\Phi^{\f}(t)\=\Phi^{\f}(t)+\big(0\enspace \tilde{v}(t)^{\top}\enspace 0\big)^{\top}, \q t\in[\t,T].
\ea$$
\vskip-6mm
\el

\proof
\vskip-4mm
For any $(\t,\varphi_1,\f_2)\1n\in\1n\widetilde{\cI}$ and $(\Theta_1,\Theta_2,\Theta_3,v+\tilde{v})\1n\in\1n\mathcal{S}[0,T]$,
let $(X_1,\sX_1,\sX_2)$ be the extended causal feedback solution corresponding to
$(\Theta_1,\Theta_2,\Theta_3,v+\tilde{v})$ (recall (\ref{E0601})).
Then, by Lemma \ref{lem2}, $\dbX$ satisfies the SVIE (\ref{E7}) replacing $\Phi^\f$ with $ \wt\Phi^\f$.
Correspondingly, so does its auxiliary equation $\cX$.
We begin to treat the first term $\dbE\langle G\sX_2(T),\sX_2(T)\rangle$ in $\wt\cJ(\t,\varphi_1,\f_2;u_1^{\Theta,v+\tilde{v}})$.
Applying the It\^{o}'s formula to
 $t\1n\mapsto \1n\sX_2(t)^{\top}G\sX_2(t)$ on $[\t, T]$ yields
\vskip-5mm
$$\ba{ll}
\ns\ds \dbE[\sX_2(T)^{\top}G\sX_2(T)]\1n=\1n\varphi_2^{\top}G\varphi_2\1n+\1n\dbE\1n\int_\t^T\2n\big[
\dbX(t)^{\top}\mathbb{G}_1(t)\dbX(t)\1n+\1n\si(T,t)^{\top}G\si(T,t)\\
\ns\ds \qq\q+2\big(\si(T,t)\1n^{\top}\1nGC(T,t)\enspace \si(T,t)\1n^{\top}\1nGD(T,t)\enspace b(T,t)\1n^{\top}\1nG\big)\dbX(t)\big]dt,
\ea
$$
\vskip-3mm
\no where
\vskip-6mm
$$\ba{ll}
\dbG_1(t)\=\left(\begin{array}{ccccc}
  C(T,t)^{\top}GC(T,t)  &   C(T,t)^{\top}GD(T,t)   & A(T,t)^{\top}G  \\
 D(T,t)^{\top}GC(T,t)   & D(T,t)^{\top}GD(T,t) & B(T,t)^{\top}G \\
  GA(T,t)  & GB(T,t)  & 0
\end{array}
\right).
\ea
$$
\vskip-2mm
\no Thus,  the extended cost functional (\ref{E1.9}) can be rewritten as
\vskip-6mm
\bel{Lifting-up}\ba{ll}
\ns\ds \2n \wt\cJ(\t,\1n\varphi_1,\1n\f_2;u_1^{\Theta,v+\tilde{v}})\2n=\2n\varphi_2\1n^{\top}G\varphi_2\2n+\2n\dbE\2n\int_\t^T\3n\Big[\mathbb{X}(t)\1n^{\top}\1n\mathbb{Q}(t)\mathbb{X}(t)dt
\2n+\2n\si(T,\1nt)\1n^{\top}\1nG\si(T,\1nt)\2n+2\dbX(t)\1n^{\top}\1n\dbG_2(t)\Big]dt,
\ea
\ee
\vskip-3mm
\no where $\dbQ$ is defined by (\ref{E3.1}) and
%
%
$\dbG_2(t)\=\bigg(\begin{array}{ccccc}
 q(t)+ C(T,t)^{\top}G\si(T,t)  \\[-0.6mm]
\rho(t)+ D(T,t)^{\top}G\si(T,t)   \\[-0.6mm]
  Gb(T,t)
\end{array}
\bigg).$
\vskip-1.5mm

In the following, we turn to calculating the term $\mathbb{X}(t)^{\top}\mathbb{Q}(t)\mathbb{X}(t)$.
By applying It\^{o}'s formula to $\th\mapsto\cX(r,\th)^{\top}P_1(r)\cX(r,\th)$ on $[\t,r]$ and
integrating it over $[\t,T]$, one has
\vskip-5mm
\bel{FG1}\ba{ll}
\ns\ds \dbE\1n\int^T_{\t}\1n\big[\dbX(r)^{\top}P_1(r)\dbX(r)\1n-\1n\wt\Phi^\f(r)^{\top}P_1(r)\wt\Phi^\f(r)\big]dr\\[-1mm]
\ns\ds=\dbE\1n\int^T_{\t}\1n\big[\cX(r,r)^{\top}P_1(r)\cX(r,r)\1n-\1n\wt\Phi^\f(r)^{\top}P_1(r)\wt\Phi^\f(r)\big]dr \\[-1mm]
\ns\ds=\dbE\int^T_{\t}\1n\int^T_{r}\1n\Big\{\big[2\cX(\theta,r)^{\top}\1nP_1(\th)\dbA(\theta,r)\2n+\1n\dbX(r)^{\top}\1n\dbC(\theta,r)^{\top}P_1(\th)\dbC(\theta,r)\big]\dbX(r) \\
\ns\ds \q+ \big[2\cX(\theta,r)\1n^{\top}\1nP_1(\th)\Phi^{b}(\theta,r)\2n+\1n\big(\Phi^{\si}(\theta,r)\1n^{\top}
\2n+\1n\dbX(r)\1n^{\top}\1n\dbC(\theta,r)\1n^{\top}\big)P_1(\th)\Phi^{\si}(\theta,r)\big]\Big\}d\theta dr.
\ea
\ee
\vskip-3mm
\no We see that the new term of the form $\cX(\theta,r)^{\top}(\cdots)\dbX(r)$ appears on the right hand. Next,
to handle this new term, for any $(\th,r)\in \triangle_{*}[\t,T]$,
using the It\^{o}'s formula to $\th'\mapsto\cX(\th,\th')^{\top}P_2(\th, r)^{\top}$
$\cX(r,\th')$ on $[\t,r]$, we have
\vskip-6mm
$$\ba{ll}
\ns\ds \dbE[\cX(\theta,r)^{\top}\1nP_2(\theta,r)^{\top}\1n\dbX(r)\1n-\1n\wt\Phi^\f(\theta)^{\top}P_2(\theta,r)^{\top}\wt\Phi^\f(r)]\\
\ns\ds \2n=\2n\dbE\2n\int^r_{\t}\3n\Big\{\frac{1}{2}\Big[\Phi^{\si}\1n(\th,\theta')\1n^{\top}\1nP_2(\theta,r)\1n^{\top}\1n\big(\dbC(r,\th')\dbX(\theta')\2n
+\2n\Phi^{\si}\1n(r,\theta')\big)\2n+\2n\dbX(\theta')\1n^{\top}\1n\dbC(r,\th')\1n^{\top}\1nP_2(\theta,r)\1n^{\top}\1n\Phi^{\si}\1n(\th,\theta')\\
\ns\ds\qq+\Phi^{\si}(r,\theta')^{\top}\1nP_2(\theta,r)\1n^{\top}\1n\big(\dbC(\th,\th')\dbX(\theta')\2n+\2n\Phi^{\si}(\th,\theta')\big)
+\1n\dbX(\theta')^{\top}\1n\dbC(\th,\th')\1n^{\top}\1nP_2(\theta,r)\1n^{\top}\1n\Phi^{\si}(r,\theta')\\
\ns\ds\qq+\dbX(\theta')^{\top}\1nG_1(\theta,\1nr,\1n\theta')\dbX(\theta')\Big]+\1n\big(\dbX(\theta')^{\top}\1n \dbA(\theta,\th')^{\top}+\Phi^{b}(\th,\theta')^{\top}\big)P_2(\theta,r)^{\top}\cX(r,\theta')\\
\ns\ds\qq+\cX(\theta,\theta')^{\top}\1nP_2(\theta,r)^{\top}\1n\big(\dbA(r,\th')\dbX(\theta')+\Phi^{b}(r,\theta')\big)d\theta',\\
\ea
$$
\vskip-3mm
\no where
\vskip-7mm
$$\ba{ll}
\ns\ds G_1(\theta,r,\theta')\1n\=\1n\dbC(\theta,\theta')^{\top}P_2(\theta,r)^{\top}\dbC(r,\theta')
\1n+\1n\dbC(r,\theta')^{\top}P_2(\theta,r)\dbC(\theta,\theta'),\;\t\1n\les\1n \th'\1n\les\1n r\1n\les\1n \th\1n\les\1n T.
\ea
$$
\vskip-1.5mm
\no  Combining Fubini theorem with $P_2(\th,r)^{\top}=P_2(r,\th), (\th,r)\in [\t,T]^2$, we have
\vskip-6mm
$$\ba{ll}
\ns\ds \dbE\int^T_{\t}\int^T_{r}[\cX(\theta,r)^{\top}P_2(\theta,r)^{\top}\dbX(r)-\wt\Phi^\f(\theta)^{\top}P_2(\theta,r)^{\top}\wt\Phi^\f(r)]d\theta dr\\
\ns\ds\2n=\2n\dbE\2n\int^T_{\t}\3n\int_{r}^{T}\3n\Big\{\1n\int_{\theta}^T\2n\frac{1}{2}\dbX(r)\1n^{\top}\1nG_1(\theta',\1n\theta,\1nr)\dbX(r)d\theta'
\2n+\3n\int_{r}^T\3n\Big[\1n\cX(\theta,r)\1n^{\top}\1nP_2(\theta,\th')\1n^{\top}\1n\big(\dbA(\th',r)\dbX(r)\2n+\1n\Phi^{b}(\theta',r)\big) \\
\ns\ds\qq+\dbX(r)^{\top}\1n\dbC(\th',r)\1n^{\top}\1n P_2(\theta',\th)^{\top}\1n\Phi^{\si}(\theta,r)
\1n+\1n \Phi^{\si}(\theta',r)^{\top}\1n P_2(\theta',\th)^{\top}\1n \Phi^{\si}(\theta,r)d\theta'\Big]\Big\}d\theta dr.
\ea
$$
\vskip-1.5mm
\no Thus, it follows from (\ref{FG1}) and the above equality that
\vskip-7mm
$$\ba{ll}
\ns\ds \wt\cJ(\t,\varphi_1,\f_2;u_1^{\Theta,v+\tilde{v}})\\
\ns\ds =\2n\varphi_2^{\top}\1nG\varphi_2\1n+\1n\dbE\2n\int_\t^T\2n\wt\Phi^{\f}(r)^{\top}P_1(r)\wt\Phi^{\f}(r)dr
\1n+2\dbE\2n\int_\t^T\2n\int_r^T\2n\wt\Phi^{\f}(\theta)^{\top}P_2(\theta,r)^{\top}\wt\Phi^{\f}(r)d\theta dr\\[-0.5mm]
\ns\ds \q\2n+\1n\dbE\2n\int_\t^T\3n\1n\Big\{\1n\si(T,r)\1n^{\top}\1nG\si(T,r) \2n+\3n\int_r^T\3n \Big[\1n\Phi^{\si}\1n(\th,r)\1n^{\top}\1nP_1(r)
\2n+\3n\int_r^T\3n\1n\Phi^{\si}\1n(\theta',r)\1n^{\top}\1nP_2(\theta',\th)\1n^{\top}\1nd\th'\1n\Big]\Phi^{\si}\1n(\th,r)d\theta\1n\Big\} dr\\[-0.5mm]
\ns\ds \q\2n+\dbE\2n \int_\t^T \3n\dbX(r)^{\top}\1n\Big[\dbQ(r)\1n+\2n\int_r^T\3n\Big(\1n\dbC(\th,r)^{\top}\1nP_1(\th)\dbC(\th,r)
\1n+\2n\int_{\theta}^T\3nG_1(\theta',\theta,r)d\theta'\Big)d\theta\1n-\1nP_1(r)\Big]\dbX(r)dr\\[-1mm]
\ns\ds \q\2n+2\dbE \1n\int_\t^T\2n\int_r^T\2n\Big[\cX(\theta,r)^{\top}\1n\Big(P_1(\th)\dbA(\th,r)\2n-\1nP_2(\theta,r)^{\top}
\1n+\2n\int_{r}^T\3nP_2(\theta,\th')\1n^{\top}\1n\dbA(\th',r)d\theta'\Big)\dbX(r)\Big]d\theta dr\\[-0.5mm]
\ns\ds \q\2n+2\dbE \1n\int_\t^T\2n \dbX(r)^{\top}\dbG_3(r)dr\1n+2\dbE \1n\int_\t^T\2n\int_r^T\2n \cX(\th,r)^{\top}\dbG_4(\th,r)d\th dr, 
\ea
$$
\vskip-2.5mm
\no where  $(\dbG_3, \dbG_4)$ is defined by (\ref{G34}).
Finally, the duality principle  (Lemma \ref{duality}) and the system (\ref{P12}) yield the representation of
 $\wt\cJ(\t,\1n\varphi_1,\1n\f_2;$ $\1nu_1^{\Theta,v+\tilde{v}})$. The proof is complete.
\endpf

\begin{remark}\label{Distinction-H-W}
As is shown in Remark \ref{Remark-difference-strategy}, both the causal feedback strategy and feedback outcome are different from \cite{Hamaguchi-Wang-I, Hamaguchi-Wang-II}. Therefore, the methodology in \cite[Theorem 4.1]{Hamaguchi-Wang-II} may not work well in our framework which prompts us to seek for new ideas. Thanks to the first equality in (\ref{E6}), one can transform the cost functional into a pure quadratic form without control by some lifting up arguments, see (\ref{Lifting-up}). Considering the fact that $\dbX$ satisfies a linear SVIE (\ref{E7}), it then becomes natural for us to borrow the arguments in \cite{Wang-Yong-2023-SICON} and come up with the above (\ref{P12}), as well as the procedures involved with It\^{o}'s formula and Fubini theorem in the previous proof.
\end{remark}

The following conclusion  can be obtained directly by Lemma \ref{the3.2}.
\begin{corollary}\label{cor3.3}
 Let $(\rm H1)$-$(\rm H2)$ hold and $(\Theta_1,\Theta_2,\Theta_3,v)\1n\in\1n\mathcal{S}[0,T]$ be fixed.
Let $(P_1,P_2)\1n\in\1n L^\infty(0,T;\mathbb{S}^{2d+l})\1n\times \1nL^2([0,T]^2;\mathbb{R}^{(2d+l)\times(2d+l)})$ be the solution of the system (\ref{P12}).
Define $(\dbG_3, \dbG_4)$ by (\ref{G34}), and let
$(\dbY, \dbZ)\1n\in\1n L^{2}_{\mathbb{F},c}(\triangle_\ast[0,T];\mathbb{R}^{2d+l})\1n\times\1n L^2_{\mathbb{F}}( \triangle_\ast[0,T];\mathbb{R}^{2d+l})$ be the solution of equation (\ref{YZ1}) with  $(\xi,\gamma)\1n=\1n(\dbG_3, \dbG_4)$.
Then,  for any $(\t,\varphi_1,\varphi_2)\in\widetilde{\cI}$,
it holds that
\vskip-6mm
$$\ba{ll}
\ns\ds  \widetilde{\cJ}(\t,\varphi_1,\varphi_2;u_1^{\Theta,v}) \2n =\2n\varphi_2^{\top}G\varphi_2\1n
+\1n\dbE\1n\int_\t^T\2n\bigg\{\Phi^{\f}(t)^{\top} \Big(P_1(t)\Phi^{\f}(t)\1n+ \3n \int_\t^T\2nP_2(s,t)\Phi^{\f}(s)ds\Big)\\
\ns\ds\qq+\si(T,t)^{\top} G\si(T,t)\1n+ \2n \int_t^T\2n \Phi^{\si}(s,t)^{\top} \Big(P_1(s)\Phi^{\si}(s,t)\1n+ \2n \int_t^T\2nP_2(\th,s)\Phi^{\si}(\th,t)d \th\Big)ds\\
\ns\ds\qq+2\Phi^{\f}(t)^{\top}\dbY(t,\t)\1n+ 2\1n \int_t^T\2n \Big( \Phi^{b}(s,t)^{\top}\dbY(s,t)\1n+ \1n \Phi^{\si}(s,t)^{\top}\dbZ(s,t) \Big)ds\bigg\}dt,
\ea
$$
\vskip-3mm
\no where $ \Phi^\f, \Phi^b$ and $\Phi^\si$ are defined by (\ref{E2.6}) and (\ref{b12}), respectively.

\end{corollary}

For the homogeneous version of Problem (LQ-SVIE), i.e., $b, \si, q, \rho\1n=\1n0$,
 we write the corresponding extended functional and extended value function by 
 $\widetilde{\cJ}^0(\t,\varphi_1,\varphi_2; \cd)$ and $\wt V^0(\t,\varphi_1,\varphi_2)$, respectively.
Then, Lemma \ref{the3.2} yields the following result.
\begin{corollary}\label{cor3.4}
 Let $(\rm H1)$-$(\rm H2)$ hold and $(\Theta_1,\Theta_2,\Theta_3,v)\1n\in\1n\mathcal{S}[0,T]$ be fixed.
Let $(P_1,P_2)\1n\in\1n L^\infty(0,T;\mathbb{S}^{2d+l})\1n\times \1nL^2([0,T]^2;\mathbb{R}^{(2d+l)\times(2d+l)})$ be the solution of the system (\ref{P12}).
Then,  for any $(\t,\varphi_1,\varphi_2)\in\widetilde{\cI}$, $\tilde{v}\1n\in\1n \mathcal{U}[\t,T]$, it holds that
\vskip-7mm
$$\ba{ll}
\ns\ds  \widetilde{\cJ}^0(\t,\varphi_1,\varphi_2;u_1^{\Theta,v+\tilde{v}}) \2n =\2n\varphi_2^{\top}G\varphi_2\1n
+\1n\dbE\1n\int_\t^T\2n \wt\Phi^{\f}(t)^{\top} \Big(P_1(t)\wt\Phi^{\f}(t)\1n+ \3n \int_\t^T\2nP_2(s,t)\wt\Phi^{\f}(s)ds\Big)dt.
\ea
$$
\vskip-2mm
\end{corollary}

To present the main result of this paper, i.e. Theorem \ref{the3}, let us give a more preliminary result.

\bl\label{lem4.1}
(i) Let $f_1\1n\in\1n L^2(0,T;\1n\mathbb{R}^{l\times d})$ and
$f_2\1n\in\1n L^2(\Box_3[0,T];\1n \mathbb{R}^{l\times d})$ be given functions such that for a.e. $\1n(s_1,\1ns_2)\1n\in\1n[0,T]^2, t\1n\mapsto\1n f_2(s_1,\1ns_2,\1nt)$ is
continuous on $[0,\1ns_1\wedge s_2]$. Assume that for any $\t\in[0,T],x\in L^2(\t,T;\mathbb{R}^d)$, $v\in\mathcal{U}[\t,T]$, it holds that
\vskip-5.5mm
$$\ba{ll}
 \ns\ds\int^T_{\t}v(s)^\top f_1(s)x(s)ds+\int^T_{\t}\int^T_{\t} v(s_2)^\top f_2(s_1,s_2,\t)x(s_1)ds_1ds_2=0.
\ea$$
\vskip-2mm
\no Then $f_1(t)\2n=\2n0,$ a.e. $\1nt\1n\in\1n[0,T]$ and $f_2(s_1,\1ns_2,\1nt)\2n=\2n0,$ a.e. $\1n(s_1,\1ns_2)\1n\in\1n[0,T]^2, \forall t\1n\in\1n[0,\1ns_1\1n\wedge\1n s_2]$.

\no(ii) Let $f:\triangle_{*}[0,T]\rightarrow \mathbb{R}^{l}$ such that for a.e. $s\in[0,T],\; t\mapsto f(s,t)$ is
continuous on $[0,s]$. Assume that
\vskip-7.5mm
$$\ba{ll}
\ns\ds \int^T_{\t} v(s)^\top f(s,\t)ds=0
\ea
$$
\vskip-1.5mm
\no for any $\t\1n\in\1n[0,T),v\1n\in\1n\mathcal{U}[\t,T]$. Then
$ f(s,t)\1n=\1n0 $ for $\ \ae \1ns\1n\in\1n[0,T]$ and any $t\1n\in\1n[0,s]$.
\el
The proof of the above lemma is a proper adjustment to \cite[Lemma 4.4]{Hamaguchi-Wang-I}.
We omit the details for the sake of limitation space.

In the following, for the sake of the sequel, we define
\vskip-5mm
\bel{f123}
 \widehat{f}_{12}(t,s)\1n\=\1n\BR(t)^\dagger \BD(t)f(t,s)\1n
 +\2n\int^T_s\2n \BR(t)^\dagger \BB_1(r,t)I_{[t,T]}(r)f(r,s)dr\1n+\1n\BR(t)^\dagger \BB_2(t)f(T,s)
\ee
\vskip-3mm
\no for $(t,s)\2n\in\2n\triangle_\ast[0,T]$ and $f\1n\=\1nA,B,C,D,b,\si$ (recall the notions $\BR, \BD,\BB_1$ and $\BB_2$ in Definition \ref{def2.11}).
 Consequently, denote by $\widehat{\Phi}^\f,\widehat{\dbA}, \widehat{\dbC}, \widehat{\Phi}^b$
and $\widehat{\Phi}^\si$ the corresponding notations to (\ref{E2.6})-(\ref{b12}) with $(\Th_1,\Th_2,\Th_3)$ replaced by $(\BR^\dagger \BD,\BR^\dagger \BB_1,\BR^\dagger \BB_2)$.

%
%
%
%
%
%
%
%
%
%
%
%
%

\ms

\begin{theorem}\label{the3}
Suppose $(\rm H1)$-$(\rm H2)$ hold. Then Problem (LQ-SVIE) admits a strongly optimal causal feedback strategy $(\widehat{\Theta}_1,\widehat{\Theta}_2,\widehat{\Theta}_3,\widehat{v})$
 if and only if   the following two conditions hold:
\begin{enumerate}[1)]
\item   The Riccati system (\ref{E4.1}) admits a regular solution $\sP\1n=\1n(\sP_1,\sP_2,\sP_3,\sP_4)\1n\in\1n \Upsilon[0,T]$.

\item  Let $(P_1,P_2)\1n\in\1n L^\infty(0,T;\mathbb{S}^{2d+l})\1n\times \1nL^2([0,T]^2;\mathbb{R}^{(2d+l)\times(2d+l)})$
be the solution of the system (\ref{P12}) with $(\Th_1,\Th_2,\Th_3)\1n=\1n(\BR^\dagger \BD,\BR^\dagger \BB_1,\BR^\dagger \BB_2)$.
Define $\widehat{G}_4\in L^2_{\mathbb{F}}( \triangle_\ast[0,T];\mathbb{R}^{d})$
\vskip-6mm
\bel{G4P12}\ba{ll}
\ns\ds \2n\widehat{G}_4(t, r)\1n\=\1n \big(0_{d\times d}\enspace 0_{d\times l} \enspace I_{d}\big)
 \Big(P_1(t)\widehat{\Phi}^b(t,r)\2n+\2n \int_r^T\3nP_2(t, \theta')\widehat{\Phi}^b(\th',r)d\th'\Big),\;(t, r)\1n\in\1n \triangle_\ast[0,T],
\ea
\ee
\vskip-3mm
\no and denote $(\h\dbG_3, \h\dbG_4)\1n\in\1n L^2_{\mathbb{F}}(0,T;\mathbb{R}^{2d+l})\1n\times\1n L^2_{\mathbb{F}}( \triangle_\ast[0,T];\mathbb{R}^{2d+l})$
%
%
\vskip-7mm
$$\ba{ll}
\ns\ds \widehat{\dbG}_3(t)\1n\=\1n\left(\2n\begin{array}{ccccc}
 q(t)\2n+\2n (C^{\top}\2n\triangleleft\sP\triangleright  \si)(t)  \\[-0.5mm]
\rho(t)\2n+\2n (D^{\top}\2n\triangleleft\sP\triangleright  \si)(t)  \\[-0.5mm]
  Gb(T,t)
\end{array}
\2n\right),  t\1n\in\1n[0,T],\;
\widehat{\dbG}_4(t, r)\1n\=\1n\left(\3n\begin{array}{ccccc}
 (\sP_{1,3,4}\triangleright  b)(t,r)  \\[-0.5mm]
0  \\[-0.5mm]
  \widehat{G}_4(t,r)
\end{array}
\3n\right), (t,r)\1n\in\1n\triangle_*[0,T].
\ea
$$
\vskip-3mm
\no  Let $(\dbY, \dbZ)\1n\in\1n L^{2}_{\mathbb{F},c}(\triangle_\ast[0,T];\mathbb{R}^{2d+l})\1n\times\1n L^2_{\mathbb{F}}( \triangle_\ast[0,T];\mathbb{R}^{2d+l})$
be the adapted solution of the following equation:
\vskip-6mm
\bel{YZ3}\left\{\2n\ba{ll}
\ns\ds \dbY(t,s)\1n=\1n\dbY(t,t)\1n+\2n\int_s^t\2n\widehat{\dbG}_4(t,r)dr\1n-\2n\int_s^t\2n \dbZ(t,r)d W(r),\;\;(t,s)\1n\in\1n\triangle_\ast[0,T],\\
\ns\ds \dbY(t,t)\1n=\1n\widehat{\dbG}_3(t)\1n+\2n\int_t^T\2n\big[\widehat{\dbA}(r,t)^{\top}\dbY(r,t)\1n+\1n\widehat{\dbC}(r,t)^{\top}\dbZ(r,t)\big] dr,\;t\1n\in\1n[0,T].
\ea\right.
\ee
\vskip-1.5mm
\no  Then the process $\h\sN\1n\in\1n L^2_{\mathbb{F}}( 0,T;\dbR^l)$ defined by (recall $\h B_{12}$ and $\h D_{12}$ in (\ref{f123}))
\vskip-6.5mm
\bel{sN}\ba{ll}
\ns\ds\h \sN(t)\1n=\1n\rho(t)\1n+\1n (D^{\top}\2n\triangleleft\sP\triangleright  \si)(t)\1n+\2n\int_t^T\2n\big(B(s,t)^{\top}\2n\enspace \widehat{B}_{12}(s,t)^{\top}\2n \enspace B(T,t)^{\top}\big)\dbY(s,t)ds\\[-1mm]
\ns\ds\qq\qq+\int_t^T\2n\big(D(s,t)^{\top}\2n\enspace \widehat{D}_{12}(s,t)^{\top}\2n \enspace D(T,t)^{\top}\big)\dbZ(s,t)d s,
\ea
\ee
\vskip-2mm
\no satisfies $\h \sN(t)\1n\subset\1n\sR(\BR(t))$ for a.e. $\1nt\1n\in\1n[0,T]$.

\end{enumerate}

\no In this case,  it holds that
\vskip-5mm
\bel{G4}\ba{ll}
\ns\ds  \int_r^T\2n\widehat{G}_4(t,r)dt =(\sP_{2,3}\triangleright  b)(r)-Gb(T,r),\;\;r\in[0,T),
\ea
\ee
\vskip-1.5mm
\no and $(\widehat{\Theta}_1,\widehat{\Theta}_2,\widehat{\Theta}_3,\widehat{v})$ admits the following representation:
\vskip-6mm
 \bel{E2403207}\ba{ll}
\ns\ds  \BD(t)\1n+\1n\BR(t)\widehat{\Theta}_1(t)\1n=\1n0,\;\;  \BB_2(t)\1n+\1n\BR(t)\widehat{\Theta}_3(t)\1n=\1n0,\;\;
    \h\sN(t)\1n+\1n\BR(t)\widehat{v}(t)\1n=\1n0, \;\; t\in[0,T],\\
\ns\ds   \BB_1(r,t)I_{[0,r]}(t)\1n+\1n\BR(t)\widehat{\Theta}_2(r,t)\1n=\1n 0,\;\;(r,t)\1n\in\1n[0,T]^2.
\ea
\ee
\vskip-2.5mm
\no   Further, for each $(\t,\f_1,\f_2)\in\wt\cI$, the extended value function is given by  
\vskip-6.5mm
\bel{E24032601}\ba{ll}
\ns\ds \wt V(\t, \f_1,\f_2)\1n=\1n\varphi_2^{\top}\1n\sP_2(\t)\varphi_2\1n+\1n\dbE\1n\int_\t^T\2n\Big\{\varphi_1(s)^{\top}\1n\sP_1(s)\varphi_1(s)
\1n+\1n2\varphi_2^{\top}\1n\sP_3(s,\t)\varphi_1(s)\\[-1.8mm]
\ns\ds\qq+\2n\int_\t^T\3n\1n\varphi_1(s)^{\top}\1n\sP_4(s,r,\t)\varphi_1(r)dr\1n+\1n\widehat{v}(s)^{\top}\BR(s)\widehat{v}(s)
\1n+\1n (\si^{\top}\2n\triangleleft\sP\triangleright  \si)(s)\\[-1mm]
\ns\ds\qq+ 2\widehat{\Phi}^{\f}(s)^{\top}\1n\dbY(s,\t)\1n+\2n \int_s^T\2n 2 \big[\langle\dbY(r,s),\widehat{\Phi}^b(r,s) \rangle\1n+\1n\langle\dbZ(r,s),\widehat{\Phi}^\si(r,s) \rangle\big] dr\Big\}ds.
\ea
\ee
\vskip-3mm

\end{theorem}
\proof
\vskip-3.5mm
\no\rm \textbf{The necessity}: Suppose $(\widehat{\Theta}_1,\widehat{\Theta}_2,\widehat{\Theta}_3, \widehat{v})\in \mathcal{S}[0,T]$ is a
strongly optimal causal feedback strategy.
Let $(P_1,P_2)\1n\in\1n L^\infty(\t,T;\mathbb{S}^{2d+l})\1n\times \1nL^2([\t,T]^2;\mathbb{R}^{(2d+l)\times(2d+l)})$
 be the solution  of the system (\ref{P12}) with $(\Th_1,\Th_2,\Th_3)\1n=\1n(\h\Th_1,\h\Th_2,\h\Th_3)$.
Recall Lemma \ref{the3.2} and  Corollary \ref{cor3.3},
 let $(\dbY, \dbZ)\1n\in\1n L^{2}_{\mathbb{F},c}(\triangle_\ast[0,T];\mathbb{R}^{2d+l})\1n\times\1n L^2_{\mathbb{F}}( \triangle_\ast[0,T];\mathbb{R}^{2d+l})$ be the adapted solution of equation (\ref{YZ1}) with  $(\xi,\gamma)\1n=\1n(\dbG_3, \dbG_4)$,
for any $(\t,\varphi_1,\varphi_2)\in\widetilde{\cI}$ and $v\in\mathcal{U}[\t,T]$, we can obtain that
\vskip-7mm
$$\ba{ll}
\ns\ds \widetilde{\cJ}(\t,\varphi_1,\varphi_2;u_1^{\widehat{\Theta},\widehat{v}+v})
- \widetilde{\cJ}(\t,\varphi_1,\varphi_2;u_1^{\widehat{\Theta},\widehat{v}})\\
  \ns\ds =\1n\dbE\2n\int_\t^T \2nv(t)^\top\1n\Big[\1nP_{22}^{(1)}(t) v(t)\1n+\2n\int_\t^T\2nP^{(2)}_{22}(s,t)v(s)ds
  \1n+\2n2 P_{22}^{(1)}(t) \widehat{v}(t)\1n+\2n2\int_\t^T\2nP^{(2)}_{22}(s,t)\widehat{v}(s)ds\\[-1.2mm]
    \ns\ds\qq+ 2 \sN(t) +  2 \sM_1(t)\varphi_1(t)+2\sM_2(t,\t)\varphi_2+2\int_\t^T\sM_3(s,t,\t)\varphi_1(s)ds\Big]dt,
\ea
$$
\vskip-3mm
\no where $P^{(1)}_{ij}$, $P^{(2)}_{ij}$, $i,j\1n=\1n\{1,2,3\}$ are components of the solution $(P_1,P_2)$ and
\vskip-6mm
$$\ba{ll}
\ns\ds  \sN(t)\=\big(0_{l\times d}\enspace I_{l} \enspace 0_{l\times d}\big)\dbY(t,\t),\;\;\;
\;\;\;\;\sM_1(t)\=P_{21}^{(1)}(t) + P_{22}^{(1)}(t) \widehat{\Theta}_1(t),\\[-1.5mm]
\ns\ds
\sM_2(t,\t)\=P_{23}^{(1)}(t)+ P_{22}^{(1)}(t)  \widehat{\Theta}_3(t)+\int_\t^T\big(P^{(2)}_{23}(s,t)+P^{(2)}_{22}(s,t)\widehat{\Theta}_3(s)\big)ds,\\[-1.5mm]
\ns\ds
\sM_3(s,t,\t)\1n\=\1n
P^{(2)}_{21}(s,t)\1n+\1n P^{(2)}_{22}(s,t)\widehat{\Theta}_1(s)
\1n+\1nP_{22}^{(1)}(t) \widehat{\Theta}_2(s,t)\1n+\2n\int_\t^T\3n P^{(2)}_{22}(s_1,t)\widehat{\Theta}_2(s_1,s)ds_1.
\ea
$$
\vskip-3mm
\no By the strongly optimality of $(\widehat{\Theta}_1,\widehat{\Theta}_2,\widehat{\Theta}_3, \widehat{v})$,
 for any $(\t,\varphi_1,\varphi_2)\1n\in\1n\widetilde{\cI}$,
we can obtain that
\vskip-6mm
$$\ba{ll}
\ns\ds
 \widetilde{\cJ}(\t,\varphi_1,\varphi_2;u_1^{\widehat{\Theta},\widehat{v}+v})
- \widetilde{\cJ}(\t,\varphi_1,\varphi_2;u_1^{\widehat{\Theta},\widehat{v}})\ges 0,\;\;\forall v\1n\in\1n\mathcal{U}[\t,T],
\ea
$$
\vskip-2.5mm
\no which implies that
\vskip-5.5mm
\bel{E4.3}
\dbE\int_\t^T v(t)^\top P_{22}^{(1)}(t) v(t) dt +\dbE\int_\t^T \int_\t^Tv(t)^\top P^{(2)}_{22}(s,t)v(s)dsdt\ges 0,
\ee
\vskip-2mm
\no and
\vskip-8mm
\begin{eqnarray}
\dbE\2n\int_\t^T \2nv(t)^\top\1n\Big[\1n P_{22}^{(1)}(t) \widehat{v}(t)\1n+\2n\int_\t^T\2nP^{(2)}_{22}(s,t)\widehat{v}(s)ds+\sM_1(t)\varphi_1(t)\notag\\[-1.8mm]
+ \sN(t)+\sM_2(t,\t)\varphi_2+\int_\t^T\sM_3(s,t,\t)\varphi_1(s)ds\Big]dt=0.\label{E4.4}
\end{eqnarray}
\vskip-3mm

Now we take $v(\cd)\2n\=\2nvI_{[\bar{\t},\bar{\t}+\e]}(\cd), \bar{\t}\2n\in\2n[\tau,T),$ with $\bar{\t}\1n+\1n\e\1n<\1nT$ and any $v\2n\in\2n\dbR^l$.
 Then (\ref{E4.3}) becomes
\vskip-5.5mm
 $$\ba{ll}
   \ns\ds\int_{\bar{\t}}^{\bar{\t}+\e} v ^\top P_{22}^{(1)}(t) v  dt+\int_{\bar{\t}}^{\bar{\t}+\e} \int_{\bar{\t}}^{\bar{\t}+\e}v^\top P^{(2)}_{22}(s,t)vdsdt\ges 0.
\ea
$$%
\vskip-2mm
\no By Lebesgue differentiation theorem,
\vskip-4mm
$$\ba{ll}
\ns\ds  \lim_{\e\to 0}\frac{1}{\varepsilon}\int_{\bar{\t}}^{\bar{\t}+\e} v ^\top P_{22}^{(1)}(t) v  dt=v^\top P_{22}^{(1)}(\bar{\t}) v,\;\;\;$a.e. $\bar{\t}\in[\t,T).
\ea
$$%
\vskip-1.5mm
\no On the other hand, by H\"{o}lder's inequality,
\vskip-6.5mm
$$\ba{ll}
\ns\ds  \lim_{\e\to 0}\frac{1}{\varepsilon}\2n\int_{\bar{\t}}^{\bar{\t}+\e}\2n \int_{\bar{\t}}^{\bar{\t}+\e}\2nv^\top P^{(2)}_{22}(s,t)  v  dsdt
\les\lim_{\e\to 0}\1n\Big(\2n\int_{\bar{\t}}^{\bar{\t}+\e}\2n\int_{\bar{\t}}^{\bar{\t}+\e}\2n |v^\top P^{(2)}_{22}(s,t)  v |^2 dsdt\1n\Big)^{\frac{1}{2}}=0.
\ea
$$%
\vskip-2.5mm
\no Therefore, the arbitrariness of $v$ leads to
\vskip-5.5mm
\bel{E4.14}
 P_{22}^{(1)}(\bar{\t})\ges0,\ \ae \bar{\t}\in[\t,T).
\ee%
\vskip-2.5mm
 By the arbitrariness of $\varphi_1$ and $\f_2$, let us take $\varphi_1=\f_2=0$. Then,
 for any $\t\in[0,T)$ and $ v\1n\in\1n\mathcal{U}[\t,T]$, (\ref{E4.4}) becomes
\vskip-5.5mm
\bel{E24032301}
\dbE\int_\t^T \1n v(t)^\top\1n\big(\sN(t)\1n+\1n P_{22}^{(1)}(t) \widehat{v}(t)\big) dt \1n+\1n \dbE\int_\t^T \int_\t^T\1nv(t)^\top P^{(2)}_{22}(s,t)\widehat{v}(s)dsdt=0.
\ee
\vskip-3mm
\no Thus, it follows from (\ref{E4.4}) and (\ref{E24032301}) that for any
$(\t,\varphi_1,\varphi_2)\1n\in\1n\widetilde{\cI}$, $v\1n\in\1n\mathcal{U}[\t,T]$,
\vskip-5mm
\bel{E24032302}
\dbE\int_\t^Tv(t)^\top\Big[\sM_1(t)\varphi_1(t)dt+\sM_2(t,\t)\varphi_2+\int_\t^T\sM_3(s,t,\t)\varphi_1(s)ds\Big]dt=0.
\ee
\vskip-2mm
\no Again by the  arbitrariness of $\varphi_1$, let us take $\varphi_1=0$. Then (\ref{E24032302}) becomes
\vskip-5.5mm
$$\ba{ll}
\ns\ds\dbE\int^T_\t  v(t)^\top\sM_2(t,\t)\varphi_2dt=0,  \ \ \t\1n\in\1n[0,T), v\1n\in\1n\mathcal{U}[\t,T], \varphi_2\1n\in\1n\mathbb{R}^d.
\ea$$
\vskip-2mm
\no For a.e. $t\in[0,T), \t\mapsto \sM_2(t,\t)$
 is absolutely continuous on $[0,t]$.
Therefore, it follows from Lemma \ref{lem4.1} that
\vskip-6.5mm
\bel{E4.6}
\sM_2(t,t_1)=0\;\;\ae t\in[0,T],\;\;\forall t_1\in[0,t].
\ee
\vskip-2.5mm
\no  Then, by (\ref{E24032302}) and (\ref{E4.6}), for any $\t\1n\in\1n[0,T),\varphi_1\1n\in\1n L^2(\t,T;\mathbb{R}^d), v\1n\in\1n\mathcal{U}[\t,T]$,
 we get
\vskip-5.5mm
$$\ba{ll}
 \ns\ds\dbE\int_\t^Tv(t)^\top\sM_1(t)\varphi_1(t)dt+\dbE\int_\t^T\int_\t^Tv(t)^\top\sM_3(s,t,\t)\varphi_1(s)dsdt=0.
\ea$$
\vskip-2.5mm
\no
For  a.e. $\1n(s,t)\1n\in\1n[0,T]^2,$ we have $\t\1n\mapsto \1n\sM_3(s,t,\t)$
 is absolutely continuous on $[0,s\wedge t]$ due to the absolutely continuity of Lebesgue integral. Thus,
by Lemma \ref{lem4.1},
\vskip-6mm
\bel{E4.7}
\sM_1(t)\1n=\1n0,\;\ae t\1n\in\1n[0,T],\;\;\sM_3(s,t,t_1)\1n=\1n0,\;\ae (s,t)\1n\in\1n[0,T]^2, \;\forall t_1\1n\in\1n[0,s\wedge t].
\ee
\vskip-2.5mm
Recall system (\ref{P12}), we can rewrite $P_{12}^{(2)}$, the components of $P_2$, as follows:
\vskip-7mm
$$\ba{ll}
\ns\ds P_{12}^{(2)}(r,t)\2n=\2nA(r,t)\1n^{\top}\2n\sM_1(t)\1n^{\top}\3n+\2nA(T,t)\1n^{\top}\2n\sM_2(r,t)\1n^{\top}
\3n+\3n\int_t^T\3n\1nA(s,t)\1n^{\top}\2n\sM_3(r,s,t)\1n^{\top}ds,(r,t)\1n\in\1n\triangle_\ast[\t,T].
\ea
$$
\vskip-2.5mm
\no  By  (\ref{E4.6}), (\ref{E4.7}) and the fact that $P_2(r,t)=P_2(t,r)^{\top}$, it is easily seen that
\vskip-5mm
\bel{E4.10}
 P_{12}^{(2)}(r,t)=P_{21}^{(2)}(t,r)^{\top}\2n=0,\;\;\; (r,t)\in\triangle_\ast[0,T].
\ee
\vskip-3mm
\no Similarly, it holds that
\vskip-7mm
\bel{E4.11}
P_{22}^{(2)}(r,t)=0,\;\;\;(r,t)\in[0,T]^2.
\ee%
\vskip-2mm
\no In addition, it follows directly from the structure of the system (\ref{P12}) that
\vskip-5mm
\bel{E4.12}
P_{31}^{(2)}(r,t)=P_{32}^{(2)}(r,t)=P_{33}^{(2)}(r,t)=0,\;\;\;(r,t)\in\triangle_\ast[0,T].
\ee
\vskip-2mm
\no Substituting (\ref{E4.10})--(\ref{E4.12}) into (\ref{E4.6}) and (\ref{E4.7}), respectively, we can see that
\vskip-5.5mm
\bel{E4.13}\left\{\ba{ll}
\ns\ds P_{21}^{(1)}(t)+ P_{22}^{(1)}(t)\widehat{\Theta}_1(t)=0,\;\;\;$a.e.$\; t\in[0,T],\\[-1mm]
\ns\ds P_{21}^{(2)}(r,t)I_{[0,r]}(t)+P_{22}^{(1)}(t)\widehat{\Theta}_2(r,t)=0,\;\;$a.e.$\;(r,t)\in [0,T]^2,\\[-1mm]
\ns\ds P_{23}^{(1)}(t) + P_{22}^{(1)}(t) \widehat{\Theta}_3(t)+\1n\int_{t}^T\1nP_{23}^{(2)}(s,t)ds=0,\;\;\;$a.e.$\; t\in[0,T].
\ea\right.
\ee
\vskip-2mm
\no By inserting (\ref{E4.10})--(\ref{E4.13}) into the system (\ref{P12}),
we can rewrite $P_1$ as
\vskip-5.5mm
\bel{E240303}\ba{ll}
\ns\ds P_1(t)\2n=\2n\left(\begin{array}{ccccc}
\2n Q(t)\1n+\1n(C^{\top}\2n\triangleleft\sP\triangleright  C)(t)\1n  &  \1n(C^{\top}\2n\triangleleft \sP\triangleright\1n D)(t) +S(t)^{\top}\1n \1n
  &  \1nA(T,t)^{\top} \1nG\2n \\
\2n (D^{\top}\2n\triangleleft \sP \triangleright  C)(t) +S(t)\1n  & \1n R(t)\1n+\1n(D^{\top}\2n\triangleleft \sP \triangleright \1n  D)(t)\1n
 &  \1n B(T,t)^{\top}\1n G\2n \\
 \2nGA(T,t) \1n  &\1n GB(T,t)\1n &\1n  0\2n
\end{array}
\right),
\ea
\ee
\vskip-2mm
\no where $\sP=(\sP_1,\sP_2,\sP_3,\sP_4)$ is given by
\vskip-6.5mm
$$\left\{\2n\ba{ll}
\ns\ds\sP_1(s)\=P_{11}^{(1)}(s)+\widehat{\Theta}_1(s)^{\top}P_{21}^{(1)}(s),\;\;s\in[0,T],\\[-1.5mm]
\ns\ds\sP_2(s)\=G-\int_s^T\widehat{\Theta}_3(r)^{\top}P_{22}^{(1)}(r)\widehat{\Theta}_3(r)dr,\;\;s\in[0,T],\\[-1.8mm]
\ns\ds\sP_3(s,t)\2n\=\2nP_{31}^{(1)}\1n(s)\2n+\2n\widehat{\Theta}_3(s)^{\top}P_{21}^{(1)}(s)
\1n+\2n\int_t^s\2n\widehat{\Theta}_3(\th)^{\top}P_{21}^{(2)}(s,\th)d\th\1n
+\3n\int_s^T\3nP_{31}^{(2)}(s,\th)d\th,(s,t)\1n\in\2n\triangle_{*}[0,T],\\
\ns\ds\sP_4(s,\th,t)\2n\=\2nP_{11}^{(2)}(\th,s)\2n+\2n\widehat{\Theta}_1(s)^{\top}P_{21}^{(2)}(\th,s)\2n
+\2nP_{12}^{(2)}(\th,s)\widehat{\Theta}_1(\th)\2n+\3n\int_t^{\th\wedge s} \3nP_{12}^{(2)}(r,s)P_{22}^{(1)}(r)^\dagger
 P_{21}^{(2)}(\theta,r)dr\\
\ns\ds\qq\qq\1n=\1n\sP_4(\th,s,t)^{\top},(s,\th,t)\1n\in\1n\Box_3[0,T].
\ea\right.
$$
\vskip-2.5mm

Next we represent the components of $P_2$ in terms of $(\sP_1,\sP_2,\sP_3,\sP_4)$.
By Fubini theorem and the structure of the system (\ref{P12}), we can get, for $(s,t)\1n\in \1n\triangle_{*}[0,T],$
\vskip-7mm
$$\ba{ll}
\ns\ds P_{11}^{(2)}(s,t)+P_{12}^{(2)}(s,t)\widehat{\Theta}_1(s)\\
\ns\ds=A(s,t)^{\top}\big[P_{11}^{(1)}(s)+\widehat{\Theta}_1(s)^{\top}P_{21}^{(1)}(s)\big]\\[-1mm]
\ns\ds\q\2n+A(T,t)^{\top}\Big[P_{31}^{(1)}(s)+\widehat{\Theta}_3(s)^{\top}P_{21}^{(1)}(s)
+\int_t^s\widehat{\Theta}_3(\th)^{\top}P_{21}^{(2)}(s,\th)d\th+\2n\int_s^T\3nP_{31}^{(2)}(s,\th)d\th\Big]\\[-1.5mm]
\ns\ds\q\2n+\3n\int_t^T\3n\2nA(r,t)^{\top}\1n\Big[
P_{11}^{(2)}\1n(s,r)\2n+\2n\widehat{\Theta}_1(r)^{\top}\2nP_{21}^{(2)}\1n(s,r)\2n+\2nP_{12}^{(2)}(s,r)\widehat{\Theta}_1(s)
\2n+\3n\int_t^{r\wedge s}\3n\widehat{\Theta}_2(r,\th)^{\top}\1nP_{21}^{(2)}\1n(s,\th)d\theta\Big]dr\\[-1.8mm]
\ns\ds=A(s,t)^{\top}\sP_1(s)+A(T,t)^{\top}\sP_3(s,t)+\int_t^T\1nA(r,t)^{\top}\sP_4(r,s,t)dr.
\ea
$$%
\vskip-2.5mm
\no Combining the above equality  with (\ref{E4.10}), we obtain that
\vskip-7.5mm
\bel{E240302}\ba{ll}
\ns\ds P_{11}^{(2)}\1n(s,t)\1n=\1n(A^{\top}\2n\triangleleft \1n\sP_{1,3,4})(s,t),
\; P_{11}^{(2)}\1n(t,s)
\1n=\1n(\sP_{1,3,4}\1n\triangleright\1n A)(s,t),\;(s,t)\1n\in \1n\triangle_{*}[0,T].
\ea
\ee
\vskip-2mm
\no Similarly, it also holds that
\vskip-6.5mm
\bel{E240304}\ba{ll}
\ns\ds P_{21}^{(2)}(s,t)\1n=\1n(B^{\top}\2n\triangleleft \1n\sP_{1,3,4})(s,t),\; P_{12}^{(2)}(t,s)
\1n=\1n(\sP_{1,3,4}\1n\triangleright\1n B)(s,t),\;(s,t)\1n\in\1n \triangle_{*}[0,T],\\
\ns\ds P_{13}^{(1)}(t)\1n+\2n\int_t^T\3nP_{13}^{(2)}(r,t)dr
\1n=\1n(A^{\top}\2n\triangleleft \1n\sP_{2,3})(t),
\; P_{31}^{(1)}(t)\1n+\2n\int_t^T\3nP_{31}^{(2)}(t,r)dr
\1n=\1n(\sP_{2,3}\1n \triangleright\1n A)(t),\\
\ns\ds P_{23}^{(1)}(t)\1n+\2n\int_t^T\3nP_{23}^{(2)}(r,t)dr\1n=\1n(B^{\top}\2n\triangleleft \1n\sP_{2,3})(t),
\;P_{32}^{(1)}(t)\1n+\2n\int_t^T\3nP_{32}^{(2)}(t,r)dr\1n=\1n(\sP_{2,3}\1n\triangleright\1n B)(t).
\ea
\ee
\vskip-2mm
\no Substituting the representations (\ref{E240303})--(\ref{E240304}) of the components of $(P_1,P_2)$ into (\ref{E4.13}),
we then obtain that  $(\widehat{\Theta}_1,\widehat{\Theta}_2,\widehat{\Theta}_3)$ satisfies  
\vskip-4mm
\bel{E24120801}\ba{ll}
\BD(t)\1n+\1n\BR(t)\widehat{\Theta}_1(t)\1n=\1n0,\;\;  \BB_2(t)\1n+\1n\BR(t)\widehat{\Theta}_3(t)\1n=\1n0, \;\;$a.e.$\;\; t\in[0,T],\\
\ns\ds   \BB_1(r,t)I_{[0,r]}(t)\1n+\1n\BR(t)\widehat{\Theta}_2(r,t)\1n=\1n 0,\;\;$a.e.$\;\;(r,t)\1n\in\1n[0,T]^2.
\ea\ee
\vskip-2mm
\no  Combining  the above equalities with (\ref{E240303})-(\ref{E240304}) we see that
$(\sP_1,\1n\sP_2,\1n\sP_3,\1n\sP_4)$ satisfies the system (\ref{E4.1}).
Furthermore, by slightly modifying \cite[Proposition A.1.5]{Sun-Yong-2020} and (\ref{E4.14})
it follows from that $\sP\2n=\2n(\sP_1,\1n\sP_2,\1n\sP_3,\1n\sP_4)\2n\in \2n\Upsilon[0,T]$ is a regular solution of the Riccati system (\ref{E4.1}).
Then,  by  (\ref{E4.12}), (\ref{E240304}), (\ref{E24120801}) and Fubini theorem, we arrive at (\ref{G4}) holds.

To determine $\widehat{v}$,  combining (\ref{E4.10})--(\ref{E240304}) yields that the equation (\ref{YZ1})
satisfied by $(\dbY, \dbZ)$ with $(\xi,\gamma)\1n=\1n(\dbG_3, \dbG_4)$ can be rewritten as (\ref{YZ3}).
Thus, combining Fubini theorem, we can get that $\sN$ is equal to $\h\sN$ defined by (\ref{sN}).
Further,
plugging (\ref{E4.11}) and (\ref{E240303}) into (\ref{E24032301}) and utilizing  the arbitrariness of $v$, we can obtain
\vskip-5mm
$$\ba{ll}
\h\sN(t)+\big(R(t)\1n+\1n(D^{\top}\2n\triangleleft \sP \triangleright \1n  D)(t)\big) \widehat{v}(t)=0,\ \ae t\in[0,T].
\ea$$
\vskip-2mm
\no  Therefore, we conclude that  $(\widehat{\Theta}_1,\widehat{\Theta}_2,\widehat{\Theta}_3, \widehat{v})$ satisfies (\ref{E2403207}).

\medskip

\textbf{The sufficiency}: Assume that $\sP\1n=\1n(\sP_1,\sP_2,\sP_3,\sP_4)\1n\in \1n\Upsilon[0,T]$ is a regular solution of the Riccati
system (\ref{E4.1}). Let $(\widehat{\Theta}_1,\widehat{\Theta}_2,\widehat{\Theta}_3,\widehat{v})$ satisfy (\ref{E2403207}).
It suffices to prove that
\vskip-5mm
\bel{E24033002}\ba{ll}
\ns\ds  \widetilde{\cJ}(\t,\varphi_1,\varphi_2;u_1^{\widehat{\Theta},\widehat{v}+\wt v})
 \ges \widetilde{\cJ}(\t,\varphi_1,\varphi_2;u_1^{\widehat{\Theta},\widehat{v}}),\;\;\;\forall \;(\t,\varphi_1,\varphi_2)\1n\in\1n\wt\cI, \tilde{v}\1n\in\1n\mathcal{U}[\t,T].
\ea
\ee
\vskip-1.5mm
\no We first give a representation of $\wt\cJ(\1n\t,\varphi_1,\f_2;u_1^{\widehat{\Theta},\widehat{v}+\wt v})$ by the Riccati
system (\ref{E4.1}) and  the duality principle. In the following, for notational simplicity, we denote
\vskip-6mm
\bel{Notation-use-0}\ba{ll}
\ns\ds \widetilde{u}(t)\=u_1^{\widehat{\Theta},\widehat{v}+\wt v}(t),\ \ \dot{\sP}_3(s,t)\=\frac{\partial \sP_3}{\partial t}(s,t), \ \ \dot{\sP}_4(s_1,s_2,t)\1n\=\1n\frac{\partial \sP_4}{\partial t}(s_1,s_2,t).
\ea
\ee%
\vskip-2.5mm
\no  Moreover, by $(\rm H1)$-$(\rm H2)$, H\"{o}lder inequality and $\sP\1n=\1n(\sP_1,\sP_2,\sP_3,\sP_4)\1n\in\1n \Upsilon[0,T]$,  it is easy to see that
$\dot{\sP}_3\1n\in \1nL^{2}(\triangle_*[0,T];\mathbb{R}^{d\times d})$
and $\dot{\sP}_4\1n\in\1n L^{2,2,1}(\Box_3[0,T];\mathbb{R}^{d\times d})$.

 To begin with, we represent the first term $\dbE\langle G \1n\wt\sX_2(T),\1n\wt\sX_2(T)\rangle$ in $\wt\cJ(\1n\t,\1n\varphi_1,\1n\f_2;\1n\wt u)$(recall (\ref{E1.9})).
  Applying the It\^{o}'s formula to
 $s\mapsto \wt\sX_2(s)^{\top}\sP_2(s)\wt\sX_2(s)$ on $[\t, T]$ yields
\vskip-6mm
$$\ba{ll}
\ns\ds \dbE[\wt\sX_2(T)^{\top}G\wt\sX_2(T)]-\varphi_2^{\top}\sP_2(\t)\varphi_2\\
\ns\ds=\dbE\1n\int_\t^T\2n\Big[\wt\sX_2(s)^{\top}\1n\dot{\sP}_2(s)\wt\sX_2(s)
\1n+\1n2\wt\sX_2(s)^{\top}\1n\sP_2(s)\big(A(T,s)\wt X_1(s)\1n+\1nB(T,s)\wt u(s)\2n+\2n  b(T,s)\big)\\[-1mm]
\ns\ds\q+\big(\wt X_1(s)^{\top}\1nC(T,s)^{\top}\2n+\1n\wt u(s)^{\top}\1nD(T,s)^{\top}\2n+\2n\si(T,s)^{\top}\1n\big)\sP_2(s)
\big(C(T,s)\wt X_1(s)\1n+\1nD(T,s)\wt u(s)\2n+\2n\si(T,s)\big)\Big]ds,
\ea
$$
\vskip-2mm
\no where $(\widetilde{X}_1,\widetilde{\sX}_1,\widetilde{\sX}_2)$
 is the extended causal feedback solution at $(\t,\varphi_1,\varphi_2)$ corresponding to
$(\widehat{\Theta}_1,\widehat{\Theta}_2,\widehat{\Theta}_3,\widehat{v}+\widetilde{v})$ (recall (\ref{E0601})).

We see that the new terms of the form $\wt X_1(s)^{\top}(\cdots)\wt X_1(s)$ and $\wt \sX_2(s)^{\top}(\cdots)\wt X_1(s)$ appear on the right hand
for $s\in[\t,T]$. Next, to deal with these new terms, we use the It\^{o}'s formula
  to $t\1n\mapsto\1n \wt \sX_1(s,t)^{\top}\1n\sP_1(s)\wt  \sX_1(s,t)$ and $t\1n\mapsto\1n \wt \sX_2(t)^{\top}\1n\sP_3(s,t)\wt \sX_1(s,t)$
  on $[\t, s]$, respectively. We have
\vskip-6.4mm
$$\ba{ll}
\ns\ds \dbE[\wt X_1(s)^{\top}\sP_1(s)\wt X_1(s)] =\dbE[\wt \sX_1(s,s)^{\top}\sP_1(s) \wt \sX_1(s,s)]\\
\ns\ds=\varphi_1(s)^{\top}\sP_1(s)\varphi_1(s)+\dbE\int_\t^s\Big[
2\wt \sX_1(s,r)^{\top}\sP_1(s)\big(A(s,r)\wt X_1(r)+B(s,r)\wt u(r)\2n+\2n  b(s,r)\big)\\[-1mm]
\ns\ds\q+\big(\wt X_1(r)^{\top}\1nC(s,r)^{\top}\2n+ \1n\wt u(r)^{\top}\1nD(s,r)^{\top}\2n+\2n\si(s,r)^{\top}\1n\big)\sP_1(s)
\big(C(s,r)\wt X_1(r)\1n+\1nD(s,r)\wt  u(r)\2n+\2n \si(s,r)\big)\Big]dr,
\ea
$$
\vskip-2.5mm
\no and (recall (\ref{Notation-use-0}) for $\dot{\sP}_3$)
\vskip-6.5mm
$$\ba{ll}
\ns\ds \dbE[\wt \sX_2(s)^{\top}\sP_3(s,s) \wt X_1(s)]=\dbE[\wt \sX_2(s)^{\top}\sP_3(s,s)\wt  \sX_1(s,s)]\\
\ns\ds=\1n\varphi_2^{\top}\1n\sP_3(s,\t)\varphi_1(s)\1n+\1n\dbE\2n\int_\t^s\2n\Big[\wt \sX_2(r)^{\top}\1n\dot{\sP}_3(s,r)\wt \sX_1(s,r)
\1n+\2n\wt \sX_2(r)^{\top}\1n\sP_3(s,r)\1n\big(A(s,r)\wt X_1(r)\2n+\2n  b(s,r)\\[-1mm]
\ns\ds\q+B(s,r)\wt  u(r)\big)+\big(\wt X_1(r)^{\top}A(T,r)^{\top} +\wt u(r)^{\top}B(T,r)^{\top}\2n+\2n  b(T,r)^{\top}\1n\big)\sP_3(s,r)\wt \sX_1(s,r)\\[-1mm]
\ns\ds\q+\1n\big(\wt X_1(r)^{\top}\1nC(T,r)^{\top}\3n+\1n \wt u(r)^{\top}\1nD(T,r)^{\top}\3n+\2n \si(T,r)^{\top}\1n\big)\sP_3(s,r)\1n\big(C(s,r)\wt X_1(r)\1n+\1nD(s,r) \wt u(r)\2n+\2n  \si(s,r)\big)\Big]dr.
\ea
$$
\vskip-2mm
\no
 Similarly,  we can utilize the It\^{o}'s formula to $\th\mapsto \wt \sX_1(s,\th)^{\top}\sP_4(s,r,\th)\wt \sX_1(r,\th)$ on $[\t, r]$
  to deal with the new term $\wt \sX_1(s,r)^{\top}(\cdots)\wt X_1(r)$.
 Finally, by Fubini theorem and Riccati
system (\ref{E4.1}), it follows from some calculations that
\vskip-6mm
$$\ba{ll}
\ns\ds \wt\cJ(\t,\f_1,\f_2;\wt  u)\1n=\1nI_1+I_2+I_3+I_4+I_5,
\ea
$$
\vskip-2mm
 \no where (recall the notations of  $\BR,\BD,\BB_1$ and $\BB_2$ in the Definition \ref{def2.11})
\vskip-7mm
$$\ba{ll}
\ns\ds I_1\1n\=\1n\varphi_2^{\top}\1n\sP_2(\t)\varphi_2\1n+\2n\int_\t^T\3n\Big[\varphi(s)\1n^{\top}\2n\sP_1(s)\varphi_1(s)
\2n+\2n2\varphi_2^{\top}\2n\sP_3(s,\t)\varphi_1(s)\1n+\2n\int_\t^T\3n\1n\varphi_1(s)^{\top}\2n\sP_4(s,r,\t)\varphi_1(r)dr\Big]ds,\\[-1mm]
\ns\ds I_2\1n\=\1n\dbE\1n\int_\t^T\2n\Big\{\wt u(s)^{\top}\Big[\BR(s)\wt  u(s)
\1n+\1n2\BD(s)\wt X_1(s)\1n+\1n 2\BB_2(s)\wt \sX_2(s)
\1n+\1n2\1n\int_s^T\3n\BB_1(r,s)\wt \sX_1(r,s)dr\Big]\Big\}ds,\\[-1mm]
\ns\ds I_3\1n\=\1n\dbE\1n\int_\t^T\2n\Big\{\wt \sX_2(s)\1n^{\top}\2n\dot{\sP}_2(s)\wt \sX_2(s)\2n+
\2n\int_s^T\3n\Big[2\wt \sX_2(s)\1n^{\top}\2n\dot{\sP}_3(r,s)
\1n+\2n\int_s^T\3n\wt \sX_1(\th,s)\1n^{\top}\2n\dot{\sP}_4(\th,r,s)d\th \Big]\wt \sX_1(r,s)dr\Big\}ds,\\[-1mm]
\ns\ds I_4\1n\=\1n\dbE\1n\int_\t^T\3n\Big\{\Big[\wt X_1(s)^{\top}\BD(s)^{\top}\2n+\1n2\wt \sX_2(s)^{\top}\BB_2(s)^{\top}
\2n+\1n2\2n\int_s^T\3n\wt\sX_1(r,s)^{\top}\BB_1(r,s)^{\top}dr\Big] \BR(s)^\dagger\BD(s)\wt X_1(s)\Big\}ds,\\[-1mm]
\ns\ds I_5\1n\=\1n\dbE\1n\int_\t^T\2n\Big\{(\si^{\top}\2n\triangleleft\sP\triangleright  \si)(s)\2n+\1n2\wt X_1(s)^{\top}\1n\big(q(s)\2n+\2n
(C^{\top}\2n\triangleleft\sP\triangleright  \si)(s)\big)\2n+\1n2\wt u(s)^{\top}\1n\big(\rho(s)\2n+\2n(D^{\top}\2n\triangleleft\sP\triangleright  \si)(s)\big)\\[-1mm]
\ns\ds \q + 2\wt \sX_2(s)^{\top}(\sP_{2,3}\triangleright  b)(s)\2n+\2n \int_s^T\2n 2 \wt\sX_1(r,s)^{\top}(\sP_{1,3,4}\triangleright  b)(r,s)dr \Big\} ds.
\ea
$$
\vskip-3mm

We first look at the term $I_2$.
By the representation (\ref{E2403207}) of $(\widehat{\Theta}_1,\widehat{\Theta}_2,\widehat{\Theta}_3,\widehat{v})$, we observe that
\vskip-7mm
$$\ba{ll}
\ns\ds -\BR(t)\widetilde{u}(t)=
\BD(t)\widetilde{X}_1(t)\1n+\2n\int^T_t\2n\BB_1(r,t)\widetilde{\sX}_1(r,t)dr
\1n+\1n\BB_2(t)\widetilde{\sX}_2(t)\1n+\1n\h\sN(t)\1n-\1n\BR(t)\widetilde{v}(t).
\ea
$$
\vskip-3mm
\no Thus the term $I_2$ can be  rewritten  as
\vskip-5.3mm
\bel{E24033001}\ba{ll}
\ns\ds I_2=\dbE\2n\int_\t^T\Big\{\wt u(s)^{\top}\BR(s)\big[2\widetilde{v}(s) -2\h\sN(s)-\wt u(s)\Big\}ds.
\ea
\ee
\vskip-3mm
\no We compute the second term $\wt u(s)^{\top}\BR(s)\wt u(s)$ in (\ref{E24033001}). It follows from  (\ref{E2403207}) that
\vskip-7mm
$$\ba{ll}
\ns\ds \wt u(s)^{\top}\BR(s)\wt u(s)\1n=\1n\Big[\wt X_1(s)^{\top}\BD(s)^{\top}\2n+\1n2\wt \sX_2(s)^{\top}\BB_2(s)^{\top}\2n
+\1n2\2n\int_s^T\3n\wt\sX_1(r,s)^{\top}\BB_1(r,s)^{\top}\1n
dr\Big] \BR(s)^\dagger \BD(s)\wt X_1(s)\\[-1.8mm]
\ns\ds\qq+\2n\int_s^T\3n\int_s^T\3n\wt\sX_1(r,\1ns)\1n^{\top}\BB_1(r,\1ns)\1n^{\top}\BR(s)^\dagger
\BB_1(\th,\1ns)\wt\sX_1(\th,\1ns)d\th dr\2n+\2n\big(\widehat{v}(s)\2n+\2n\widetilde{v}(s)\big)\1n^{\top}\1n\BR(s)\big(2\wt u(s)\1n-\1n
\widehat{v}(s)\2n-\2n\widetilde{v}(s)\big) \\[-1.8mm]
\ns\ds\qq+\1n\wt\sX_2(s)^{\top}\1n\BB_2(s)^{\top}\BR(s)^\dagger
\Big[\BB_2(s)\wt\sX_2(s)\1n
+\1n2\2n\int_s^T\2n\BB_1(r,s)\wt\sX_1(r,s)dr\Big].
\ea
$$
\vskip-3mm
\no Plugging the above into (\ref{E24033001}) and recalling the Riccati system (\ref{E4.1}), we can get
\vskip-5.5mm
\bel{E24120803}\ba{ll}
\ns\ds I_2\1n+\1nI_3\1n+\1nI_4\1n=\1n I_{234}-\dbE\1n\int_\t^T\1n2\widetilde{v}(s)^{\top}\h \sN(s)ds,
\ea
\ee
\vskip-3mm
\no where $I_{234}$ is defined by
\vskip-6.5mm
$$\ba{ll}
\ns\ds I_{234}\1n\=\1n \dbE\1n\int_\t^T\2n\big(\widetilde{v}(s)^{\top}\BR(s)\widetilde{v}(s)+\h v(s)^{\top}\BR(s)\h v(s)\big)ds.
\ea
$$
\vskip-3mm

On the other hand, we observe  the fact of
\vskip-6mm
\bel{E24121001}\ba{ll}
 \ns\ds I_1+I_{234}=\widetilde{\cJ}^0(\t,\varphi_1,\varphi_2;u_1^{\widehat{\Theta},\widehat{v}+\wt v}).
\ea
\ee
\vskip-2mm
\no By (\ref{E24121001}) and $\BR(\cd)\1n\ges\1n 0$,  it is easy to see that, for any $(\t,\varphi_1,\varphi_2)\in\wt\cI$ and $\tilde{v}\in\mathcal{U}[\t,T]$,
\vskip-5mm
\bel{E24121002}\ba{ll}
 \ns\ds \widetilde{\cJ}^0(\t,\varphi_1,\varphi_2;u_1^{\widehat{\Theta},\widehat{v}+\wt v})
\1n-\1n \widetilde{\cJ}^0(\t,\varphi_1,\varphi_2;u_1^{\widehat{\Theta},\widehat{v}})\1n=\1n\dbE\1n\int_\t^T\2n\tilde{v}(s)^{\top}\BR(s)\tilde{v}(s)ds\1n\ges\1n 0.
\ea
\ee
\vskip-2.5mm
\no In addition, combining Corollary \ref{cor3.4} and (\ref{E24121002}), we can obtain that (\ref{E24032302}) holds.
Using the same argument to the proof of the necessity part, we can obtain the relationship between the system (\ref{P12})
 and  Riccati system (\ref{E4.1}), i.e. (\ref{E4.6})--(\ref{E240304}). Consequently,  it follows that (\ref{G4}) holds.

Thus, by  (\ref{E6}), (\ref{G4}) and the definition of $(\widehat{\dbG}_3, \widehat{\dbG}_4)$, $I_5$ can be  rewritten  as
\vskip-5.5mm
$$\ba{ll}
\ns\ds I_5=\dbE\1n\int_\t^T\2n\Big\{(\si^{\top}\2n\triangleleft\sP\triangleright  \si)(s)\2n
+2 \wt\dbX(s)^{\top}\widehat{\dbG}_3(s) \1n+\2n \int_s^T\2n 2  \wt\cX(r,s)^{\top}\widehat{\dbG}_4(r,s) dr\Big\} ds,
\ea
$$
\vskip-3mm
\no where
\vskip-8mm
$$\ba{ll}
\ns\ds \wt \dbX(s)^{\top}\=\Big(\3n\begin{array}{ccccc}
\wt X_1(s)^{\top}\2n &   \wt u(s)^{\top}\2n& \wt \sX_2(s)^{\top}
\end{array}\3n\Big),\;\;
\wt \cX(r,s)^{\top}\=\Big(\3n\begin{array}{ccccc}
\wt \sX_1(r,s)^{\top}\2n &   \wt u(s)^{\top}\2n& \wt \sX_2(s)^{\top}
\end{array}\3n\Big).
\ea$$
\vskip-2.5mm
\no  Then applying the duality principle  (Lemma \ref{duality}) yields that
\vskip-5.5mm
\bel{I4}\ba{ll}
\ns\ds I_5=\dbE\1n\int_\t^T\2n\Big\{(\si^{\top}\2n\triangleleft\sP\triangleright  \si)(s)\1n
+2 \wt v(s)^{\top}\h\sN(s)\1n+ 2\widehat{\Phi}^{\f}(s)^{\top}\dbY(s,\t)\\[-1mm]
\ns\ds\qq\qq\q+\1n \int_s^T\2n 2 \big[\langle\dbY(r,s),\widehat{\Phi}^b(r,s) \rangle\1n+\1n\langle\dbZ(r,s),\widehat{\Phi}^\si(r,s) \rangle\big] dr\Big\}ds,
\ea
\ee
\vskip-2mm
\no where $(\dbY, \dbZ)$ solves (\ref{YZ3}) and $\h\sN$ is defined by (\ref{sN}). 

Combining  the above and  (\ref{E24120803}), we can obtain that
\vskip-6mm
\bel{E24032102}\ba{ll}
\ns\ds
 \widetilde{\cJ}(\t,\varphi_1,\varphi_2;\wt u)=I_1+ \dbE\1n\int_\t^T\2n\Big\{\widehat{v}(s)^{\top}\BR(s)\widehat{v}(s)\1n+\1n\widetilde{v}(s)^{\top}\BR(s)\widetilde{v}(s)
\1n+\1n (\si^{\top}\2n\triangleleft\sP\triangleright  \si)(s)\1n\\[-0.5mm]
\ns\ds\q+ 2\widehat{\Phi}^{\f}(s)^{\top}\1n\dbY(s,\t)\1n+\2n \int_s^T\2n 2 \big[\langle\dbY(r,s),\widehat{\Phi}^b(r,s) \rangle\1n+\1n\langle\dbZ(r,s),\widehat{\Phi}^\si(r,s) \rangle\big] dr\Big\}ds.
\ea
\ee
\vskip-3mm
\no Therefore, for any $(\t,\varphi_1,\varphi_2)\in\wt\cI$ and $\tilde{v}\in\mathcal{U}[\t,T]$, it is easy to see that
\vskip-6mm
$$\ba{ll}
 \ns\ds \widetilde{\cJ}(\t,\varphi_1,\varphi_2;u_1^{\widehat{\Theta},\widehat{v}+\wt v})
\1n-\1n \widetilde{\cJ}(\t,\varphi_1,\varphi_2;u_1^{\widehat{\Theta},\widehat{v}})=\dbE\1n\int_\t^T\2n\tilde{v}(s)^{\top}\BR(s)\tilde{v}(s)ds.
\ea
$$
\vskip-2mm
\no Since $\sP\1n=\1n(\sP_1,\sP_2,\sP_3,\sP_4)\1n\in\1n \Upsilon[0,T]$ is a regular solution of the Riccati
 system (\ref{E4.1}), that is, $\BR(t)\1n\ges\1n 0, t\1n\in\1n[0,T]$, we see that
(\ref{E24033002}) holds
for any $(\t,\varphi_1,\varphi_2)\1n\in\1n\wt\cI$ and $\tilde{v}\1n\in\1n\mathcal{U}[\t,T]$,
which implies that
 $(\widehat{\Theta}_1,\widehat{\Theta}_2,\widehat{\Theta}_3,\widehat{v})$ is a strongly optimal causal feedback strategy.
Further, it follows from (\ref{E24032102}) that the representation of $\wt V(\t, \f_1,\f_2)$ satisfies (\ref{E24032601}).
The proof is complete.
\endpf

\begin{remark}\label{remark3.8}
We point out that the above proof, especially the necessity part, indicates one important reason of introducing extended notions in Section 2. In fact, by the previous procedures we see the arbitrariness of $\varphi_1(\cd)$ and $\f_2$ is frequently used and becomes very crucial in obtaining (\ref{E4.6}),  (\ref{E4.7}) and eventually the representation of
 $(\widehat{\Theta}_1,\1n\widehat{\Theta}_2,\1n\widehat{\Theta}_3,\1n\widehat{v})$, the Riccati system (\ref{E4.1}).
Nevertheless, if we follow the \it original \rm \it framework by using $(\t,\f)\in\cI$, we then obtain
(\ref{E24032302}) with $\f_1(\cd)\=\f(\cd)$, $\f_2\=\f(T)$. The point is we may fail to obtain
(\ref{E4.6}), (\ref{E4.7}) due to the lack of enough freedom for $\varphi$, not to mention the other desired conclusions.
Therefore, the set of input condition $\widetilde{\cI}$ allows us to have more \it test conditions \rm \it than $\cI$. We point out that
such an obstacle does not appear in the Lagrange cost functional of \cite{Hamaguchi-Wang-I, Hamaguchi-Wang-II}.
\end{remark}

\begin{remark}\label{remark3.12}
Let us point out an interesting fact derived from Theorem \ref{the3},
i.e., the value of $\h\Th_2(s,t)$ with $s\les t$ only makes sense
if when
$R(t)\1n+\1n(D^{\top}\1n\triangleleft \sP \triangleright   D)(t)\1n=\1n0$. In other words, if
$D\1n=\1n0$, and $R(\cd)\1n>\1n0$, then this term vanishes.
\end{remark}

 The following result shows the uniqueness of regular solution to Riccati system (\ref{E4.1})
 in Theorem \ref{the3}.

\begin{corollary}\label{cor4.2}
The Riccati system (\ref{E4.1}) has at most one regular solution.
\end{corollary}

\proof
 Assume that $(\sP_1,\sP_2,\sP_3,\sP_4),\big(resp.(\wt\sP_1,\wt\sP_2,\wt\sP_3,\wt\sP_4)\big)\1n\in\1n \Upsilon[0,T]$
is a regular solution of the Riccati system (\ref{E4.1}). We consider  the homogeneous case, i.e., $b\1n=\1n\si\1n=\1nq\1n=\1n\rho\1n=\1n0$.
For any $(\t,\varphi_1,\varphi_2)\1n\in\1n\wt\cI$, the representation (\ref{E24032601}) of the extended value function  and the uniqueness of the infimum
yield that
\vskip-7mm
\bel{E2413}\ba{ll}
\ns\ds \1n\f_2^{\top}\2n\widehat{\sP}_2(\t)\f_2\2n+\3n\int_\t^T\3n\Big[\f_1\1n(s)^{\top}\2n\widehat{\sP}_1\1n(s)\2n+\2n2\f_2^{\top}\2n\widehat{\sP}_3(s,\t)
\2n+\3n\int_\t^T\3n\2n\f_1\1n(\th)^{\top}\2n\widehat{\sP}_4(\th,s,\t)d\th\Big]\1n\f_1\1n(s)ds\2n=\2n0,
\ea
\ee
\vskip-2.5mm
\no where $\widehat{\sP}_i\=\sP_i-\wt\sP_i,i\in\{1,2,3,4\}$.

By the  arbitrariness of $\varphi_2$, let $\varphi_2=0$.  From Lemma \ref{lem4.1}, we get that
\vskip-6.5mm
\bel{E2415}\ba{ll}
\ns\ds \widehat{\sP}_1(t)\1n=\1n0,\;$a.e.$\; t\in[0,T];\;\;\widehat{\sP}_4(s,t,r)\1n=\1n0,\;$a.e.$\;(s,t)\in[0,T]^2,\;\forall r\in[0,s\wedge t].
\ea
\ee
\vskip-2.5mm
\no Choosing $\varphi_1=0$ in (\ref{E2413}), it is clear that
\vskip-6mm
\bel{E2416}\ba{ll}
\ns\ds \widehat{\sP}_2(t)=0,\;\;$a.e.$\; t\in[0,T].
\ea
\ee
\vskip-2mm
\no From (\ref{E2413}), (\ref{E2415}) and (\ref{E2416}), it follows from Lemma \ref{lem4.1} that
\vskip-6mm
\bel{E2417}\ba{ll}
\ns\ds \widehat{\sP}_3(s,t)=0,\;\;$a.e.$\; s\in[0,T],\;\forall t\in[0,s].
\ea
\ee
\vskip-2mm
\no Combining (\ref{E2415})--(\ref{E2417}), we complete our proof.
\endpf

\begin{remark}\label{remark3.9}
By the above proof, we see that the arbitrariness of $\f_1$ and $\varphi_2$ is essential to ensure the uniqueness of the regular solution. In other words, if we follow the framework of $\cI$ rather $\wt\cI$, we may not derive (\ref{E2415}) and (\ref{E2416}) due to the lack of enough \it test \rm initial conditions. This gives another reason of introducing the \it extended language \rm \it in Section 2.
\end{remark}

Next we show that under proper conditions, there exists $(\t,\varphi_1,\varphi_2)\in\widetilde{\cI}$ such that the optimal control, if it exist, can not be the following form,
\vskip-4mm
$$\ba{ll}
\ns\ds \h u(t)\=\h\Th_1(t)\h X(t)+\int_\t^t\Th_2(s,t)\h X(s)ds+\h v(t), \ \ t\in[\t,T],
\ea
$$
\vskip-1.5mm
\no which gives an answer to $(\rm Q3)$ in the Introduction from another standpoint.

\begin{corollary}\label{corollary3.10}
Suppose $(\rm H1)$-$(\rm H2)$ hold s.t. $A\1n\equiv\1n B$ in (\ref{E1}) is invertible, $S\1n\equiv\1n 0$ and
\vskip-6mm
\bel{Counter-coefficients}\ba{ll}
\ns\ds \1nQ(t)\2n+\2nC(T,\1nt)^{\top}\1nGC(T,\1nt)\2n-\2nC(T,\1nt)^{\top}\1nGD(T,\1nt)
\Big[R(t)\2n+\2nD(T,\1nt)\1n^{\top}\1nGD(T,\1nt)\2n\Big]^\dagger D(T,t)^{\top}\1nGC(T,t)\2n\neq 0
\ea\ee
\vskip-2mm
\no for $t\in[\t,T].$ Assume the 4-tuple $(\widehat{\Theta}_1,\widehat{\Theta}_2,\h \Th_3,\widehat{v})\1n\in\1n \mathcal{S}[0,T]$ is
strongly optimal. Then it must not be the form of $(\widehat{\Theta}_1,\widehat{\Theta}_2,0,\widehat{v})$
with $\h\Th_2(r,t)=0$ $(r>t)$.
\end{corollary}

\proof We prove the conclusion by contradiction.
Suppose Problem (LQ-SVIE) admits a strongly optimal $(\widehat{\Theta}_1,\widehat{\Theta}_2,0,\widehat{v})$ such that
$\h\Th_2(r,t)=0$ with $r>t$.
Then, following the necessity proof in Theorem \ref{the3} we can obtain
\vskip-6mm
$$\ba{ll}
\ns\ds (B^{\top}\1n\triangleleft \sP_{2,3})(t)=0,\; t\1n\in\1n[0,T],\;\;(B^{\top}\1n\triangleleft\sP_{1,3,4})(r,t)=0,\;(r,t)\1n\in\1n\triangle_*[0,T],
\ea
$$
\vskip-1.5mm
\no where   $\sP=(\sP_1,\sP_2,\sP_3,\sP_4)$ is a regular solution of  the following system
\vskip-5mm
$$\3n\3n\left\{\3n\1n\ba{ll}
\ns\ds \sP_1(t)=Q(t)+(C^{\top}\1n\triangleleft \sP \triangleright  C)(t)\\
\ns\ds\q- (C^{\top}\2n\triangleleft\1n \sP \1n\triangleright \1n D)(t)
\big(R(t)\1n+\1n(D^{\top}\2n\triangleleft\1n \sP \1n\triangleright  \1n D)(t)\big)^\dagger\1n
(D^{\top}\2n\triangleleft\1n \sP \triangleright\1n  C)(t),\;t\1n\in\1n[0,T],\\
\ns\ds \sP_2=G,\;\;\sP_3=\sP_4=0.
\ea\right.
$$
\vskip-2mm
\no
Clearly, in this case, we see that
\vskip-7mm
$$\ba{ll}
\ns\ds \sP_1(t)\1n=\1nQ(t)\1n+\1nC(T,\1nt)^{\top}\1nGC(T,\1nt)\2n+\3n\int_t^T\3nC(s,\1nt)^{\top}\1n\sP_1(s)C(s,\1nt)ds
-\Big[\1nC(T,\1nt)^{\top}\1nGD(T,\1nt)\\[-1mm]
\ns\ds\qq+\3n\int_t^T\3nC(s,t)\1n^{\top}\1n\sP_1(s)\1nD(s,\1nt)ds\Big]\Big[R(t)\1n
+\1nD(T,\1nt)\1n^{\top}\1nGD(T,\1nt)\2n+\3n\int_t^T\3nD(s,\1nt)\1n^{\top}\1n\sP_1(s)\1nD(s,\1nt)ds\Big]^\dagger\\[-1mm]
\ns\ds\qq\times\Big[D(T,t)^{\top}\1nGC(T,t)\2n+\3n\int_t^T\3nD(s,t)^{\top}\1n\sP_1(s)C(s,t)ds\Big].
\ea
$$
\vskip-2mm
\no By the above (\ref{Counter-coefficients}), we see that $\sP_1$ is nonzero.
On the other hand, using the facts that $(B^{\top}\1n\triangleleft\sP_{1,3,4})(r,t)=0$ and $\sP_3=\sP_4=0$, and $B$ is invertible,
we obtain that $\sP_1=0$, which is a contradiction.
\endpf

We point out that (\ref{Counter-coefficients}) is easy to check. For example, in the one-dimensional case, if $Q(\cd)\ges0$,
$R(\cd)\ges0$, $G\ges0$ such that either $R(\cd)>0$ or $G>0$, then (\ref{Counter-coefficients}) hold true.

\section{Well-posedness of Riccati system}

In this section, let us discuss the solvability of Riccati system (\ref{E4.1}).
We call its regular solution $\sP\1n=\1n(\sP_1,\1n\sP_2,\1n\sP_3,\1n\sP_4)$ \emph{strongly regular} if there exists a constant $\lambda\1n>\1n 0$ such that
$ R(t)\1n+\1n(D^{\top}\1n\triangleleft \sP \triangleright   D)(t)\1n\ges \1n\lambda I_l$ for a.e. $\1nt\1n\in\1n[0,T]$.
 For notational simplicity, we consider  the   case $b\1n=\1n\si\1n=\1nq\1n=\1n\rho\1n=\1nS\1n=\1n0$. 
To begin with, for any $(\Theta_1,\1n\Theta_2,\1n\Theta_3)\1n\in\1n L^\infty(0,T;\mathbb{R}^{l\times d})\1n\times \1nL^2(\triangle_{*}[0,T];\mathbb{R}^{l\times d})\1n\times \1n L^\infty(0,T;\mathbb{R}^{l\times d})$, define
\vskip-6mm
\bel{E241103}\ba{ll}
\ns\ds \2n Q_1(t)\1n\=\1n\Th_1(t)\1n^{\top}\1nR(t)\Th_1(t)\1n+\1nQ(t),\;\;Q_2(t)\1n\=\1n\Th_3(t)\1n^{\top}\1nR(t)\Th_3(t),\;t\in[0,T],\\
\ns\ds\2n Q_4(s,t)\1n\=\1n\Th_3(t)\1n^{\top}\1nR(s)\Th_2(s,t),\;\;Q_6(s,t)\1n\=\1n\Th_2(s,t)\1n^{\top}\1nR(s)\Th_1(t),\;(s,t)\1n\in\1n\triangle_{*}[0,T],\\
\ns\ds \2n Q_3(t)\2n\=\2n\Th_3\1n(t)\1n^{\top}\2nR(t)\Th_1\1n(t),t\2n\in\2n[0,\1nT],
Q_5(s,\1nt,\1n\th)\2n\=\2n\Th_2(\1ns,\1n\th)\1n^{\top}\2nR(\th)\Th_2(\1nt,\1n\th),(\1ns,\1nt,\1n\th)\2n\in\2n\Box_3[0,\1nT].
\ea\ee
\vskip-1.5mm
\no It is easy to see that $\big(Q_1, Q_2, Q_3, Q_4, Q_5, Q_6\big)\1n\in\1n \cY[0,T]\1n\=\1n L^\infty(0,T;\mathbb{S}^{d})\times L^\infty(0,T;\mathbb{S}^{d})\times L^\infty(0,T;$
 $\mathbb{R}^{d\times d})\times L^2(\triangle_{*}[0,T];\mathbb{R}^{d\times d})\times L^{2,2,1}_{sym}(\Box_3[0,T];\mathbb{R}^{d\times d})\times
 L^2(\triangle_{*}[0,T];\mathbb{R}^{d\times d})$.

Now we introduce the following Lyapunov system:
\vskip-6mm
\bel{E241101}\left\{\3n\ba{ll}
\ns\ds \2n\sP_1(t)\1n=\1nf_1[\Th;\sP](t)+Q_1(t),t\in[0,T],\\
\ns\ds \2n \sP_2(t)\1n=G+\2n\int_t^T\2n\big(f_2[\Th;\sP](s)+Q_2(s)\big)ds,\;t\in[0,T],\\[-1.5mm]
\ns\ds\2n \sP_3(t,r)\2n=\2nf_3[\Th;\1n\sP](t)\2n+\2nQ_3(t)
\2n+\3n\int_r^t\2n\big(f_4[\Th;\1n\sP](t,s)\2n+\2nQ_4(t,s)\big)ds,(t,\1nr)\1n\in\1n\triangle_{*}[0,T],\\[-1mm]
\ns\ds \1n\sP_4(s,t,r)\1n=\1n\sP_4(t,s,r)^{\top} \\
\ns\ds\q =\1n\Big(f_6[\Th;\1n\sP](s,t)\2n+\2nQ_6(s,t)\Big)I_{[t,T)}(s)\1n+\1n\Big(f_6[\Th;\1n\sP](t,s)^{\top}\2n+\2nQ_6(t,s)^{\top}\1n\Big)I_{[0,t)}(s) \\
\ns\ds\qq +\2n\int_r^{s\wedge t}\2n\big(f_5[\Th;\sP](t,s,\th)+Q_5(s,t,\th)\big)d\th,\;(s,t,r)\in\Box_3[0,T],
\ea\right.
\ee
\vskip-1.5mm
\no where 
\vskip-1.5mm
\vskip-6mm
$$\ba{ll}
\ns\ds \2n f_1[\Th;\sP](t)\1n\=(C^{\top}\1n\triangleleft \sP \triangleright  C)(t)\2n+\2n \BD(t)^{\top}\1n\Th_1(t)
\2n+\2n \Th_1(t)^{\top}\1n\BD(t)\2n+\2n\Th_1(t)^{\top}\1n(D^{\top}\1n\triangleleft\1n \sP \1n\triangleright  \1n D)(t)\Th_1(t),\\
\ns\ds \2nf_2[\Th;\sP](t)\1n\=\BB_2(t)^{\top}\1n\Th_3(t)\1n+\1n\Th_3(t)^{\top}\1n\BB_2(t)
\1n+\1n\Th_3(t)^{\top}\1n(D^{\top}\2n\triangleleft \1n\sP \1n\triangleright  \1n D)(t)\Th_3(t),\\
\ns\ds\2nf_3[\Th;\sP](t)\1n\=(\1n\sP_{2,3}\1n\triangleright\1n A)(t)\2n+\2n\BB_2(t)^{\top}\1n\Th_1(t)
\2n+\2n\Th_3(t)^{\top}\1n\BD(t)\2n+\2n\Th_3(t)^{\top}\1n(D^{\top}\3n\triangleleft \1n\sP \1n\triangleright  \1n D)(t)\Th_1(t),\\
\ns\ds\2n f_4[\Th;\sP](t,s)\1n\=\BB_2(s)^{\top}\1n\Th_2(t,s)
\1n+\1n\Th_3(s)^{\top}\1n\BB_1(t,s)\1n+\1n\Th_3(s)^{\top}\1n(D^{\top}\3n\triangleleft \1n\sP \1n\triangleright  \1n D)(s)\Th_2(t,s),\\
\ns\ds\2n f_5[\Th;\sP](s,t,\th)\1n\=\BB_1(s,\1n\th)\1n^{\top}\1n\Th_2(t,\1n\th)
\2n+\2n\Th_2(s,\1n\th)\1n^{\top}\1n\BB_1(t,\1n\th)\2n+\2n\Th_2(s,\1n\th)\1n^{\top}\1n(D^{\top}\3n\triangleleft \1n\sP \1n\triangleright  \1n D)(\th)\Th_2(t,\1n\th),\\
\ns\ds\2n f_6[\Th;\1n\sP](s,\1nt)\2n\=\2n(\1n\sP_{1,3,4}\1n\triangleright\1n A\1n)(t)\2n+\1n\BB_1\1n(\1ns,\1nt)\1n^{\top}\1n\Th_1\1n(t)
\2n+\2n\Th_2(\1ns,\1nt)\1n^{\top}\1n\BD(t)\2n+\2n\Th_2(\1ns,\1nt)\1n^{\top}\1n(\1nD^{\top}\3n\triangleleft \1n\sP \1n\triangleright  \1n D\1n)(s)\Th_1\1n(t).
\ea$$
\vskip-1.5mm
\no Here,  due to notational simplicity, for the solution $\sP\1n=\1n(\sP_1,\1n\sP_2,\1n\sP_3,$ $\1n\sP_4)$ of the above Lyapunov system (\ref{E241101}), we still denote
$\BD\1n\=\1n (D^{\top}\1n\triangleleft\1n \sP\1n \triangleright \1n C)$,
$\BB_1\1n\=\1n(B^{\top}\1n\triangleleft\1n \sP_{1,3,4})$,
$\BB_2\1n\=\1n(B^{\top}\1n\triangleleft \1n\sP_{2,3})$.

Thus we have the following result.

\bl\label{lem5-1}
Let $(\rm H1)$-$(\rm H2)$ hold and $(\Theta_1,\1n\Theta_2,\1n\Theta_3,\1nv)\1n\in\1n\mathcal{S}[0,T]$. For each $\big(Q_1,\1n Q_2,\1n Q_3,$
$ Q_4,\1n Q_5,\1n Q_6\big)\1n\in\1n \cY[0,T]$ defined by (\ref{E241103}), there exists a unique solution
$\sP\1n=\1n(\sP_1,\1n\sP_2,\1n\sP_3,$ $\1n\sP_4)\1n\in\1n\Upsilon[0,T]$ to the Lyapunov system (\ref{E241101}).
Furthermore, for any $(\t,\1n\varphi_1,\1n\varphi_2)\1n\in\1n\widetilde{\cI}$,  the following hold:
\vskip-6mm
$$\ba{ll}
\ns\ds \dbE\1n\bigg\{\1n\sX_2(T)\1n ^{\top}\1n G\sX_2(T)\1n+\3n\int_\t^T\3n\1n\Big[X_1(t)\1n ^{\top}\1n Q_1(t)X_1(t)\1n
+\1n\sX_2(t)\1n ^{\top}\1n  Q_2(t)\sX_2(t)\2n+\1n2\sX_2(t)\1n ^{\top}\1n  Q_3(t)X_1(t)\\[-1mm]
\ns\ds\q\2n+ \3n\int_t^T\3n\2n\Big( 2\sX_1(s,\1nt)\1n^{\top}\2n \big(Q_4(s,\1nt) \1n^{\top}\2n\sX_2(t)\2n+\2nQ_6(s,\1nt)X_1(t) \big)
\2n+\3n\int_t^T\3n\2n\sX_1(s,\1nt)\1n^{\top}\2n  Q_5(s,\1nr,\1nt)\sX_1(r,\1nt) dr\1n \Big) ds\Big] dt\1n\bigg\}\\
\ns\ds =\2n\varphi_2\1n^{\top}\2n\sP_2(\t)\varphi_2\2n+\3n\int_\t^T\2n\Big[\varphi_1\1n(s)\1n^{\top}\3n\sP_1(s)\varphi_1\1n(s)
\2n+\1n2\varphi_2\1n^{\top}\2n\sP_3(s,\1n\t)\varphi_1\1n(s)\2n+\3n\int_\t^T\3n\2n\varphi_1(s)\1n^{\top}\3n\sP_4(\1ns,\1nr,\1n\t)\varphi_1(r)dr\1n\Big]\1nds\\[-1mm]
\ns\ds\q+\dbE\1n\int_\t^T\3n v(t)\1n^{\top}\1n \big(R(t)\1n+\1n(D^{\top}\2n\triangleleft\1n \sP \1n\triangleright\1n   D)(t)\big)v(t) dt,
\ea$$
\vskip-1.5mm
\no where $(X_1,\sX_1,\sX_2)$ is the extended causal feedback solution at $(\t,\varphi_1,\varphi_2)$ corresponding to
$(\Theta_1,\Theta_2,\Theta_3,v)$ (recall (\ref{E0601})).

\el

\proof

First, we take
$\cP^{(1)}(\cd)\2n\=\2n\Big(\2n\begin{array}{ccccc}
  \sP_1(\cd)\2n &\2n 0\\[-0.7mm]
  0\2n &\2n0
\end{array}
\2n\Big)$, and
\vskip-4mm
$$
\ba{ll}
\ns\ds   \cP^{(2)}(s,r,t)\2n\=\2n\Big(\2n\begin{array}{ccccc}
 \sP_4(s,r,t)\1n & \3n\sP_3(s,t)^{\top}I_{(t,T)}(r)   \\[-0.3mm]
\sP_3(r,t)I_{(t,T)}(s)\1n & \3n\sP_2(t)
\end{array}
\2n\Big),\;\;(s,r,t)\in \Box_3[0,T],
\ea$$
\vskip-3mm
\no  and
\vskip-6.5mm
$$
\ba{ll}
\ns\ds \Xi(t)\1n\=\1n\left(\2n\begin{array}{ccccc}
\Th_1(t) & 0
\end{array}
\2n\right),t\in[0,T],\;\;\;
\Gamma(s,t)\1n\=\1n\left(\2n\begin{array}{ccccc}
\Th_2(s,t) & \Th_3(t)
\end{array}
\2n\right),(s,t)\1n\in\1n \triangle_*[0,T].
\ea$$
\vskip-1.5mm
\no  It's easy to see that $(\Xi,\1n\Gamma)\2n\in\2n L^\infty(0,T;\1n\mathbb{R}^{l\1n\times\1n 2d})\1n\times\1n L^2(\triangle_{*}[0,T];\1n\mathbb{R}^{l\1n\times\1n 2d})$.
Denote $\cP\1n=\1n(\cP^{(1)},\1n\cP^{(2)})$,

%
%
%
%
%
  %
%
  %
%
%
%
  %
%
%
%
%
%
  %
%
%
%
  %
%
%
%
%
\vskip-5mm
$$
\ba{ll}
\cA(s,t)\1n\=\1n\Big(\2n\begin{array}{ccccc}
  A(s,t)   \2n  & \2n0 \\[-0.3mm]
\frac{1}{T-t}A(T,t) \2n  &\2n 0
\end{array}
\2n\Big),\;\;\cB(s,t)\1n\=\1n\Big(\2n\begin{array}{ccccc}
 B(s,t)     \\[-0.3mm]
\frac{1}{T-t}B(T,t)
\end{array}
\2n\Big),\;(s,t)\in \triangle_*[0,T],\\
\cC(s,t)\1n\=\1n\Big(\2n\begin{array}{ccccc}
  C(s,t)  \2n   & \2n0 \\[-0.3mm]
\frac{1}{T-t}C(T,t) \2n  & \2n0
\end{array}
\Big),\;\;\cD(s,t)\1n\=\1n\Big(\2n\begin{array}{ccccc}
D(s,t)     \\[-0.3mm]
\frac{1}{T-t}D(T,t)
\end{array}
\2n\Big),\;(s,t)\in \triangle_*[0,T].
\ea$$
\vskip-1.5mm
\no Here, we recall the notations $\ltimes$ and $\rtimes$  in \cite{Hamaguchi-Wang-I,Hamaguchi-Wang-II}, i,e.,
 for each $M:\1n\triangle_*[0,T]\1n\rightarrow \dbR^{2d\times d_1}$ and $N:\triangle_*[0,T]\rightarrow \dbR^{2d\times d_2}$ with $d_1,d_2\in\dbN$,
\vskip-6mm
\bel{E241104}\ba{ll}
\ns\ds (M\1n\ltimes\1n\cP)(s,t)\1n\=\1nM(s,t)\cP^{(1)}(s)\1n+\2n\int_t^T\3nM(r,t)\cP^{(2)}(r,s,t)dr,(s,t)\1n\in\1n \triangle_{*}\1n[0,T],\\[-0.5mm]
\ns\ds(\cP\1n\rtimes\1n N)(s,t) \1n\=\1n\cP^{(1)}(s)N(s,t)\1n+\2n\int_t^T\3n\cP^{(2)}(s,r,t)N(r,t)dr,(s,t)\1n\in\1n \triangle_{*}\1n[0,T],\\[-1.5mm]
\ns\ds(M\1n\ltimes\1n\sP\1n\rtimes \1nN)(t)\1n\=\1n\int_t^T\3nM(s,t)\1n\Big[\cP^{(1)}(s)N(s,t)\1n+\2n\int_t^T\3n\cP^{(2)} (s,\th,t)N (\theta,t)d\th \Big]ds.
\ea
\ee\vskip-2mm
\no Thus,  the Lyapunov system (\ref{E241101}) can be rewritten as:
\vskip-6mm
\bel{E241102}\left\{\3n\ba{ll}
\ns\ds \cP^{(1)}(t)\1n=F_1[\Xi;\cP](t)+\cQ_1(t),t\in(0,T),\\
\ns\ds  \cP^{(2)}(s,t,t)\1n=\cP^{(2)}(t,s,t)^{\top}\2n=F_2[\Xi,\Gamma;\cP](s,t)+\cQ_2(s,t),(s,t)\1n\in\1n\triangle_{*}[0,T],\\
\ns\ds \dot{\cP^{(2)}}(s,t,r)\1n=-F_3[\Xi,\Gamma;\cP](s,t,r)\1n-\1n\cQ_3(s,t,r),\;(s,t,r)\in \Box_3[0,T],
\ea\right.
\ee
\vskip-1.5mm
\no where
\vskip-1.5mm
\vskip-6mm
$$\ba{ll}
\ns\ds \2n F_1[\Xi;\cP](t)\1n\=(\cC^{\top}\1n\ltimes \cP \rtimes  \cC)(t)\2n+\2n (\cD^{\top}\1n\ltimes\1n\cP\1n\rtimes \1n\cC)(t)\Xi(t)
\2n+\2n \Xi(t)^{\top}\1n(\cC^{\top}\ltimes\1n\cP\1n\rtimes \1n\cD)(t)\\
\ns\ds\qq\qq\q+\Xi(t)^{\top}\1n(\cD^{\top}\1n\ltimes\1n\cP\1n\rtimes \1n\cD)(t)\Xi(t)\1n=\1n\Big(\2n\begin{array}{ccccc}
 f_1[\Th;\sP](t)\1n &\1n 0\\[-1mm]
  0\1n &\1n0
\end{array}
\2n\Big),t\in(0,T),\\[-1mm]
\ns\ds \2nF_2[\Xi,\Gamma;\cP](s,t)\1n\=(\cP\1n\rtimes\1n \cA)(s,t)\2n+\2n(\cP\1n\rtimes\1n \cB)(s,t)\Xi(t)\2n
+\2n\Gamma(s,t)^{\top}\1n(\cD^{\top}\1n\ltimes\1n\cP\1n\rtimes \1n\cC)(t)\\
\ns\ds\qq\qq\qq+\Gamma(s,t)^{\top}\1n(\cD^{\top}\1n\ltimes\1n\cP\1n\rtimes \1n\cD)(t)\Xi(t)\1n=\1n\Big(\2n\begin{array}{ccccc}
f_6[\Th;\sP](s,t) \2n   & \2n0 \\[-0.3mm]
f_3[\Th;\sP](t) \2n  & \2n0
\end{array}
\2n\Big),(s,t)\1n\in\1n\triangle_{*}[0,T],\\[-1mm]
\ns\ds\2nF_3[\Gamma;\1n\cP](\1ns,\1nt,\1nr)\2n\=\2n\Gamma(s,\1nr)\1n^{\top}\2n(\cB\1n^{\top}\2n\ltimes\1n\cP)(t,\1nr)\2n+\2n(\cP\1n\rtimes\1n \cB)(s,\1nr)\Gamma(t,\1nr)
\2n+\2n\Gamma(s,\1nr)\1n^{\top}\2n(\cD^{\top}\2n\ltimes\1n\cP\1n\rtimes \1n\cD)(t)\Gamma(t,\1nr)\\[-0.5mm]
\ns\ds\qq\qq\qq=\1n\Big(\3n\begin{array}{ccccc}
f_5[\Th;\sP](s,t,r) \3n   & \3n f_4[\Th;\sP](s,r)\1n^{\top}\1nI_{(r,T)}(t)   \\[-0.3mm]
f_4[\Th;\sP](t,r)I_{(r,T)}(s) \3n  & \3n f_2[\Th;\sP](r)
\end{array}
\3n\Big),(s,\1nt,\1nr)\1n\in\1n \Box_3[0,T],
\ea$$
\vskip-2.5mm
\no and $\wt Q(\cd)\1n\=\1n\Big(\3n\begin{array}{ccccc}
Q(\cd) \1n   &  \1n  0   \\[-0.3mm]
0  \1n   &  \1n  0
\end{array}
\3n\Big)$,
\vskip-6mm
%
%
%
%
  %
%
%
%
%
%
  %
%
%
%
%
  %
%
 %
%
\vskip-1.5mm
$$\ba{ll}
\ns\ds \2n \cQ_1(t)\2n\=\2n\Xi(t)\1n^{\top}\1nR(t)\Xi(t)\2n+\2n \wt Q(t),\;
\cQ_2(s,t)\2n\=\2n\Gamma(s,t)\1n^{\top}\1nR(t)\Xi(t),\;
 \cQ_3(s,t,r)\2n\=\2n\Gamma(s,r)\1n^{\top}\1nR(r)\Gamma(t,r).
\ea$$
\vskip-1.5mm
\no  Clearly, we can obtain that
\vskip-6mm
$$\ba{ll}
\ns\ds \big(F_1,F_2,F_3\big),\big(\cQ_1,\cQ_2,\cQ_3\big)\2n\in\2n L^\infty(0,T;\1n\mathbb{S}^{2d})\2n\times\2n L^2(\triangle_\ast[0,T];\1n\mathbb{R}^{2d\1n\times\1n2d})
\2n\times\2n L^{2,2,1}_{sym}(0,T;\1n\mathbb{R}^{2d\1n\times\1n 2d}).
\ea$$
\vskip-1.5mm
\no In this case, \cite[Theorem 4.12]{Hamaguchi-Wang-I} guarantees that the existence and uniqueness of the solution $\cP\1n=\1n(\cP^{(1)},\1n\cP^{(2)})$ to (\ref{E241102}),
which implies  the  well-posedness of Lyapunov system (\ref{E241101}).
The final assertion can be obtained directly from a proof similar to the sufficiency of Theorem \ref{the3}.
The proof is complete.
\endpf

%


An argument similar to \cite[Lemma 6.2]{Hamaguchi-Wang-II}  leads to the following result.
\bl\label{lem5-2}
For any $\Theta_1,\1n\Theta_3\1n\in\1nL^\infty(0,T;\mathbb{R}^{l\times d}),\1n\Theta_2\1n\in\1nL^2(\triangle_{*}[0,T];\mathbb{R}^{l\times d})$, let $\sP\1n=\1n(\sP_1,\1n\sP_2,\1n\sP_3,\1n\sP_4)\1n\in\1n\Upsilon[0,T]$ be the solution of Lyapunov system (\ref{E241101}).
Assume that the functional  $u\mapsto\widetilde{\cJ}^0(0,0,0;u)$ is uniformly convex, i.e., there  exists a  $\lambda'\1n>\1n 0$ such that
\vskip-5mm
\bel{E241105}\ba{ll}
\ns\ds \wt \cJ^0(0,0,0;u)\ges \lambda' \dbE\int_0^T|u(s)|^2ds,\;\;\forall u\in\mathcal{U}[0,T].
\ea\ee
\vskip-1mm
\no Then,  there exist $\lambda>0$ and $\alpha\in\mathbb{R}$ such that $\BR(t)\1n\ges\1n \lambda I_l$ for a.e.$t\1n\in\1n[0,T],$
(recall $\BR$ in Definition \ref{def2.11}), and for any $(\t,\1n\varphi_1,\1n\varphi_2)\1n\in\1n\widetilde{\cI}$,
\vskip-6mm
\bel{E241107}\ba{ll}
\ns\ds \varphi_2^{\top}\1n\sP_2(\t)\varphi_2\1n+\2n\int_\t^T\2n\Big[\varphi_1(s)\1n^{\top}\1n\sP_1(s)\varphi_1(s)
\2n+\2n\int_\t^T\3n\varphi_1(s)\1n^{\top}\1n\sP_4(s,r,\t)\varphi_1(r)dr\\[-1mm]
\ns\ds \qq\qq+2\varphi_2\1n^{\top}\1n\sP_3(s,\t)\varphi_1(s)\Big]\1nds
\1n\ges\1n \alpha \Big(|\f_2|^2+\2n\int_\t^T\3n |\f_1(s)|^2ds\Big).
\ea\ee
\vskip-3mm

\el

The following is the main theorem of this section.

\begin{theorem}\label{Riccati}
 Let $(\rm H1)$-$(\rm H2)$ hold. Then, the functional  $u\mapsto\widetilde{\cJ}^0(0,0,0;u)$ is uniformly convex  if and only if
  the Riccati system (\ref{E4.1}) admits a strongly regular solution $\sP\1n=\1n(\sP_1,\1n\sP_2,\1n\sP_3,\1n\sP_4)\1n\in\1n \Upsilon[0,T]$.

\end{theorem}

\proof

\rm \textbf{The necessity}: Suppose the functional  $u\mapsto\widetilde{\cJ}^0(0,0,0;u)$ is uniformly convex. We now prove the existence of a strongly regular solution to
 Riccati system (\ref{E4.1}) by a iterative method. To begin with, for simplicity, we define some suitable operators to represent Lyapunov system (\ref{E241101}).
Denote $\cT\1n:\1n \cY_0[0,T]\1n\rightarrow\1n \Upsilon[0,T]$ by $\cT(Q[\Th])\1n\=\1n \big(\cT_1(Q[\Th]),\cT_2(Q[\Th]),\cT_3(Q[\Th]),\cT_4(Q[\Th])\big)$ with
\vskip-6mm
$$\ba{ll}
\ns\ds \cT_1(Q[\Th])(t)\1n=\1nQ_1(t),\;\cT_2(Q[\Th])(t)\1n=G+\2n\int_t^T\2nQ_2(s)ds,\;t\in[0,T],\\[-1.5mm]
\ns\ds\cT_3(Q[\Th])(t,r)\2n=\2nQ_3(t)\2n+\3n\int_r^t\3nQ_4(t,s)ds,(t,\1nr)\1n\in\1n\triangle_{*}[0,T],\\[-1mm]
\ns\ds \cT_4(Q[\Th])(s,\1nt,\1nr)\2n=\2nQ_6(s,\1nt)I_{[t,T)}(s)\2n+\2nQ_6(t,\1ns)\1n^{\top}\1nI_{[0,t)}(s) \2n
+\3n\int_r^{s\wedge t}\3nQ_5(s,\1nt,\1n\th)d\th,(s,\1nt,\1nr)\2n\in\2n\Box_3[0,T],
\ea
$$
\vskip-3mm
\no where  $Q[\Th]\1n\=\1n\big(G, Q_1,\1n Q_2,\1n Q_3,\1n Q_4,\1n Q_5,\1n Q_6\big)$ (recall (\ref{E241103})), and $\cY_0[0,T]\1n\=\1n\mathbb{S}^d\1n\times\1n\cY[0,T]$.
 Similarly, we define $f[\Th](\sP): \Upsilon[0,T]\rightarrow \cY_0[0,T]$ by
\vskip-6mm
$$\ba{ll}
\ns\ds f[\Th](\sP)\1n\=\1n\big(0,f_1[\Th;\sP],\1nf_2[\Th;\sP],\1nf_3[\Th;\sP],\1nf_4[\Th;\sP],\1nf_5[\Th;\sP],\1nf_6[\Th;\sP]\big).
\ea
$$
\vskip-1.5mm
\no Obviously, $\cT$ and $f[\Th]$ are bounded linear operators. Then, let $\Theta_1,\1n\Theta_3\2n\in\2nL^\infty(0,T;\mathbb{R}^{l\times d}),$ $\Theta_2\1n\in\1nL^2(\triangle_{*}[0,T];\mathbb{R}^{l\times d})$,
 for each $Q[\Th]\1n\in\1n\cY_0[0,T]$, $\sP\1n\in\1n \Upsilon[0,T]$  is the solution to the Lyapunov system (\ref{E241101}) if
 and only if  the following equation holds:
\vskip-5.5mm
\bel{E241109}\ba{ll}
\ns\ds \sP=\cT\big(f[\Th](\sP)+Q[\Th]\big).
\ea
\ee
\vskip-1.5mm

 By Lemma \ref{lem5-1}, with $(\Th^{(0)}_1,\Th^{(0)}_2,\Th^{(0)}_3)=(0,0,0)$, the following equation admits a unique solution
 $\sP^{(0)}\1n\in\1n \Upsilon[0,T]$:
%
%
%
%
%
%
%
 %
%
%
\vskip-4mm
$$\ba{ll}
\ns\ds \sP^{(0)}=\cT\big(f[\Th^{(0)}](\sP^{(0)})+Q[\Th^{(0)}]\big).
\ea
$$
\vskip-1.5mm
\no For $i=1,2,\cdots$, we set
\vskip-6mm
\bel{E241108}\ba{ll}
\ns\ds\Th^{(i)}_1(t)\1n\=\1n-\big(R(t)\1n+\1n(D^{\top}\2n\triangleleft \1n\sP^{(i-1)}\1n \triangleright \1n  D)(t)\big)\1n^{-1}
(D^{\top}\1n\triangleleft \sP^{(i-1)} \triangleright   C)(t),\; t\1n\in\1n[0,T],\\[-1mm]
\ns\ds\Th^{(i)}_2(r,t)\2n\=\2n-\big(\1nR(t)\2n+\2n(D^{\top}\2n\triangleleft \1n\sP^{(i-1)} \1n\triangleright \1n  D)(t)\big)\1n^{-1}
(B\1n^{\top}\2n\triangleleft \1n\sP^{(i-1)}_{1,3,4})(r,\1nt),(r,\1nt)\1n\in\1n\triangle_{*}[0,T],\\[-1mm]
\ns\ds\Th^{(i)}_3(t)\1n\=\1n-\big(R(t)\1n+\1n(D^{\top}\2n\triangleleft\1n \sP^{(i-1)} \1n\triangleright  \1n D)(t)\big)\1n^{-1}
(B^{\top}\1n\triangleleft \sP^{(i-1)}_{2,3})(t),\; t\1n\in\1n[0,T],
\ea
\ee
\vskip-2mm
\no and  let $\sP^{(i)}$ be the solution to the following Lyapunov system:
%
%
%
%
%
%
%
%
%
 %
%
%
\vskip-5.5mm
\bel{E241110}\ba{ll}
\ns\ds \sP^{(i)}=\cT\big(f[\Th^{(i)}](\sP^{(i)})+Q[\Th^{(i)}]\big).
\ea
\ee
\vskip-1.5mm
\no Lemma \ref{lem5-2} implies that
$\BR^{(i)}(\cd)\1n\=\1nR(\cd)\2n+\2n(D^{\top}\3n\triangleleft\1n \sP^{(i)} \1n\triangleright  \1n D)(\cd)\1n\ges\1n \lambda I_l$ for a.e. $t\1n\in\1n[0,T],$
and
\vskip-6mm
\bel{E241111}\ba{ll}
\ns\ds \varphi_2^{\top}\1n\sP^{(i)}_2(\t)\varphi_2\1n+\2n\int_\t^T\2n\Big[\varphi_1(s)\1n^{\top}\1n\sP^{(i)}_1(s)\varphi_1(s)
\2n+\2n\int_\t^T\3n\varphi_1(s)\1n^{\top}\1n\sP^{(i)}_4(s,r,\t)\varphi_1(r)dr\\[-1mm]
\ns\ds \qq+2\varphi_2\1n^{\top}\1n\sP^{(i)}_3(s,\t)\varphi_1(s)\Big]\1nds\ges \alpha \Big(|\f_2|^2+\2n\int_\t^T\3n |\f_1(s)|^2ds\Big),\;\; \forall (\t,\1n\varphi_1,\1n\varphi_2)\1n\in\1n\widetilde{\cI}.
\ea\ee
\vskip-3mm

In the following, we turn to proving that $\{\sP^{(i)}\}_{i\in \mathbb{N}}$ converges to the strongly regular solution of the Riccati system (\ref{E4.1}).
For each $i\in \mathbb{N},$ denote
\vskip-5mm
$$\ba{ll}
\ns\ds \bar{\sP}^{(i)}\1n\=\1n\big(\sP^{(i)}_1\2n-\2n\sP^{(i+1)}_1,\sP^{(i)}_2\2n-\2n\sP^{(i+1)}_2,
\sP^{(i)}_3\2n-\2n\sP^{(i+1)}_3,\sP^{(i)}_4\2n-\2n\sP^{(i+1)}_4\big),\\[-1mm]
\ns\ds \bar{\Th}^{(i)}\1n\=\1n\big(\Th^{(i)}_1\1n-\1n\Th^{(i+1)}_1,\Th^{(i)}_2\1n-\1n\Th^{(i+1)}_2,
\Th^{(i)}_3\1n-\1n\Th^{(i+1)}_3\big).
\ea
$$
\vskip-1.5mm
\no By some straightforward calculations, we can get $\bar{\sP}^{(i)}$ satisfies the following equation:
\vskip-5.5mm
$$\ba{ll}
\ns\ds \bar{\sP}^{(i)}=\cT\big(f[\bar{\Th}^{(i+1)}](\bar{\sP}^{(i)})+\bar{Q}[\bar{\Th}^{(i)}]\big),
\ea
$$
\vskip-1.5mm
\no where $\bar{Q}[\bar{\Th}^{(i)}]\2n\=\2n\big(0,\bar{Q}_1,\bar{Q}_2,\1n\bar{Q}_3,\1n\bar{Q}_4,\1n\bar{Q}_5,\1n\bar{Q}_6\big)$ and
\vskip-6mm
$$\ba{ll}
\ns\ds \2n \bar{Q}_1(t)\1n\=\1n\bar{\Th}^{(i)}_1\1n(t)\1n^{\top}\1n\BR^{(i)}(t)\bar{\Th}^{(i)}_1(t),\;\;\;\;\;
\bar{Q}_2(t)\1n\=\1n\bar{\Th}^{(i)}_3(t)\1n^{\top}\1n\BR^{(i)}(t)\bar{\Th}^{(i)}_3(t),\;t\in[0,T],\\
\ns\ds\2n\bar{Q}_4(s,t)\2n\=\2n\bar{\Th}^{(i)}_3\1n(t)\1n^{\top}\1n\BR^{(i)}(s)\bar{\Th}^{(i)}_2\1n(s,t),\;
\bar{Q}_6(s,t)\2n\=\2n\bar{\Th}^{(i)}_2\1n(s,t)\1n^{\top}\1n\BR^{(i)}(s)\bar{\Th}^{(i)}_1\1n(t),\;(s,t)\1n\in\1n\triangle_{*}[0,T],\\
\ns\ds \2n\bar{Q}_3(t)\2n\=\2n\bar{\Th}^{(i)}_3(t)^{\top}\1n\BR\1n^{(i)}(t)\bar{\Th}^{(i)}_1(t),t\2n\in\2n[0,\1nT],
\bar{Q}_5(s,\1nt,\1n\th)\2n\=\2n\bar{\Th}^{(i)}_2(s,\th)^{\top}\1n\BR^{(i)}(\th)\bar{\Th}^{(i)}_2(t,\th),(\1ns,\1nt,\1n\th)\2n\in\2n\Box_3[0,\1nT].
\ea$$
\vskip-1.5mm
\no Using Lemma \ref{lem5-1} and the definition of $\bar{Q}[\bar{\Th}^{(i)}]$, we can get that 
\vskip-6mm
\bel{E241112}\ba{ll}
\ns\ds \2n\varphi_2\1n^{\top}\2n\bar{\sP}\1n^{(i)}_2\1n(\t)\varphi_2\2n+\3n\int_\t^T\2n\Big[\varphi_1\1n(s)\1n^{\top}\3n\bar{\sP}\1n^{(i)}_1\1n(s)\varphi_1\1n(s)
\2n+\1n2\varphi_2\1n^{\top}\2n\bar{\sP}\1n^{(i)}_3\1n(s,\1n\t)\varphi_1\1n(s)\2n
+\3n\int_\t^T\3n\2n\varphi_1(s)\1n^{\top}\3n\bar{\sP}\1n^{(i)}_4\1n(\1ns,\1nr,\1n\t)\varphi_1(r)dr\1n\Big]\1nds\\[-1mm]
%
%
%
%
 \ns\ds=\dbE\2n\int_\t^T\2n\big\langle\BR^{(i)}(t)u^{\bar{\Theta}^{(i)}\1n,0}(t), u^{\bar{\Theta}^{(i)}\1n,0}(t)\big\rangle dt\ges 0,\;\;\;\forall (\t,\1n\varphi_1,\1n\varphi_2)\1n\in\1n\widetilde{\cI},
\ea\ee
\vskip-3mm
\no where $u^{\bar{\Theta}^{(i)}\1n,0}(t)\2n\=\2n \bar{\Th}^{(i)}_1(t)\bar{X}_1(t)\2n+\2n\int_t^T\2n\bar{\Th}^{(i)}_2(r,t)\bar{\sX}_1(r,t)dr
\2n+\2n \bar{\Th}^{(i)}_3(t)\bar{\sX}_2(t), t\2n\in\2n[0,T]$, and
$(\bar{X}_1,\1n\bar{\sX}_1,\1n\bar{\sX}_2)$ is the extended causal feedback solution corresponding to
$(\bar{\Theta}^{(i)}_1,\1n\bar{\Theta}^{(i)}_2,\1n\bar{\Theta}^{(i)}_3,0)$ (recall (\ref{E0601})).

For each $i\1n\in\1n \mathbb{N}$ and any $(\t,\1n\varphi_1,\1n\varphi_2)\1n\in\1n\widetilde{\cI}$,
we define $p^{(i)}_3(\t)\1n:\1n  L^2(\t,T;\mathbb{R}^{d})\1n \mapsto\1n  \mathbb{R}^{d}$
and $p^{(i)}(\t)\1n:\1n  L^2(\t,T;\mathbb{R}^{d})\1n \mapsto\1n  L^2(\t,T;\mathbb{R}^{d})$ by:
\vskip-6mm
$$\ba{ll}
\ns\ds p^{(i)}_3(\t)\f_1\2n=\2n \int_\t^T\3n\1n\sP^{(i)}_3\1n(s,\t)\varphi_1(s)ds,\;
\big(p^{(i)}(\t)\f_1\big)(t)\2n=\2n \sP^{(i)}_1\1n(t)\varphi_1(t)\2n+\2n \int_\t^T\3n\1n\sP^{(i)}_4\1n(t,\1ns,\1n\t)\varphi_1(s)ds.
\ea$$
\vskip-2mm
\no The adjoint operator $p^{(i)}_3(\t)^*$ of $p^{(i)}_3(\t)$ from $L^2(\t,T;\mathbb{R}^{d})$ into $\mathbb{R}^{d}$ takes the form
\vskip-4mm
$$\ba{ll}
\ns\ds p^{(i)}_3(\t)^*\f_2\1n=\2n \sP^{(i)}_3(s,\t)^{\top}\1n\varphi_2.
\ea$$
\vskip-1mm
\no If we take $\mathbb{P}_i^{\t}\1n\=\1n\bigg(\2n\begin{array}{ccccc}
\sP^{(i)}_2(\t) \2n&\2n  p^{(i)}_3(\t)\\
p^{(i)}_3(\t)^*  \2n&\2n  p^{(i)}(\t)
\end{array}\2n\bigg)$, then it is easy to see that $\mathbb{P}^{\t}_i$ is a self-adjoint bounded linear operator on $\mathbb{R}^{d} \1n\times\1n L^2(\t,T;\mathbb{R}^{d})$.
In addition, it follows that
\vskip-6mm
\bel{E241113}\ba{ll}
\ns\ds \Big\langle \dbP_i^{\t} \Big(\2n\begin{array}{ccccc}
\f_2 \\[-1.4mm]
\f_1
\end{array}\2n\Big), \1n \Big(\2n\begin{array}{ccccc}
\f_2 \\[-1.4mm]
\f_1
\end{array}\2n\Big) \Big\rangle_{\mathbb{R}^{d} \1n\times\1n L^2(\t,T;\mathbb{R}^{d})}
%
\1n=\2n\varphi_2^{\top}\1n\sP^{(i)}_2(\t)\varphi_2\1n+\2n\int_\t^T\2n\Big[\varphi_1(s)\1n^{\top}\1n\sP^{(i)}_1(s)\varphi_1(s) \\
\ns\ds \qq\q+2\varphi_2\1n^{\top}\1n\sP^{(i)}_3(s,\t)\varphi_1(s)
\1n+\2n\int_\t^T\3n\varphi_1(s)\1n^{\top}\1n\sP^{(i)}_4(s,r,\t)\varphi_1(r)dr\Big]ds.
\ea\ee
\vskip-2mm
\no Hence, combining (\ref{E241111})-(\ref{E241113}), we obtain that for any $\t\1n\in\1n [0,T), \{\dbP_i^{\t}\}_{i\in\mathbb{N}}$ is a bounded and
monotone sequence of   self-adjoint operators. Thus, $\{\dbP_i^{\t}\}_{i\in\mathbb{N}}$ is strongly convergent.
Further, the operator norm $\|\dbP_i^{\t}\|_{\cL}$ of $\dbP_i^{\t}$ is estimated as
 $\|\dbP_i^{\t}\|_{\cL}\les \|\dbP_0^{\t}\|_{\cL}$, which yields the uniform boundedness of  $\{\dbP_i^{\t}\}_{i\in\mathbb{N}}$.
%

 In particular, we note that $(D^{\top}\1n\triangleleft\1n \sP^{(i)} \1n\triangleright \1n C)\1n\in\1n L^\infty(0,T;\mathbb{R}^{l\times d})$ can be rewritten as:
\vskip-6mm
$$\ba{ll}
%
%
%
  %
%
%
%
  %
\ns\ds  (D^{\top}\1n\triangleleft\1n \sP^{(i)} \1n\triangleright \1n C)(t)\1n=\1n \Big(\Big\langle \Big(\2n\begin{array}{ccccc}
D_k(T,t)\\[-0.6mm]
D_k(\cd,t)
\end{array}\2n\Big),  \dbP_i^{t}\Big(\2n\begin{array}{ccccc}
C_j(T,t) \\[-0.6mm]
C_j(\cd,t)
\end{array}\2n\Big)\Big\rangle_{\mathbb{R}^{d} \1n\times\1n L^2}\Big)_{k,j},\;t\1n\in\1n [0,T],
\ea
$$
\vskip-2mm
\no where $D_{k}$, $C_{j}$, $k\1n\in\1n \{1,2,\cdots,\1n l\},j\1n\in\1n \{1,2,\cdots,\1n d\}$ are components of $D$ and $C$.  Similarly, we have

\vskip-6mm
$$\ba{ll}
\ns\ds  \bigg(\2n\begin{array}{ccccc}
(\sP^{(i)}_{2,3}\1n\triangleright\1n B)(t) \\
(\sP^{(i)}_{1,3,4}\1n\triangleright\1n B )(s,t)
\end{array}\2n\bigg)\2n=\2n \bigg(\1n\Big(\1n\dbP_i^{t}\Big(\1n\2n\begin{array}{ccccc}
B_{1}(T,t) \\[-0.6mm]
B_{1}(\cd,t)
\end{array}\2n\1n\Big)\Big)(s),\cdots,\Big(\dbP_i^{t}\Big(\2n\begin{array}{ccccc}
B_{l}(T,t) \\[-0.6mm]
B_{l}(\cd,t)
\end{array}\2n\Big)\Big)(s)   \1n\bigg),(s,t)\2n\in\2n \triangle_{*}[0,T].
\ea
$$
\vskip-2mm
\no  The above equalities and the dominated convergence theorem yield that the sequences
$\Th^{(i)}_{1}$, $\Th^{(i)}_{2}$ and $\Th^{(i)}_{3}$ defined in (\ref{E241108})  converge,
and then consequently the sequences $f[\Th^{(i)}](\sP^{(i)})$ and $Q[\Th^{(i)}]$ in (\ref{E241110}) also converge.
Combining this  with the Lyapunov system (\ref{E241110}) we get that  the limit of $\{\sP^{(i)}\}_{i\in\mathbb{N}}$ exists as $i\rightarrow\infty$, i.e.,
\vskip-4mm
$$\ba{ll}
\ns\ds \lim\limits_{i\rightarrow\infty} \1n\sP^{(i)}_1\1n(t)\2n=\2n\sP_1(t),
\lim\limits_{i\rightarrow\infty} \1n\sP^{(i)}_2\1n(t)\2n=\2n\sP_2(t), $a.e.$ t\1n\in\1n [0,T],\;\;
\sup\limits_{i\in\mathbb{N}}\big(\|\sP^{(i)}_1\|_{L^\infty}\2n+\2n\|\sP^{(i)}_2\|_{L^\infty}\1n\big)\2n<\2n\infty,\\
\ns\ds \lim\limits_{i\rightarrow\infty}\1n \sP^{(i)}_3\1n(t,s)\1n=\1n\sP_3(t,s)\;\; $in$ \;\;\sL^2(\triangle_\ast[0,T]),\;\;\;\;\;
\lim\limits_{i\rightarrow\infty} \1n\sP^{(i)}_4\1n=\1n\sP_4\;\; $in$ \;\;L^{2,2,1}_{sym}(\Box_3[0,T]).
\ea
$$
\vskip-1mm
\no Further, we have $\lim\limits_{i\rightarrow\infty}\1n \BR^{(i)}(t)\2n=\2nR(t)\2n+\2n(D^{\top}\2n\triangleleft \1n\sP \triangleright \1n  D)(t)
\1n\ges\1n \lambda I_l$ for a.e. $t\1n\in\1n[0,T],$ and
\vskip-6mm
\bel{E24121110}\ba{ll}
\ns\ds
\widehat{\Theta}_1(t)=-\big(R(t)+(D^{\top}\1n\triangleleft \sP \triangleright   D)(t)\big)^{-1}(D^{\top}\1n\triangleleft \sP \triangleright   C)(t),\;\; t\in[0,T],\\[-1mm]
\ns\ds\widehat{\Theta}_2(r,t)\1n
=\1n-\big(R(t)+(D^{\top}\1n\triangleleft \sP \triangleright   D)(t)\big)^{-1}(B^{\top}\2n\triangleleft \1n\sP_{1,3,4})(r,t),\; (r,t)\1n\in\1n\triangle_{*}[0,T],\\[-1mm]
\ns\ds
\widehat{\Theta}_3(t)=-\big(R(t)+(D^{\top}\1n\triangleleft \sP \triangleright   D)(t)\big)^{-1}(B^{\top}\1n\triangleleft \sP_{2,3})(t),\;\; t\in[0,T].
\ea
\ee
\vskip-2mm
\no  Therefore, we see that $\sP\1n=\1n(\sP_1,\sP_2,\sP_3,\sP_4)\1n\in\1n \Upsilon[0,T]$ satisfies the Lyapunov system
\vskip-5.5mm
\bel{E241118}\ba{ll}
\ns\ds \sP=\cT\big(f[\widehat{\Theta}](\sP)+Q[\widehat{\Theta}]\big).
\ea
\ee
\vskip-1.5mm
\no Finally, substituting (\ref{E24121110}) into the above system (\ref{E241118})  yields that  $\sP\1n=\1n(\sP_1,\sP_2,\sP_3,$
$\sP_4)\1n\in\1n \Upsilon[0,T]$ is the strongly regular solution of the Riccati system (\ref{E4.1}).

\medskip

\textbf{The sufficiency}: Let $\sP\1n\in\1n \Upsilon[0,T]$ be  the strongly regular solution of the Riccati
system (\ref{E4.1}) and $(\widehat{\Theta}_1,\widehat{\Theta}_2,\widehat{\Theta}_3)$ satisfy  (\ref{E24121110}). Then, it is easy to observe that
 $\sP$ satisfies the Lyapunov system (\ref{E241101}) with $(\Th_1,\Th_2,\Th_3)\1n=\1n(\widehat{\Theta}_1,\widehat{\Theta}_2,\widehat{\Theta}_3)$.
 By Lemma \ref{lem5-1}, for any $v\in\mathcal{U}[0,T]$, we can obtain that
\vskip-5.5mm
$$\ba{ll}
\ns\ds \wt\cJ^0(0,0,0;u^{\widehat{\Theta},v})\1n=\1n\dbE\1n\int_\t^T\3n v(t)\1n^{\top}\1n \big(R(t)
\1n+\1n(D^{\top}\2n\triangleleft\1n \sP \1n\triangleright\1n   D)(t)\big)v(t) dt\1n\ges\1n
\lambda\dbE\1n\int_\t^T\3n| v(t)|^2 dt.
\ea
$$
\vskip-1.5mm
\no Hence, it follows from  \cite[Lemma 6.2]{Hamaguchi-Wang-II} that the functional  $u\1n\mapsto\1n\widetilde{\cJ}^0(0,0,0;u)$ is uniformly convex.

\endpf

\section{Five special cases}

To demonstrate our unified treatment, we revisit several relevant literature as special cases and give the detailed comparisons, as well as some interesting phenomena.

\subsection{ The case of $G=0$}

In this subsection, we revisit the particular case (see e.g. \cite{Hamaguchi-Wang-I,Hamaguchi-Wang-II}) when
the terminal cost functional disappear (i.e., $G=0$).
We will show that both the causal feedback strategy and the Riccati system (\ref{E4.1}) reduce to
those in \cite{Hamaguchi-Wang-II} under proper conditions.

Suppose the functional  $u\1n\mapsto\1n\cJ^0(0,0,0;u)$ is uniformly convex.
By Theorem \ref{Riccati} and Corollary \ref{cor4.2}, it is easy to see that the Riccati system (\ref{E4.1}) admits a unique strongly regular solution
$(\sP_1,0,0,\sP_4)\1n\in\1n \Upsilon[0,T]$.
Then, for each $M:\1n\triangle_*[0,T]\1n\rightarrow \dbR^{d_1\times d}$ and $N:\triangle_*[0,T]\rightarrow \dbR^{d\times d_2}$ with $d_1,d_2\in\dbN$,
the notations in Definition \ref{def2.10} can be rewritten as
\vskip-6mm
$$\ba{ll}
\ns\ds (M\triangleleft \sP_{2,3})(t)=0,\; (\sP_{2,3}\triangleright N)(t)=0,\;
(\1nM\1n\triangleleft \1n\sP\1n\triangleright\1n N)(t)\2n=\2n(M\1n\ltimes\1n\sP\1n\rtimes \1nN)(t) ,\;t\in(0,T),\\[-0.3mm]
\ns\ds(\1nM\1n\triangleleft \1n\sP_{1,3,4}\1n)(s,t) \2n=\2n(M\1n\ltimes\1n\sP)(s,t),
\;(\1n\sP_{1,3,4}\1n\triangleright \1nN)\2n=\2n(\sP\1n\rtimes\1n N)(s,t),\;(s,t)\1n\in\1n \triangle_{*}\1n[0,T].
\ea
$$\vskip-2mm
\no Here, we recall the notations $\ltimes$ and $\rtimes$  in  (\ref{E241104}).

As to Riccati system (\ref{E4.1}), it reduces to (recall $\dot{\sP_4}$ in (\ref{Notation-use-0})):
\vskip-6mm
$$\left\{\3n\ba{ll}
\ns\ds \sP_1(t)\1n=\1nQ(t)\1n+\1n(C^{\top}\2n\ltimes\1n \sP\1n\rtimes\1n C)(t)\1n-\1n\big(S(t)^{\top}\2n+\1n(C^{\top}\2n\ltimes \1n\sP\1n\rtimes\1n D)(t)\big)\\
\ns\ds\qq\qq\times\big(R(t)\1n+\1n(D^{\top}\2n\ltimes \1n\sP\1n\rtimes \1nD)(t)\big)^{-1}\1n
\big(S(t)\1n+\1n(D^{\top}\2n\ltimes \1n\sP\1n\rtimes\1n C)(t)\big),t\in(0,T),\\
\ns\ds  \sP_4(s,t,t)\1n=\1n(\sP\1n\rtimes\1n A)(s,t)\1n-\1n(\sP\1n\rtimes \1n B)(s,t)
\big(R(t)\1n+\1n(D^{\top}\1n\ltimes \1n\sP\1n\rtimes \1nD)(t)\big)^{-1}\1n(D^{\top}\1n\ltimes \1n\sP\1n\rtimes\1n C)(t)\\
\ns\ds\qq\qq=\sP_4(t,s,t)^{\top}\1n,\;\;(s,t)\in\triangle_\ast[0,T],\\
\ns\ds \dot{\sP}_4(s_1,s_2,t)\2n=\2n(\sP\1n\rtimes  \1nB)(s_1,t)\big(R(t)\1n+\1n(D^{\top}\1n\ltimes
\1n\sP\1n\rtimes\1n D)(t)\big)^{-1}\1n(B^{\top}\1n\ltimes\1n\sP)(s_2,t),(s_1,s_2,t)\1n\in\1n \Box_3(0,T),
\ea\right.
$$
\vskip-1.5mm
\no which is exactly the \emph{Riccati-Volterra equation} in \cite{Hamaguchi-Wang-II}.
 In this case, it follows that
\vskip-5mm
$$\ba{ll}
\ns\ds  \widehat{\Theta}_1(t)\1n=\1n-\big(R(t)\1n+\1n(D^{\top}\2n\ltimes\1n \sP\1n\rtimes\1n D)(t)\big)^{-1}
(D\1n^{\top}\2n\ltimes \1n\sP\1n\rtimes\1n C)(t),\;\;\widehat{\Theta}_3(t)\1n=\1n0,\; \;\;t\1n\in\1n(0,T),\\
\ns\ds\widehat{\Theta}_2(r,t)\1n=\1n-\big(R(t)\1n+\1n(D^{\top}\2n\ltimes \1n\sP\1n\rtimes\1n D)(t)\big)\1n^{-1}\1n(B^{\top}\2n\ltimes \1n\sP)(r,t),\;\; (r,t)\1n\in \1n\triangle_\ast[0,T],\\
%
%
 %
\ea
$$
\vskip-2mm

\no By Theorem \ref{the3}, we can write the solution of (\ref{YZ3}) as
$(\dbY,\dbZ)\1n=\1n\bigg(\bigg(\3n\begin{array}{ccccc}
Y_1     \\[-1.2mm]
Y_2    \\[-1.4mm]
0
\end{array}\3n\bigg),\bigg(\3n\begin{array}{ccccc}
Z_1     \\[-1.2mm]
Z_2    \\[-1.4mm]
0
\end{array}\3n\bigg)\bigg).$ Define
\vskip-6mm
$$\ba{ll}
\ns\ds  \eta (t,s)\1n\=\1n Y_1(t,s)\1n+\1n\h \Th_1(t)^{\top} Y_2(t,s)\1n+\2n\int^T_t \2n\h \Th_2(t,r)^{\top} Y_2(r,s)dr,\\
\ns\ds  \zeta (t,s)\1n\=\1n Z_1(t,s)\1n+\1n\h \Th_1(t)^{\top} Z_2(t,s)\1n+\2n\int^T_t \2n\h \Th_2(t,r)^{\top} Z_2(r,s)dr.
\ea
$$
\vskip-2mm
\no Then (\ref{YZ3}) can be rewritten as
\vskip-6mm
$$\ba{ll}
\ns\ds  \bigg(\begin{array}{ccccc}
\2n Y_1(t,s) \\[-1.2mm]
\2n Y_2(t,s)\\[-1.4mm]
0
\end{array}
\2n\bigg)\1n=\1n \bigg(\begin{array}{ccccc}
\2n q(t)\1n+\1n(C^{\top}\2n\ltimes\1n \sP\1n\rtimes\1n \si)(t) \\[-1.2mm]
\2n \rho(t)\1n+\1n(D^{\top}\2n\ltimes\1n \sP\1n\rtimes\1n \si)(t)\\[-1.4mm]
0
\end{array}\2n\bigg)\1n+\2n\int^T_t \2n\bigg(\begin{array}{ccccc}
\2n A(r,t)^{\top}\1n\eta(r,t) \\[-1.2mm]
\2n B(r,t)^{\top}\1n\eta(r,t)\\[-1.4mm]
0
\end{array}\2n\bigg)dr\1n+\2n\int^T_t \2n\bigg(\begin{array}{ccccc}
\2n C(r,t)^{\top}\1n\zeta(r,t) \\[-1.2mm]
\2n D(r,t)^{\top}\1n\zeta(r,t)\\[-1.4mm]
0
\end{array}\2n\bigg)dr\\
\ns\ds \qq\qq\qq-\2n\int^t_s \2n\bigg(\begin{array}{ccccc}
\2n ( \sP\1n\rtimes\1n b)(t,r)\\[-1.2mm]
\2n 0\\[-1.4mm]
0
\end{array}\2n\bigg)dr-\2n\int^t_s \2n\bigg(\begin{array}{ccccc}
\2n Z_1(t,r)\\[-1.2mm]
\2n Z_2(t,r)\\[-1.2mm]
0
\end{array}\2n\bigg)dW(r).
\ea
$$
\vskip-2mm
\no Therefore, it follows that $(\eta, \zeta)$ satisfies the following equation
\vskip-6mm
$$\ba{ll}
\ns\ds  \eta (t,s)\2n=\2n q(t)\2n+\2n\h \Th_1(t)^{\top}\1n\big[\rho(t)\2n+\2n(D^{\top}\2n\ltimes\1n \sP\1n\rtimes\1n \si)(t)\big]
\2n+\2n(C^{\top}\2n\ltimes\1n \sP\1n\rtimes\1n \si)(t)\2n+\3n\int^T_t \3n\1n\Big\{ \big[A(r,t)^{\top}\1n\eta(r,t)\2n+\2n C(r,t)^{\top}\1n\zeta(r,t)\big] \1n\\
 \ns\ds\qq\qq + \h \Th_1(t)^{\top}\1n\big[B(r,t)^{\top}\1n\eta(r,t)\1n+\1nD(r,t)^{\top}\1n\zeta(r,t)\big] \Big\}dr
 \1n+\1n\int^T_s \2n\h \Th_2(t,r)^{\top}\Big\{ \big[\rho(r)\1n+\1n(D^{\top}\2n\ltimes\1n \sP\1n\rtimes\1n \si)(r)\big]\\
 \ns\ds\qq\qq +\2n\int^T_r\2n\big[B(\th,r)^{\top}\1n\eta(\th,r)\1n+\1nD(\th,r)^{\top}\1n\zeta(\th,r)\big] d\th \Big\}dr\1n-\1n\int^t_s \2n ( \sP\1n\rtimes\1n b)(t,r) dr\1n-\1n\int^t_s \2n \zeta(t,r)dW(r),
\ea
$$
\vskip-2mm
\no which is also the \emph{Type-II EBSVIE} in \cite{Hamaguchi-Wang-II}.
Then the process $\h\sN$ defined by (\ref{sN}) can be rewritten as
\vskip-8mm
$$\ba{ll}
\ns\ds\h \sN(t)\1n=\1n\rho(t)\1n+\1n (D^{\top}\2n\ltimes\sP\rtimes \si)(t)\1n+\2n\int_t^T\2nB(s,t)^{\top}\eta (s,t)ds\1n
+\2n\int_t^T\2nD(s,t)^{\top}\zeta(s,t)d s,\;\;t\in(0,T).
\ea
$$
\vskip-3mm

In this case,  the optimal causal feedback strategy $(\widehat{\Theta}_1,\widehat{\Theta}_2,\widehat{\Theta}_3,\widehat{v})$ admits
\vskip-6mm
$$\ba{ll}
\ns\ds  \widehat{\Theta}_1(t)\1n=\1n-\big(R(t)\1n+\1n(D^{\top}\2n\ltimes\1n \sP\1n\rtimes\1n D)(t)\big)^{-1}
(D\1n^{\top}\2n\ltimes \1n\sP\1n\rtimes\1n C)(t),\;\;\widehat{\Theta}_3(t)\1n=\1n0,\; \;\;t\1n\in\1n(0,T),\\
\ns\ds\widehat{\Theta}_2(r,t)\1n=\1n-\big(R(t)\1n+\1n(D^{\top}\2n\ltimes \1n\sP\1n\rtimes\1n D)(t)\big)\1n^{-1}\1n(B^{\top}\2n\ltimes \1n\sP)(r,t),\;\; (r,t)\1n\in \1n\triangle_\ast[0,T],\\
\ns\ds \widehat{v}(t)\1n=\1n-\big(R(t)\1n+\1n(D^{\top}\2n\ltimes \1n\sP\1n\rtimes\1n D)(t)\big)^{-1}\h \sN(t)\;\;t\1n\in\1n(0,T).
\ea
$$
\vskip-2.5mm
\no Then, it is easy to see that the above results are consistent with \cite{Hamaguchi-Wang-I, Hamaguchi-Wang-II}.

\subsection{The case of SDEs}

Let us consider one particular case with $f(s,t)\equiv f(t)$, $f(\cd,\cd)\= A(\cd,\cd), B(\cd,\cd), C(\cd,\cd), D(\cd,\cd)$, $\f(t)\equiv \f$, $\f\in\dbR^d$. The state equation (\ref{E1}) becomes a SDE and Problem (LQ-SVIE) reduces to classical
stochastic LQ problem. Next we show that our obtained result fully covers those in the SDEs case.

To begin with, let us define $\cP(t)$ as follows:
\vskip-5.5mm
$$\ba{ll}
\ns\ds \cP(t)\1n\=\1n\sP_2(t)\2n+\3n\int_t^T\3n\Big(\1n\sP_1(s)\2n+\2n
\sP_3(s,t)\2n+\2n\sP_3(s,t)^{\top}\3n+\3n\int_t^T\3n\2n\sP_4(s,\th,t)d\th\1n\Big)ds,\;t\in[0,T].
\ea
$$
\vskip-2.5mm
 \no Clearly, inserting the  Riccati system (\ref{E4.1}) into $\dot{\cP}(t)$, it follows  that
\vskip-6.5mm
$$\left\{\3n\ba{ll}
\ns\ds \dot{\cP}(t)=-A(t)^{\top}\cP(t)-\cP(t)A(t)
-C(t)^{\top}\cP(t)C(t)-Q(t)\1n+\1n\big[\cP(t)B(t)\\[-1mm]
\ns\ds\q\1n+C(t)\1n^{\top}\1n\cP(t)D(t)\big]\big[R(t)\2n+\2nD(t)\1n^{\top}\1n\cP(t)D(t)\big]^\dagger\1n
\big[B(t)\1n^{\top}\1n\cP(t)\2n+\2nD(t)\1n^{\top}\1n\cP(t)C(t)\big],t\1n\in\1n[0,T],\\[-1mm]
\ns\ds \cP(T)=G.
\ea\right.
$$
\vskip-2.5mm
\no Obviously, it is just the traditional Riccati equation in stochastic LQ problem.
Further, it follows from Theorem \ref{the3} that
\vskip-6.5mm
$$\ba{ll}
\ns\ds \h\Th(t)\1n\=\1n\h\Theta_1(t)\1n+\2n\int^T_t\3n \h\Theta_2(s,t)ds\1n+\1n\h\Theta_3(t)\1n
=\1n-\big[R(t)\1n+\1nD(t)\1n^{\top}\cP(t)D(t)\big]^\dagger
\big[B(t)^{\top}\cP(t)\1n+\1nD(t)^{\top}\cP(t)C(t)\big],
\ea
$$
\vskip-2.5mm
\no which is  consistent with  the result in the SDEs case.


\subsection{The case of stochastic VIDE}

In this subsection, we will investigate the connections between the matrix-valued  optimal causal feedback strategy in this paper and the optimal feedback operator in \cite{Wang-Yong-Zhou-2023-SICON}. Hence, we consider $b\1n=\1n\si\1n=\1nq\1n=\1n\rho\1n=\1nS\1n=\1n0$ here.

To begin with, we give  the following standard condition:

\ss

\textbf{(H3)}. $Q(t)\1n\ges \1n 0$, $R(t)\1n\ges \1n\lambda I_l$ for $t\1n\in\1n [0,T]$ and some $\lambda\1n>\1n 0$, $G\1n\ges \1n0$.

\ss

From \cite[Corollary 6.3]{Hamaguchi-Wang-II},  we can obtain the following result.
\begin{corollary}\label{cor2701}
\vskip-1mm
Suppose  $(\rm H3)$ hold.  Then  the functional  $u\1n\mapsto\1n\widetilde{\cJ}^0(0,0,0;u)$ is uniformly convex.
\end{corollary}

Combining Theorem \ref{the3}, Theorem \ref{Riccati} and Corollary \ref{cor2701},  it is easy to show the following result.
\begin{corollary}\label{cor5.2}
\vskip-1.5mm
Suppose $(\rm H1)$--$(\rm H3)$ hold. Then Riccati system (\ref{E4.1}) admits  a unique strongly regular solution $\sP\1n=\1n(\sP_1,\sP_2,\sP_3,\sP_4)\1n\in\1n\Upsilon[0,T]$.
 Consequently,  there exists a unique optimal 4-tuple $(\widehat{\Theta}_1,\widehat{\Theta}_2,\widehat{\Theta}_3,0)$ for  Problem (LQ-SVIE) that satisfies
\vskip-6mm
\bel{E241117}\ba{ll}
\ns\ds
\widehat{\Theta}_1(t)=-\big(R(t)+(D^{\top}\1n\triangleleft \sP \triangleright   D)(t)\big)^{-1}(D^{\top}\1n\triangleleft \sP \triangleright   C)(t),\;\; t\in[0,T],\\[-1mm]
\ns\ds\widehat{\Theta}_2(r,t)\1n
=\1n-\big(R(t)+(D^{\top}\1n\triangleleft \sP \triangleright   D)(t)\big)^{-1}(B^{\top}\2n\triangleleft \1n\sP_{1,3,4})(r,t),\; (r,t)\1n\in\1n\triangle_{*}[0,T],\\[-1mm]
\ns\ds
\widehat{\Theta}_3(t)=-\big(R(t)+(D^{\top}\1n\triangleleft \sP \triangleright   D)(t)\big)^{-1}(B^{\top}\1n\triangleleft \sP_{2,3})(t),\;\; t\in[0,T].
\ea
\ee
\vskip-2mm
\no Further, for each  $(\t,\varphi)\in\cI$, the value function is given by
$ V^0(\t, \f)\1n=\1n(\varphi^{\top}\1n\triangleleft\sP\triangleright\f)(\t).$
\end{corollary}

By  Corollary \ref{cor5.2} and the fact that
\vskip-7mm
$$\ba{ll}
\ns\ds (B^{\top}\1n\triangleleft \sP_{2,3})(t)\widehat{\sX}(T,t)\1n
+\2n\int^T_t\2n(B^{\top}\2n\triangleleft \1n\sP_{1,3,4})(s,t)\widehat{\sX}(s,\1nt)ds=(B^{\top}\1n\triangleleft\1n \sP\triangleright\1n \widehat{\sX})(t),t\1n\in\1n [\t,T],
\ea
$$
\vskip-2.5mm
\no  the following defined control process $\hat{u}$ is optimal
\vskip-5mm
\bel{e003}\ba{ll}
\ns\ds
\hat{u}(t)\=\1n-\1n\big(\1nR(t)\2n+\2n(D^{\top}\3n\triangleleft\1n \sP\1n\triangleright\1nD)(t)\big)^{-1}
\1n\big[(\1nB^{\top}\3n\triangleleft \1n\sP\1n\triangleright\1n\widehat{\sX})\1n(t)
\2n+\2n(D^{\top}\3n\triangleleft\1n\sP\1n\triangleright\1n C)(t)\widehat{\sX}(t,t)\big],
\ea
\ee
\vskip-1mm
\no where  $(\widehat{X},\widehat{\sX})\1n\in\1n  L^2_{\mathbb{F}}(\t,T;\mathbb{R}^{ d})\times L^2_{\mathbb{F},c}(\triangle_\ast[\t,T];\mathbb{R}^{ d})$
is causal feedback solution of (\ref{E1})
at $(\t,\varphi)$ corresponding to the optimal $(\widehat{\Theta}_1,\widehat{\Theta}_2,\widehat{\Theta}_3,0)$.

\ss

Next, in order to compare with the relevant results in \cite{Wang-Yong-Zhou-2023-SICON},  we take $\Lambda_\t \1n\=\1n C([\t, T];\dbR^d)$ for given $\t\in[0,T)$,
and consider the state equation in \cite{Wang-Yong-Zhou-2023-SICON}, i.e., (\ref{E1}) with the free term $\f(\cd,\t)\1n\in\1n\Lambda_\t$ and the coefficients satisfy the following stronger assumption:

\textbf{(H4)}. $A,C:\triangle_\ast[\t,T]\rightarrow\mathbb{R}^{d\times d}, B,D:\triangle_\ast[\t,T]\rightarrow\mathbb{R}^{d\times l}$
  are bounded and  partially differentiable with respect the two variables, with bounded derivatives.

Then, to define the optimal feedback  operator similar to \cite{Wang-Yong-Zhou-2023-SICON}, we first introduce the bilinear functional $\dbP(\t):\Lambda_\t\times \Lambda_\t\rightarrow\dbR$:
\vskip-5mm
\bel{052302}\ba{ll}
\ns\ds
\dbP(\t)\big(N(\cd,\t),M(\cd,\t)\big)\=(M^{\top}\1n\triangleleft \sP \triangleright   N)(\t),\;\;\forall M(\cd,\t), N(\cd,\t)\in \Lambda_\t,
\ea
\ee
\vskip-1mm
\no where $\sP$ is the strongly regular solution of the Riccati system (\ref{E4.1}).
 Let $\sS(\Lambda_\t)$ be the set of all bounded $\dbR$-valued bilinear functionals on $\Lambda_\t\times \Lambda_\t$ such that
\vskip-5mm
$$\ba{ll}
\ns\ds P(\t)\big(M(\cd,\t), N(\cd,\t)\big)\1n=\1nP(\t)\big(N(\cd,\t),M(\cd,\t)\big),\;\;\forall M(\cd,\t), N(\cd,\t)\in \Lambda_\t,
\ea
$$
\vskip-1mm
\no under the norm
\vskip-5mm
$$\ba{ll}
\ns\ds \|P(\t)\|_{\sS}\=\sup_{\|M\|\les 1}\big|P(\t)\big(M(\cd,\t),M(\cd,\t)\big)\big|,\;\;\forall  P(\t)\1n\in\1n \sS(\Lambda_\t).
\ea
$$
\vskip-1.5mm
\no It follows from the Riccati system (\ref{E4.1}) that $\dbP(\t)\in \sS(\Lambda_\t)$. Further, for any $t\in [\t,T]$,
$D(\cd,t)\1n=\1n\big(D_1(\cd,t),\cdots,D_l(\cd,t)\big)\1n\in\1n C([t,T];\dbR^{d\times l})$,
we define a bilinear functional $\dbP_1(t):C([t,T];\dbR^{d\times l})\1n\times\1n C([t,T];\dbR^{d\times l})\1n\rightarrow\1n\dbS^{ l}$ by:
\vskip-5mm
$$\ba{ll}
\ns\ds \dbP_1(t)\big(D(\cd,t),D(\cd,t)\big)\=\Big(\dbP(t)\big(D_i(\cd,t),D_j(\cd,t)\big)\Big),\;\;i,j\1n\in\1n\{1,2,\cdots,l\}.
\ea
$$
\vskip-1.5mm
\no Similarly, for $B(\cd,t)\1n=\1n\big(B_1(\cd,t),\cdots,B_l(\cd,t)\big)\1n\in\1n C([t,T];\dbR^{d\times l})$ and any $M(\cd,t)\1n\in\1n\Lambda_t$,
 we introduce a linear functional $\dbP_2(t):\Lambda_t\rightarrow\dbR^{ l}$:
\vskip-5mm
$$\ba{ll}
\ns\ds \dbP_2(t)M(\cd,t)\1n\=\1n\Big(\dbP(t)\big(B_i(\cd,t),M(\cd,t)\big)\2n+\2n\dbP(t)\big(D_i(\cd,t),C(\cd,t)M(t,t)\big)\Big)
\ea
$$
\vskip-1.5mm
\no for $i\1n\in\1n\{1,2,\cdots,l\}$. Hence, the value function can be rewritten as:
\vskip-5mm
\bel{051501}\ba{ll}
\ns\ds V^0(\t, \f)=\dbP(\t)\big(\f(\cd),\f(\cd)\big),\;\;\; \forall (\t,\varphi)\1n\in\1n\cI,
\ea
\ee
\vskip-1mm
\no  and the optimal control defined by (\ref{e003}) can be denoted as:
\vskip-3.5mm
$$\ba{ll}
\ns\ds\hat{u}(t)\=-\big[R(t)+\dbP_1(t)\big(D(\cd,t),D(\cd,t)\big)\big]^{-1}\dbP_2(t)\widehat{\sX}(\cd,t),\;t\in [\t,T],
\ea
$$
\vskip-1.5mm
\no which inspires us to define the optimal feedback  operator as follows:
\vskip-5mm
\bel{051502}\ba{ll}
\ns\ds
\h\Gamma(t)\=-\big[R(t)+\dbP_1(t)\big(D(\cd,t),D(\cd,t)\big)\big]^{-1}\dbP_2(t).
\ea
\ee
\vskip-1.5mm

Therefore, utilizing (\ref{051501}), the arbitrariness of $ \f$ and the result of \cite[Lemma 2.8]{Wang-Yong-Zhou-2023-SICON}, it follows that the bilinear functional $\dbP(\t)$ defined by (\ref{052302}) is  consistent with
 the strongly regular solution of so-called \emph{path-dependent Riccati equation} in \cite{Wang-Yong-Zhou-2023-SICON}.
Consequently, we can obtain  that the optimal feedback operator (\ref{051502}) coincides with that in \cite{Wang-Yong-Zhou-2023-SICON}.

\subsection{ The case of VIEs}

In this subsection, we look at one particular case with $A\1n=\1nC\1n=\1nD\1n=\1nb\1n=\1n\si\1n\equiv\1n0$ in (\ref{E1}). Hence (\ref{E1}) reduces to the deterministic VIEs, i.e.,
\vskip-5mm
\bel{E5.9}
X(t)=\varphi(t)+\int_\t^tB(t,s)u(s)ds,\;t\in[\t,T],
\ee
\vskip-2.5mm
\no where $\varphi$ and $B$ satisfy $(\rm H1)$. The corresponding $\sX$ becomes
\vskip-3.5mm
$$\ba{ll}
\ns\ds
\sX(r,t)=\varphi(r)+\int_\t^tB(r,s)u(s)ds,\;\;(r,t)\in\triangle_\ast[\t,T].
\ea
$$
\vskip-1.5mm
\no  For any $u\in\mathcal{U}[\t,T]$,
we see that $\sX(r,t)$ is equivalent to the so-called \emph{$t$-causal trajectory} $X_t(r)$ for $r> t$ in \cite{Pritchard-You-1996},
where
\vskip-7mm
\bel{E2603}\ba{ll}
\ns\ds X_t(r)\1n\=\1n\varphi(r)+\int_\t^rB(r,s)(\pi_tu)(s)ds, \;\; r\in[\t,T],
\ea
\ee
\vskip-3mm
\no and $\pi_t$ is a truncation operator defined as $(\pi_t u)(s)\=u(s)I_{[0,t]}(s)$.

 In this case, the Riccati system (\ref{E4.1}) can be simplified as
\vskip-6.5mm
\bel{E5.12}\left\{\3n\ba{ll}
\ns\ds \1n\sP_1(t)\1n=\1nQ(t),\;\sP_2(t)\1n=\1nG\1n-\2n\int_t^T\2n(\sP_{2,3}\triangleright B)(s)R(s)^{-1}(B^{\top}\1n\triangleleft \sP_{2,3})(s)ds,\;t\in[0,T],\\[-1mm]
\ns\ds\1n\sP_3(t,r)\1n=\1n-\1n\int_r^t\2n(\sP_{2,3}\1n\triangleright\1n B)(\th)R(\th)^{-1}\1n
(B^{\top}\2n\triangleleft\1n \sP_{1,3,4})(t,\th)d\th,\;(t,r)\1n\in\1n\triangle_{*}[0,T],\\[-1mm]
\ns\ds\1n\sP_4(s,t,r)\1n=\1n-\2n\int_r^{s\wedge t}\3n(\sP_{1,3,4}\1n\triangleright\1n B)(s,\th)R(\th)^{-1}
(B^{\top}\3n\triangleleft \1n\sP_{1,3,4})(t,\th)d\th,(s,t,r)\1n\in\1n\Box_3[0,\1nT].
\ea\right.
\ee
\vskip-3mm
\no
Assume that $(\rm H1)$--$(\rm H3)$ hold.
Then, Corollary \ref{cor4.2}, Theorem \ref{Riccati}  and  Corollary \ref{cor2701} yield the above Riccati system (\ref{E5.12})
admits a unique strongly regular solution $\sP\1n=\1n(\sP_1,\1n\sP_2,\1n\sP_3,\1n\sP_4)$.
Further,  Theorem \ref{the3} implies that
 the LQ problem admits the unique optimal 4-tuple $(0,\1n\widehat{\Theta}_2,\1n\widehat{\Theta}_3,\1n0)$ where
\vskip-3.5mm
$$\ba{ll}
\ns\ds
\widehat{\Theta}_3(t)=-R(t)^{-1}(B^{\top}\1n\triangleleft \sP_{2,3})(t),\;\; t\in[0,T],\\
\ns\ds\widehat{\Theta}_2(r,t)\1n
=\1n-R(t)^{-1}(B^{\top}\2n\triangleleft \1n\sP_{1,3,4})(r,t)I_{[0,r]}(t),\; (r,t)\1n\in\1n [0,T]^2.
\ea
$$
\vskip-1mm
\no  In addition, for any $(\t,\varphi)\in\cI$, the following defined $\hat{u}(\cd)$ is optimal
\vskip-5.5mm
\bel{E2601}\ba{ll}
\ns\ds
\hat{u}(t)\=\1n\int^T_t \2n\h\Th_2(r,t)\hat{\sX}(r,t)dr+\h\Th_3(t)\hat{\sX}(T,t),\;t\in[\t,T],
\ea
\ee
\vskip-3mm
\no where $(\hat{X},\hat{\sX})$ is causal feedback solution of (\ref{E5.9})  corresponding to
 $(0,\widehat{\Theta}_2,\widehat{\Theta}_3,0)$.

\ss

 We recall the relevant result in \cite{Pritchard-You-1996}.
Under $(\rm H1)$-$(\rm H4)$,  $\hat{u}(\cd)$ is the optimal control
\vskip-5mm
\bel{E2602}\ba{ll}
\ns\ds \hat{u}(t)=\1n\int^T_t \2n\Xi_1(r,t)X_t(r)dr+\Xi_2(t)X_t(T),\;t\in[0,T],
\ea
\ee
\vskip-2mm
\no where $X_t(\cd)$ is the corresponding $t$-causal trajectory defined by (\ref{E2603}), and $(\Xi_1, \Xi_2)$ is  determined by
the solution of a Fredholm integral equation. In addition, utilizing $(\rm H1)$-$(\rm H4)$  and the result of \cite[Theorem 5.3]{Pritchard-You-1996},
it follows that the 4-tuple $(0,\Xi_1, \Xi_2,0)$ is a causal feedback strategy.

Notice that Theorem \ref{the3}  and the strongly regular solvability of the Riccati system (\ref{E5.12}) guarantee the uniqueness of
the optimal feedback strategy $(0,\1n\widehat{\Theta}_2,\1n\widehat{\Theta}_3,\1n0)$.
Therefore, combining (\ref{E2601}), (\ref{E2602}) and $\hat{\sX}(r,t)\1n=\1nX_t(r), r\1n>\1nt$, it follows that
$(\widehat{\Theta}_2,\widehat{\Theta}_3)\1n=\1n(\Xi_1,\Xi_2)$, i.e.,
the optimal causal feedback strategy in this paper is consistent with
 that in \cite{Pritchard-You-1996}.
 Based on this finding, we answered (Q4) posed in the introduction.

\subsection{ The case of VIDEs}

In this subsection, let us revisit the particular VIDEs case in \cite{Pandolfi-2018-IEEETAC} and make careful comparisons.
Consider the following controlled VIDE:
\vskip-5mm
\bel{E12101}
\dot{X}(t)=\int_0^tN(t-s)X(s)ds +Bu(t),\;\;\;X(0)=\f,
\ee
\vskip-2mm
\no where $\dot{X}$ denotes the derivative of $X$. It corresponds to the case of
\vskip-3mm
$$
A(t,s)\1n\equiv\2n\int_s^t\1nN(\th\1n-\1ns)d\th,\;\; B(s,t)\1n\equiv \1nB,\;\; C (s,t)\1n=\1n D(s,t)\1n\equiv\1n0,\;\; \f(t)\1n\equiv\1n \f,
$$
\vskip-1mm
\no with $B,\f\in\dbR^d$. Moreover, let $Q(\cd)\1n\equiv \1nQ, R(\cd)\1n\equiv\1n1$ with any $Q\in\dbS^d$ in (\ref{E1.2}).
Then Problem (LQ-SVIE) reduces to the VIDEs case in \cite{Pandolfi-2018-IEEETAC}.

It is clear that the corresponding $\sX$ can be rewritten as:
\vskip-5.5mm
\bel{E2605}\ba{ll}
\ns\ds \sX(s,t)=X(t)+\int_0^t\int_t^sN(\theta-r)d\theta X(r)dr,\;\;\;(s,t)\in\triangle_{*}[0,T].
\ea
\ee
\vskip-2.5mm
\no Also the Riccati system (\ref{E4.1}) becomes
\vskip-6.5mm
\bel{E24040901}\left\{\3n\ba{ll}
\ns\ds \sP_1(t)\1n=\1nQ,\;\;\;\;\sP_2(t)\1n=\1nG\1n-\2n\int_t^T\2n\sP_{2,3}(s)^{\top} BB^{\top}\1n \sP_{2,3}(s)ds,\;t\in[0,T],\\[-1mm]
\ns\ds\sP_3(t,r)\1n=\1n(\sP_{2,3}\1n\triangleright\1n A)(t)
\2n-\3n\int_r^t\3n\1n\sP_{2,3}(\th)^{\top}BB^{\top}\1n \sP_{1,3,4}(t,\th)d\th,\;(t,r)\1n\in\1n\triangle_{*}[0,T],\\[-1mm]
\ns\ds\sP_4(s,t,t)\1n=\1n\sP_4(t,s,t)^{\top}=(\sP_{1,3,4}\1n\triangleright\1n A)(s,t),\;\;(s,t)\1n\in\1n\triangle_{*}[0,T],\\
\ns\ds\dot{\sP}_4(s,t,r)=\sP_{1,3,4}(s,r)^{\top}BB^{\top}\sP_{1,3,4}(t,r),\;\;(s,t,r)\in\Box_3[0,T],
\ea\right.
\ee
\vskip-2mm
\no where
\vskip-8mm
$$\ba{ll}
\ns\ds\sP_{2,3}(t)\=\sP_2(t)+\2n\int_t^T\sP_3(s,t)^{\top}ds,\;\;t\in[0,T],\\[-2mm]
\ns\ds\sP_{1,3,4}(s,t)\=\sP_1(s)+\sP_3(s,t)+\1n\int_t^T\1n\sP_4(r,s,t)dr,\;\;\;(s,t)\in \triangle_{*}[0,T].
\ea
$$
\vskip-3mm
\no  Assume that  $(\rm H3)$ hold.
Then, Corollary \ref{cor4.2}, Theorem \ref{Riccati}  and  Corollary \ref{cor2701} yield the above Riccati system (\ref{E24040901})
admits a unique strongly regular solution $(\1n\sP_1,\1n\sP_2,\1n\sP_3$, $\1n\sP_4)$.
Therefore, by Theorem \ref{the3}, the LQ problem admits an optimal 4-tuple $(0,\widehat{\Theta}_2,\widehat{\Theta}_3,0)$ and optimal control
\vskip-5.5mm
\bel{E24040902}\ba{ll}
\ns\ds\hat{u}(t)\1n=-B^{\top}\Big[
\1n\int^T_t \3n\sP_{1,3,4}(s,t)\sX(s,t)ds
\1n+\1n \sP_{2,3}(t)\sX(T,t)\Big],\;t\in[\t,T],
\ea
\ee
\vskip-3mm
\no where $(\1nX,\1n\sX\1n)$ is causal feedback solution of (\ref{E12101}) at $(\t,\varphi)$ corresponding to
$(0,\widehat{\Theta}_2,\widehat{\Theta}_3,0)$.
Plugging (\ref{E2605}) into (\ref{E24040902}), we can obtain  relationship w.r.t. control
\vskip-5.5mm
\bel{Feedback-compare}\ba{ll}
\ns\ds \hat{u}(t)=-B^{\top}\Big[p_0(t)X(t)+\int_0^tp_1(t,s)X(s)ds\Big],\;t\in[\t,T],
\ea
\ee
\vskip-3mm
\no where
\vskip-7.5mm
$$\ba{ll}
\ns\ds p_0(t)\=\sP_{2,3}(t)+\int^T_t \sP_{1,3,4}(s,t)ds,\\[-1.8mm]
\ns\ds p_1(t,s)\=\int^T_t\Big[\sP_{2,3}(t)\1n+\2n\int^T_\th \3n\1n\sP_{1,3,4}(r,t)dr\Big]N(\theta-s)d\theta.
\ea
$$
\vskip-1.5mm
Let us show the coincidence between (\ref{E24040901})
  and the \emph{Riccati differential equation} in \cite{Pandolfi-2018-IEEETAC}.
First, by the definition of $p_0(\cd)$,
\vskip-6mm
\bel{E24040905}\ba{ll}
\ns\ds \dot{p}_0(t)\2n=\2n\Big[\sP_{2,3}(t)^{\top}\3n+\3n\int_t^T\2n \2n\sP_{1,3,4}(r,t)^{\top}\2ndr\Big]
\1n BB^{\top}\1n \Big[\sP_{2,3}(t)\2n+\3n\int_t^T\3n\1n\sP_{1,3,4}(s,t)ds\Big]\1n-\1nQ\\[-1mm]
\ns\ds\qq\q\1n-(A^{\top}\2n\triangleleft\1n\sP_{2,3})(t)\1n
-\2n\int_t^T\3n\Big((A^{\top}\2n \triangleleft\1n\sP_{1,3,4})(s,t)\1n+\1n(\sP_{1,3,4}\1n\triangleright\1n A)(s,t)\Big)ds
\1n-\1n(\sP_{2,3}\1n\triangleright \1nA)(t).
\ea
\ee
\vskip-3mm
\no Using Fubini theorem and the fact of $A(t,s)\1n\equiv\2n\int_s^t\1nN(\th\1n-\1ns)d\th,$ we obtain that
\vskip-3mm
$$\ba{ll}
\ns\ds p_1(t,t)\1n=\1n(\sP_{2,3}\1n\triangleright \1nA)(t)\1n+\2n\int_t^T\3n(\sP_{1,3,4}\1n\triangleright\1n A)(s,t)ds.
\ea
$$
\vskip-2mm
\no Plugging it into (\ref{E24040905}),  we find that
\vskip-6mm
\bel{E24040907}\ba{ll}
\ns\ds \dot{p}_0(t)=p_0(t)BB^{\top}p_0(t)-Q-p_1(t,t)-p_1(t,t)^{\top}, \ \ p_0(T)=G.
\ea
\ee
\vskip-1.5mm
 Next, let us calculate the derivative of $p_1(\cd,\cd)$ w.r.t the first variable,
\vskip-6mm
$$\ba{ll}
\ns\ds \frac{\partial}{\partial t}p_1(t,s)\1n=\1np_0(t)N(t\1n-\1ns)\1n+\1np_0(t)BB^{\top}\1np_1(t,s)\2n
-\3n\int^T_t\3n\Big[(A^{\top}\2n\triangleleft\1n\sP_{2,3})(t)\\[-2mm]
\ns\ds\qq\qq\q+\2n\int^T_\th\2n(A^{\top}\2n \triangleleft\1n\sP_{1,3,4})(r,t)dr\Big]N(\theta\1n-\1ns)d\theta,\ \ (t,s)\in \triangle_{*}[0,T].
\ea
$$
\vskip-2mm
\no To calculate the last term on the right-hand, we define
\vskip-6mm
$$\ba{ll}
\ns\ds K(t,s,\rho)\1n\=\2n\int^T_t\3n\bigg\{\2n\int^T_t\3nN(\th\1n-\1n\rho)\1n^{\top}\1n d\th \sP_2(t)\1n
+\2n\int^T_t\3n\int^r_t\3nN(\th\1n-\1n\rho)\1n^{\top}\1n d\th\sP_3(rt)\1n^{\top}\1n dr
\1n+\2n\int^T_{\th_1}\3n\Big[\1n\int^r_t\3nN(\th\1n-\1n\rho)\1n^{\top}\1n d\th \sP_1(t)\\
\ns\ds\qq\qq+\2n\int^T_t\3nN(\th\1n-\1n\rho)^{\top}\1n d\th \sP_3(r,t)\1n
+\2n\int^T_{t}\3n\int^{r_1}_t\2nN(\th\1n-\1n\rho)^{\top}\1n d\th  \sP_4(r_1,r,t)  dr_1\Big]dr\bigg\}N(\th_1\1n-\1ns)d\th_1
\ea
$$
\vskip-1.5mm
\no for any $0\les s,\rho \les t\les T$. It is easy to see that
\vskip-3mm
$$\ba{ll}
\ns\ds K(t,s,t)=\2n\int^T_t\2n\Big[(A^{\top}\2n\triangleleft\1n\sP_{2,3})(t)
\1n+\2n\int^T_\th\3n(A^{\top}\2n \triangleleft\1n\sP_{1,3,4})(r,t)dr\1n\Big]N(\theta\1n-\1ns)d\theta,
\ea
$$
\vskip-1.5mm
\no which implies that
\vskip-5mm
\bel{E24040908}\ba{ll}
\ns\ds \frac{\partial}{\partial t}p_1(t,s)\1n=\1np_0(t)N(t\1n-\1ns)\1n+\1np_0(t)BB^{\top}\1np_1(t,s)-K(t,s,t),\;\;p_1(T,s)=0.
\ea
\ee
\vskip-2mm
At last, it is time for us to treat the new term $K(\cd,\cd,\cd)$. By the definition of $p_1$ and Fubini theorem, a direct computation shows that
\vskip-7mm
$$\ba{ll}
\ns\ds \frac{\partial}{\partial t}K(t,s,\rho)\1n=\1n-\1nI_1\2n-\2nI_2\1n+\3n\int^T_t\2n\bigg\{\1n\int^T_t\3n\1nN(\th\2n-\2n\rho)^{\top}\1n d\th \dot{\sP}_2(t)\2n
+\3n\int^T_t\3n\1n\int^r_t\3nN(\th\2n-\2n\rho)^{\top}\1n d\th\dot{\sP}_3(r,t)^{\top}\1n dr\\[-1mm]
\ns\ds\qq\2n+\3n\int^T_{\th_1}\2n\Big[\1n\int^T_t\3n\1nN(\th\2n-\2n\rho)^{\top}\2n d\th \dot{\sP}_3(r,t)\2n
+\3n\int^T_{t}\3n\1n\int^{r_1}_t\3nN(\th\2n-\2n\rho)^{\top}\2n d\th \dot{\sP}_4(r_1,r,t)  dr_1\1n\Big]dr\1n\bigg\}N(\th_1\2n-\2ns)d\th_1,
\ea
$$
\vskip-1.5mm
\no where
\vskip-6mm
\bel{E24040909}\ba{ll}
 I_1=p_1(t,\rho)^{\top}N(t-s),\;\;I_2=N(t-\rho)^{\top}p_1(t,s).
\ea
\ee
\vskip-3mm
\no By (\ref{E24040909}), Fubini theorem, the definition of $p_1$ and the Riccati system (\ref{E24040901}),
\vskip-6mm
\bel{E24040910}\ba{ll}
\ns\ds \frac{\partial}{\partial t}K(t,s,\rho)\2n=\2n-\1np_1(t,\rho)\1n^{\top}N(t\2n-\2ns)\2n-\2nN(t\2n-\2n\rho)\1n^{\top}
p_1(t,s)\2n+\2np_1(t,\rho)\1n^{\top}\1nBB\1n^{\top}p_1(t,s),K(T,s,\rho)\2n=\2n0.
\ea
\ee
\vskip-1.5mm
Combining  (\ref{E24040907}), (\ref{E24040908}), (\ref{E24040910}) and (\ref{Feedback-compare}), we can see that our Riccati system (\ref{E4.1}) reduces to the so-called \emph{Riccati differential equation} in \cite{Pandolfi-2018-IEEETAC} and our outcome process $\h u(\cd)$ also coincides with their feedback control. To sum up the above arguments, we point out several differences/advantages of the current study. First, one key notion here is optimal causal feedback strategy while \cite{Pandolfi-2018-IEEETAC} is concerned with the closed-loop/feedback representation of optimal controls. These two notions may be essentially different in the games framework.
Second, our stochastic framework is a nontrivial extension of their deterministic one. Third, we provide the necessary and sufficient conditions of causal feedback strategy, not just the sufficiency as in \cite{Pandolfi-2018-IEEETAC}.

\section{Concluding remarks}

In this paper we introduce a new kind of optimal causal feedback strategy and characterize it in terms of a novel Riccati system
with the solvability being carefully discussed.
From our study, various forms of causal feedback controls in the existing papers can be unified properly.
In addition, under SVIEs framework, many issues concerning controllability, stability and games (especially mean-field games) remain open and challenging.
 We will present some relevant results in the future.

\end{document}